\documentclass[12pt]{amsart}
\usepackage{amsfonts,latexsym,amsthm,amssymb,graphicx}
\usepackage[all]{xy}

\CompileMatrices

\newtheorem{theorem}{Theorem}[section]
\newtheorem{lemma}[theorem]{Lemma}
\newtheorem{corol}[theorem]{Corollary}
\newtheorem{prop}[theorem]{Proposition}

{\theoremstyle{definition} \newtheorem{defin}[theorem]{Definition}}
{\theoremstyle{remark} 
\newtheorem{example}[theorem]{Example}}

\newcommand{\Abb}{{\mathbb{A}}}
\newcommand{\Cbb}{{\mathbb{C}}}
\newcommand{\Qbb}{{\mathbb{Q}}}

\newcommand{\Pbb}{{\mathbb{P}}}
\newcommand{\Zbb}{{\mathbb{Z}}}

\newcommand{\Til}[1]{{\widetilde{#1}}}

\DeclareMathOperator{\rk}{rk}

\DeclareMathOperator{\im}{im}

\DeclareMathOperator{\Spec}{Spec}

\DeclareMathOperator{\Hom}{Hom}

\newcommand{\PGL}{\text{\rm PGL}}
\newcommand{\Sym}{\text{\rm Sym}}
\newcommand{\nG}{\overline{\Pbb}}
\newcommand{\oD}{\overline{D}}
\newcommand{\oE}{\overline{E}}
\newcommand{\oalpha}{\overline{\alpha}}
\newcommand{\obeta}{\overline{\beta}}
\newcommand{\nor}{n}
\newcommand{\onu}{\overline{n}}

\title{Limits of translates of plane curves ---\\
on a paper of Aldo Ghizzetti}
\author{Paolo Aluffi, Carel Faber}
\address{Max-Planck-Institut f\"ur Mathematik, Postfach 7280,
D-53072 Bonn, Germany}
\address{Dept.~of Mathematics, Florida State University, Tallahassee
FL 32306, U.S.A.}
\email{aluffi@math.fsu.edu}
\address{Inst.~f\"or Matematik, Kungliga Tekniska H\"ogskolan, 
S-100 44 Stockholm, Sweden}
\email{faber@math.kth.se}

\makeindex

\begin{document}

\begin{abstract}
We study the limits of $\PGL(3)$ translates of an arbitrary plane
curve, giving a description of all possible limits of a given curve
and computing the multiplicities of corresponding components in the
normal cone to the base scheme of a related linear system. This
information is a key step in the computation of the degree of the
closure of the {\em linear orbit\/} of an arbitrary plane curve.

Our analysis recovers and extends results obtained by
Aldo Ghizzetti in the 1930's.
\end{abstract}

\maketitle


\section{Introduction}\label{intro}

Let $\mathcal C$ be an arbitrary complex plane curve of degree $d$.
Consider $\mathcal C$ together with all its translates: the orbit
of $\mathcal C$ for the natural action of $\PGL(3)$ on the
projective space~$\Pbb^n$ of all plane curves of degree $d$. Which
plane curves appear in the orbit closure of~$\mathcal C$? Or in other
words, what are the limits of translates of $\mathcal C$?
In this article we answer a refined form of this question.

Harris and Morrison (\cite{MR99g:14031}, p.138) define the flat
completion problem for embedded families of curves as the determination of
all curves in $\Pbb^n$ that can arise as flat limits of a family of 
embedded stable curves over the punctured disc. The problem mentioned
in the first paragraph contains the isotrivial case of the flat
completion problem for plane curves, and a solution to it can in fact
be found in the marvelous article~\cite{ghizz} by the Italian
mathematician Aldo Ghizzetti (a summary of the results
is contained in~\cite{ghizzLin}). However, as we will explain below,
our main application requires a more refined type of information; thus
our aim is somewhat different than Ghizzetti's, and we cannot simply
lift his results. Consequently, our work in this paper is independent
of \cite{ghizz}. In any case, Ghizzetti's approach has substantially
influenced ours; see \S\ref{Ghizz} for a description of his work and a
comparison with ours. 

The enumerative geometry of families of plane curves with prescribed
singularities presents a notoriously difficult problem. Spectacular
progress was made in the last decade in several special cases; we mention
the work of Kontsevich \cite{MR97d:14077}
and of Caporaso-Harris \cite{MR99i:14064}. Consider the special case
where the family consists of a completely arbitrary plane curve and all
its translates. In our paper~\cite{MR2001h:14068} we explained how the
degree of this family, in other words, the number of curves in the family
passing through the appropriate number of general points in the plane,
can be computed.
For example, for a nonsingular curve $\mathcal C$ this family is the
set of all possible embeddings of $\mathcal C$ in $\Pbb^2$; our
motivation in \cite{MR94e:14032}, \cite{MR2002d:14083},
\cite{MR2002d:14084}, \cite{MR2001h:14068}, and the present article is
the study of this set, and of its natural generalization for arbitrary
plane curves.

The starting point of our method is to view the action map on $\mathcal C$
as a rational map $c$ from the $\Pbb^8$ of $3\times3$ matrices to the
$\Pbb^n$ of plane curves of degree $d$. We require a precise description
of the closure of the graph of $c$, specifically of the scheme-theoretic
inverse image of the locus of indeterminacy. This is the exceptional
divisor of the blow-up of $\Pbb^8$ along the base scheme of $c$, the
projective normal cone (PNC). A set-theoretic description of the
components of the PNC amounts to a solution of the (isotrivial)
flat completion problem, together with careful bookkeeping of the
different arcs in $\Pbb^8$ used to obtain each limit. The PNC can be
viewed as an arc space associated to the rational map $c$; it is
probably possible to recast our analysis in \S\ref{setth} in the
light of recent work on arc spaces (cf.~for example \cite{MR1905328}),
and it would be interesting to do so.

In fact a set-theoretic description of the PNC does not suffice for
our enumerative application in \cite{MR2001h:14068}. This requires
the full knowledge of the PNC {\em as a cycle,\/} that is, the determination
of the multiplicity of its different components.
Thus, we determine not only the limits of 
one-dimensional families of translates of $\mathcal C$, but we also
classify such families up to a natural notion of equivalence and we keep
track of the behavior of a family near the limit. The determination
of the limits and the classification are contained in \S\ref{setth}. As may be
expected, the determination of the multiplicities is quite delicate;
this is worked out in \S\ref{multi}. Preliminaries, and a more
detailed introduction, can be found in \S\ref{prelim}.

The final result of our analysis is stated in \S2 of
\cite{MR2001h:14068}, in the form of five `Facts'. The proofs of these
facts are spread over the present text; we recommend comparing
loc.cit.~with \S\ref{statement} and \S\ref{multstat} to establish a
connection between the more detailed statements proved here and the
summary in \cite{MR2001h:14068}.

Caporaso and Sernesi use our determination of the
limits in \cite{MR2003k:14035} (Theorem 5.2.1). Hassett
\cite{MR2000j:14045} and Hacking \cite{math.AG/0310354} study the
limits of the family of nonsingular plane curves of a given degree, by
methods different from ours: they allow the plane to degenerate
together with the curve. It would be interesting to compare their
results to ours.

\vskip6pt

\noindent{\bf Acknowledgments.} We thank an anonymous referee of our
first article on the topic of linear orbits of plane curves,
\cite{MR94e:14032}, for bringing the paper of Aldo Ghizzetti to our
attention. A substantial part of the work leading to the results in
this paper was performed while we enjoyed the peaceful atmosphere at
Oberwolfach during a `Research in Pairs' stay, and thanks are due to
the Volkswagen Stiftung for supporting the R.i.P.~program.

The first author thanks the Max-Planck-Institut f\"ur Mathematik in
Bonn, Germany, for the hospitality and support, and Florida State
University for granting a sabbatical leave in 2001-2.  The second
author thanks Princeton University for hospitality and support during
the spring of 2003.  The first author's visit to Stockholm in May 2002
was made possible by support from the G\"oran Gustafsson
foundation. The second author's visit to Bonn in July 2002 was made
possible by support from the Max-Planck-Institut f\"ur Mathematik.

\newpage


\section{Preliminaries}\label{prelim}

\subsection{}\label{rough}
We work over $\Cbb$. Roughly speaking, the question we (and Ghizzetti)
address is the following: given a plane curve $\mathcal C$, what plane
curves can be obtained as limits of translates of $\mathcal C$? By
a translate we mean the action on $\mathcal C$ of an invertible linear
transformation of $\Pbb^2$, that is, an element of $\PGL(3)$. We view
$\PGL(3)$ as an open set in the space $\Pbb^8$ parametrizing $3\times
3$ matrices $\alpha$; if $F(x,y,z)$ is a generator of the homogeneous
ideal of $\mathcal C$, $\alpha$ acts on $\mathcal C$ by composition:
we denote by $\mathcal C\circ\alpha$ the curve with ideal
$(F\circ\alpha)=(F(\alpha(x,y,z)))$, and we call the set of all
translates $\mathcal C\circ\alpha$ the {\em linear orbit\/} of
$\mathcal C$. Our guiding question here concerns the limits of families
$\mathcal C\circ g(a)$, for $g:\mathcal A \to \Pbb^8$ any map
from a smooth curve $\mathcal A$ to $\Pbb^8$, centered at a point
mapping to a singular transformation. 

Since the flat limit is determined by the completion of the local ring
of $\mathcal A$ at the center, we may replace $\mathcal A$ with
$\Spec\Cbb[[t]]$. Thus, a `curve germ in $\Pbb^8$' ({\em germ\/} for
short) in this article will simply be a $\Cbb[[t]]$-valued point
$\alpha(t)$ of $\Pbb^8$. 
Our `germs' are often called `arcs' in the literature.

The limit of $\mathcal C\circ\alpha(t)$ as $t \to 0$:
$$\lim_{t\to 0}\,\mathcal C\circ \alpha(t)$$ 
is the flat limit over the punctured
$t$-disk; concretely, this is obtained by clearing common powers of
$t$ in the expanded expression $F(\alpha(t))$ and then setting $t=0$.

It will always be assumed that the center $\alpha=\alpha(0)$ of a germ
$\alpha(t)$ is a singular transformation. Further, we may and will
assume that $\alpha(t)$ is invertible for some $t\ne 0$: indeed, this
condition may be achieved by perturbing every $\alpha(t)$ achieving a
limit, using terms of high enough power in $t$ so as not to affect the
limit. 

We note that, by the same token, every limit attained by one of our
`germs' can conversely be realized as the flat limit of a family
parametrized by a curve $\mathcal A$ mapping to $\Pbb^8$, as above;
and in fact we can even assume $\mathcal A=\Abb^1$. Indeed, truncating
$\alpha(t)$ at a high enough power does not affect $\lim_{t\to
  0}\mathcal C\circ \alpha(t)$, and a polynomial $\alpha(t)$ describes
a map $\Abb^1 \to \Pbb^8$.

\subsection{}\label{7ic}
Here is one example showing that rather interesting limits may occur:
let $\mathcal C$ be the 7-ic curve with equation
$$x^3 z^4-2 x^2 y^3 z^2+x y^6-4 x y^5 z-y^7=0$$
and the family
$$\alpha(t)=\begin{pmatrix}
1 & 0 & 0\\
t^8 & t^9 & 0\\
t^{12} & \frac 32 t^{13} & t^{14}
\end{pmatrix}\quad;$$
then
$$F\circ\alpha(t)=-\frac 1{16}\,t^{52}(x^3(8 x^2+3 y^2-8 x z)(8 x^2-3
y^2+8 x z))+\text{ higher order terms.}$$
That is,
$$\lim_{t\to 0}F\circ\alpha(t)=-\frac 1{16}\,x^3(8 x^2+3 y^2-8 x z)(8
x^2-3 y^2+8 x z)\quad,$$
a pair of quadritangent conics (see \cite{MR2002d:14083}, \S4.1),
union the distinguished tangent line taken with multiplicity~3. Note
that the connected component of the identity in the stabilizer of this
curve is the additive group.

\subsection{}
Our primary objective is essentially to describe the possible limits
of a plane curve $\mathcal C$, starting from a description of certain
features of $\mathcal C$. We now state this goal more precisely.

In the process of computing the degree of the closure of the linear
orbit of an arbitrary curve, \cite{MR2001h:14068}, we are led to
studying the closure of the graph of the rational map
$$\Pbb^8 \dashrightarrow \Pbb^N$$
mapping an invertible $\alpha$ to $\mathcal C\circ\alpha$, viewed as a
point of the space $\Pbb^N$ parametrizing plane curves of degree $\deg
\mathcal C$. This graph may be identified with the blow-up of $\Pbb^8$
along the base scheme $S$ of this rational map. Much of the enumerative
information we seek is then encoded in the {\em exceptional divisor\/}
$E$ of this blow-up, that is, in the projective normal cone of $S$ in
$\Pbb^8$. In fact, the information can be obtained from a description
of the components of $E$ (viewed as 7-dimensional subsets of
$\Pbb^8\times \Pbb^N$) and from the multiplicities with which these
components appear in $E$. A more thorough discussion of the relation
between this information and the enumerative results, as well as of
the general context underlying our study of linear orbits of plane
curves, may be found in \S1 of \cite{MR2001h:14068}.

Our goal in this paper is the description of the components of the
projective normal cone of $S$, and the computation of the
multiplicities with which the components appear in the projective
normal cone.

\subsection{}
This goal relates to the one stated more informally in \S\ref{rough}
in the sense that the components of the projective normal cone
dominate subsets of the boundary of the linear orbit of $\mathcal C$. 
Our technique will consist of studying an arbitrary $\alpha(t)$,
aiming to determine whether the limit $(\alpha(0),\lim_{t\to 0}
\mathcal C\circ \alpha(t))$ is a general point of the support of a
component of the normal cone; we will thus obtain a description of all
components of the normal cone. We should warn the reader that we will
often abuse the language and refer to the point $(\alpha(0),\lim_{t\to 0}
\mathcal C\circ \alpha(t))$ (in $\Pbb^8\times \Pbb^N$) by the
typographically more convenient $\lim_{t\to 0} \mathcal C\circ
\alpha(t)$ (which is a point of $\Pbb^N$).

We should also point out that our analysis will not exhaust the
boundary of a linear orbit: one component of this boundary may arise as 
the closure of the set of translates $\mathcal C\circ\alpha$ with $\alpha$ 
a rank-2 transformation, and a general such $\alpha$ does not belong to $S$.
Indeed, the rational map mentioned above {\em is\/} defined
at the general $\alpha$ of rank~2. To be more precise, if $\alpha$ is
a rank-2 matrix whose image is not contained in $\mathcal C$ (for
example, if $\mathcal C$ has no linear components), then $\mathcal C
\circ\alpha$ may be described as a `star' of lines through
$\ker\alpha$, reproducing projectively the tuple of points cut out by
$\mathcal C$ on the image of $\alpha$.
It would be interesting to provide a
precise description of the set of stars arising in this manner.
As this set does not contribute components to $E$, this study is not
within the scope of this paper.

\subsection{}\label{rank2obs}
One of our main tools in the set-theoretic determination of
the components of $E$ will rely precisely on the fact that limits
along rank-2 transformations do not contribute components to $E$. We
will argue that if a limit obtained by a germ $\alpha(t)$ can also be
obtained as a limit by a germ $\beta(t)$ contained in the rank-2
locus, then we can `discard' $\alpha(t)$. Indeed, such limits will
have to lie in the exceptional divisor of the blow-up of the rank-2
locus; as the rank-2 locus has dimension $7$, such limits will span
loci of dimension at most $6$. We will call such limits `rank-2
limits' for short. The form in which this observation will be applied
is given in Lemma~\ref{rank2lemma}.

Incidentally, germs $\alpha(t)$ centered at a rank-2 transformation
whose image is contained in $\mathcal C$ {\em do\/} contribute a
component to $E$ (cf.~\S\ref{globallin}); by the argument given in the
previous paragraph, however, contributing $\alpha(t)$ will necessarily
be invertible for $t\ne 0$.

\subsection{}
In \S\ref{setth} we will determine the components of $E$
set-theoretically, as subsets of $\Pbb^8\times\Pbb^N$. As a
preliminary observation (cf.~\cite{MR2001h:14068}, p.~8) we can
describe the whole of $E$ set-theoretically in terms of limits, as
follows. Recall that $S$ denotes the base locus of the rational map
$c:\Pbb^8 \dashrightarrow \Pbb^N$ defined by $\alpha\mapsto\mathcal C
\circ\alpha$.

\begin{lemma}\label{excdivsupp}
As a subset of $\Pbb^8\times\Pbb^N$, the support of $E$ is
\begin{multline*}
|E|=\{(\alpha,X)\in \Pbb^8\times\Pbb^N : \text{$X$ is a limit of
  $c(\alpha(t))$}\\
\text{for some curve germ $\alpha(t)\subset \Pbb^8$ centered at
  $\alpha\in S$ and not contained in $S$}\}\quad.
\end{multline*}
\end{lemma}

\begin{proof}
Let $\Til \Pbb^8$ be the closure of the graph of the rational map $c$
defined above. This is an 8-dimensional irreducible variety, mapping
to $\Pbb^8$ by the restriction $\pi$ of the projection on the first
factor of $\Pbb^8\times\Pbb^N$, and identified with the blow-up of
$\Pbb^8$ along the base scheme $S$ of $c$. The set $E$ is the inverse
image $\pi^{-1}(S)$ in $\Til \Pbb^8$.

Any curve germ $\alpha(t)$ in $\Pbb^8$ centered at $\alpha\in S$ and
not contained in $S$ lifts to a germ in $\Pbb^8$ centered at a point
of $E$; this yields the $\supset$ inclusion.

For the other inclusion, let $\tilde\alpha(t)$ be a germ centered at a
point $\tilde\alpha$ of $E$, and such that $\tilde\alpha(t_0)\not\in
E$ for $t_0$ near $0$; such a germ may be obtained (for example) by
successively intersecting $\Til\Pbb^8$ with general divisors of type $(1,1)$
through $\tilde\alpha$. As $\pi$ is 1-to-1 in the complement of $E$,
$\tilde\alpha(t)$ is the lift of a (unique) curve germ $\alpha(t)$ in
$\Pbb^8$, giving the other inclusion.
\end{proof}

\subsection{}
As mentioned in \S\ref{intro}, our application in \cite{MR2001h:14068}
requires the knowledge of $E$ as a cycle, that is, the computation of
the multiplicities of the components of $E$. This information is
obtained in \S\ref{multi}. Our method will essentially consist of a
local study of the families determined in \S\ref{setth}: the
multiplicity of a component will be computed by analyzing certain
numerical information carried by germs $\alpha(t)$ `marking' that
component (cf.~Definition~\ref{markers}). 

We will in fact study the inverse image $\oE$ of $E$ in the
normalization $\nG$ of $\Til \Pbb^8$: roughly speaking, the
multiplicity of the components of $\oE$ is determined by the order of
vanishing of $\mathcal C\circ\alpha(t)$ for corresponding germs
$\alpha(t)$. The number of components of $\oE$ dominating a given
component $D$ of $E$ is computed by distinguishing the contribution of
different marking germs.

Finally, the last key numerical information consists of the degree of
the components of $\oE$ over the corresponding components of $E$; this
will be obtained by studying the $\PGL(3)$ action on $\Til \Pbb^8$, and in
particular the stabilizer of a general point of each component $D$ of
$E$. Given a component $\oD$ of $\oE$ dominating $D$, we
identify a subgroup of this stabilizer (the {\em
 inessential\/} subgroup, cf.~\S\ref{inessential}), determined by the
interaction between different parametrizations of a corresponding
marking germ; the degree of $\oD$ over $D$ is the index
of this subgroup in the stabilizer (Proposition~\ref{degreetool}). 


\section{Set-theoretic description of the normal cone}\label{setth}

\subsection{}\label{preamble}
In this section we determine the different components of the
projective normal cone (PNC for short) $E$ described in
Lemma~\ref{excdivsupp}, for a given, arbitrary plane curve $\mathcal
C$ with homogeneous ideal $(F)$, where $F\in \Cbb[x,y,z]$ is a
homogeneous polynomial of degree~$d$.

The PNC can be embedded in $\Pbb^8\times\Pbb^N$, where $\Pbb^N$
parametrizes all plane curves of degree~$d$. A typical point of a
component of the PNC is in the form 
$$(\alpha(0),\lim_{t\to 0}\mathcal C\circ\alpha(t))$$
where $\alpha(t)$ is a curve germ centered at a point $\alpha(0)$ such
that $\im\alpha(0)$ is contained in $\mathcal C$. We will determine
the components of the PNC by determining a list of germs
$\alpha(t)$ which exhaust the possibilities for pairs
$(\alpha(0),\lim_{t\to 0}\mathcal C\circ\alpha(t))$ for a given
curve. Roughly speaking, we will say that two germs are equivalent
if they determine the same data $(\alpha(0),\lim_{t\to 0}\mathcal
C\circ\alpha(t))$ (see Definition~\ref{equivgerms} for the precise
notion). For a given curve $\mathcal C$ and a given germ $\alpha(t)$,
we will construct an equivalent germ in a standardized form; we will
determine which germs $\alpha(t)$ in these standard forms `contribute'
components of the PNC, in the sense that $(\alpha(0),\lim_{t\to 0}\mathcal
C\circ\alpha(t))$ belongs to exactly one component (and that a
sufficiently general such $\alpha(t)$ yields a general point of that
component), and describe that component.

\subsection{}\label{statement}
The end-result of the analysis can be stated without reference to
specific germs $\alpha(t)$. We will do so in this subsection, by listing
general points $(\alpha,\mathcal X)$ on the components of the PNC, for
a given $\mathcal C$. 

We will find five types of components: the first two will depend on
global features of $\mathcal C$, while the latter three will depend on
features of special points of $\mathcal C$ (inflection points and
singularities of the support of $\mathcal C$). The terminology
employed here matches the one in \S2 of \cite{MR2001h:14068}. In four
of the five components $\alpha$ is a rank-1 matrix, and the line
$\ker\alpha$ plays an important r\^ole; we will call this `the kernel
line'.

\begin{itemize}

\item{\em Type I.\/}
\begin{itemize} 
\item $\alpha$: a rank-2 matrix whose image is a linear component
  $\ell$ of $\mathcal C$;
\item $\mathcal X$: a fan consisting of a star of lines through the
  kernel of $\alpha$ and cutting out on the residual line $\ell'$ a
  tuple of points projectively equivalent to the tuple cut out on
  $\ell$ by the residual to $\ell$ in $\mathcal C$. The multiplicity
  of $\ell'$ in the fan is the same as the multiplicity of $\ell$ in
  $\mathcal C$.
\end{itemize}
Fans and stars are studied in \cite{MR2002d:14084}; they are items (3)
and (5) in the classification of curves with small linear orbits, in
\S1 of loc.~cit.

\item{\em Type II.\/} 
\begin{itemize}
\item $\alpha$: a rank-1 matrix whose image is a nonsingular point of
  the support $\mathcal C'$ of a nonlinear component of $\mathcal C$;
\item $\mathcal X$: a nonsingular conic tangent to the kernel line,
  union (possibly) a multiple of the kernel line. The multiplicity of
  the conic component in $\mathcal X$ equals the multiplicity of
  $\mathcal C'$ in $\mathcal C$.
\end{itemize}
Such curves are items (6) and (7) in the classification of curves with small
orbit. The extra kernel line is present precisely when $\mathcal C$ is
not itself a multiple nonsingular conic.

\item{\em Type III.\/}
\begin{itemize}
\item $\alpha$: a rank-1 matrix whose image is a point $p$ at which
  the tangent cone to $\mathcal C$ is supported on at least three
  lines;
\item $\mathcal X$: a fan with star reproducing projectively the
  tangent cone to $\mathcal C$ at $p$, and a multiple residual kernel
  line.
\end{itemize}
These limit curves are also fans, as in type~I components; but note
that type~I and type~III components are different, since for the
typical $(\alpha,\mathcal X)$ in type~I components $\alpha$ has
rank~2, while it has rank~1 for type~III components.

\item{\em Type IV.\/}
\begin{itemize}
\item $\alpha$: a rank-1 matrix whose image is a singular or
  inflection point $p$ of the support of $\mathcal C$.
\item $\mathcal X$: a curve determined by the choice of a line in the
  tangent cone to $\mathcal C$ at $p$, and by the choice of a side of
  a corresponding Newton polygon. This procedure is explained more in
  detail below.
\end{itemize}
The curves $\mathcal X$ arising in this way are items (7) through (11)
in the classification in \cite{MR2002d:14084}, and are studied
enumeratively in \cite{MR2002d:14083}.

\item{\em Type V.\/}
\begin{itemize}
\item $\alpha$: a rank-1 matrix whose image is a singular point $p$ of
  the support of $\mathcal C$.
\item $\mathcal X$: a curve determined by the choice of a line $\ell$
  in the tangent cone to $\mathcal C$ at $p$, the choice of a formal
  branch for $\mathcal C$ at $p$ tangent to $\ell$, and the choice of
  a certain `characteristic' rational number. This procedure is
  explained more in detail below.
\end{itemize}
The curves $\mathcal X$ arising in this way are item (12) in the
classification in \cite{MR2002d:14084}, and are studied enumeratively
in \cite{MR2002d:14083},~\S4.1.
\end{itemize}

Here are the details of the determination of the limit curves
$\mathcal X$ for components of type IV and V.

{\bf Type IV:\/} Let $p=\im\alpha$ be a singular or inflection point
of the support of $\mathcal C$; choose a line in the tangent cone to
$\mathcal C$ at $p$, and choose coordinates $(x:y:z)$ so that $x=0$ is
the line $\ker\alpha$, $p=(1:0:0)$, and that the selected line in the
tangent cone has equation $z=0$.
 The {\em Newton polygon\/} for $\mathcal C$ in the
chosen coordinates is the boundary of the convex hull of the union of
the positive quadrants with origin at the points $(j,k)$ for which the
coefficient of $x^iy^jz^k$ in the generator $F$ for the ideal of
$\mathcal C$ in the chosen coordinates is nonzero (see
\cite{MR88a:14001}, p.380). The part of the Newton polygon consisting
of line segments with slope strictly between $-1$ and $0$ does not
depend on the choice of coordinates fixing the flag $z=0$, $p=(1:0:0)$.

The possible limit curves $\mathcal X$ determining components of
type~IV are then obtained by choosing a side of the polygon with slope
strictly between $-1$ and $0$, and setting to~$0$ the coefficients of
the monomials in $F$ {\em not\/} on that side. These curves are
studied in \cite{MR2002d:14083}; typically, they consist of a union of
cuspidal curves. The kernel line is part of the distinguished triangle
of such a curve, and in fact it must be one of the distinguished tangents.

This procedure determines a component of the PNC, 
unless the limit curve $\mathcal X$ is supported on a conic union
(possibly) the kernel line.

{\bf Type V:\/} Let $p=\im\alpha$ be a singular point of the support
of $\mathcal C$, and let $m$ be the multiplicity of $\mathcal C$ at
$p$. Again choose a line in the tangent cone to $\mathcal C$ at $p$,
and choose coordinates $(x:y:z)$ so that $x=0$ is the kernel line,
$p=(1:0:0)$, and $z=0$ is the selected line in the tangent cone. 

We may describe $\mathcal C$ near $p$ as the union of $m$ `formal
branches'; those that are tangent to $z=0$ may be written 
$$z=f(y)=\sum_{i\ge 0} \gamma_{\lambda_i} y^{\lambda_i}$$
with $\lambda_i\in \Qbb$, $1<\lambda_0<\lambda_1<\dots$, and
$\gamma_{\lambda_i}\ne 0$.

The choices made above determine a finite set of rational numbers, which we
call the `characteristics' for $\mathcal C$ (w.r.t.~the line $z=0$):
these are the numbers $C$ such that at least two of the branches
tangent to $z=0$ agree modulo $y^C$, differ at $y^C$, and have
$\lambda_0<C$.

For a characteristic $C$, the initial exponents $\lambda_0$ and the
  coefficients $\gamma_{\lambda_0}$, $\gamma_{\frac {C+\lambda_0}2}$
  for the corresponding branches must agree. Let
  $\gamma_C^{(1)},\dots,\gamma_C^{(S)}$ be the coefficients of $y^C$
  in these $S$ branches (so that at least two of these numbers are
  distinct, by the choice of~$C$). Then $\mathcal X$ is defined by 
$$x^{d-2S}\prod_{i=1}^S\left(zx-\frac {\lambda_0(\lambda_0-1)}2
\gamma_{\lambda_0}y^2 -\frac{\lambda_0+C}2
\gamma_{\frac{\lambda_0+C}2}yx-\gamma_C^{(i)}x^2\right)\quad.$$
This is a union of `quadritangent' conics---that is, nonsingular conics
meeting at exactly one point---with (possibly) a multiple of the
distinguished tangent, which must be supported on the kernel line.

\subsection{}\label{example}
The following simple example illustrates the components described in
\S\ref{statement}: all five types are present for the curve
$$y((y^2+xz)^2-4xyz^2)=0\quad.$$
We will list five germs $\alpha(t)$, and the corresponding
limits. (This is not an exhaustive list of all the components of the
PNC for this curve.)
The specific germs used here were obtained by applying the procedures
explained in the rest of the section.

\begin{itemize}

\item{\em Type I.\/} A germ $\alpha(t)$ centered at a rank-2 matrix with
  image the linear component of the curve:
$$\alpha(t)=\begin{pmatrix}
1 & 0 & 0\\
0 & t & 0\\
0 & 0 & 1
\end{pmatrix}$$
This yields the limit
$$y\,x^2z^2\quad,$$
a fan consisting of the line $y=0$ and a star through the point
$(0:1:0)$, that is, the kernel of $\alpha(0)$.

\item{\em Type II.\/} We `aim' a one-parameter subgroup with weights
  $(1,2)$ at a nonsingular point of the curve and its tangent line:
$$\alpha(t)=
\begin{pmatrix} 
1 & 0 & 0 \\
1 & 1 & 0 \\
1 & 1 & 1
\end{pmatrix} \begin{pmatrix}
1 & 0 & 0 \\
0 & t & 0 \\
0 & 0 & t^2
\end{pmatrix} \begin{pmatrix}
1 & 0 & 0 \\
-1 & 1 & 0 \\
0 & -1 & 1
\end{pmatrix}=\begin{pmatrix}
1 & 0 & 0 \\
1-t & t & 0 \\
1-t & t-t^2 & t^2
\end{pmatrix}\quad.$$
(The curve is nonsingular at $(1:1:1)=\im\alpha(0)$, and tangent to the
line $y=z$.)
This germ yields a limit
$$x^3((x+y)^2-4xz)\quad,$$
that is, a nonsingular conic union a (multiple) tangent line supported
on the kernel line $\ker\alpha(0)$.

\item{\em Type III.\/} We aim a one-parameter subgroup with weights
  $(1,1)$ at $(0:0:1)$:
$$\alpha(t)=\begin{pmatrix}
t & 0 & 0\\
0 & t & 0\\
0 & 0 & 1
\end{pmatrix}$$
obtaining a limit of
$$xy(x-4y)z^2\quad:$$
a fan consisting of the tangent cone to $\mathcal C$ at $p$, union a
multiple kernel line.

\item{\em Type IV.\/} Considering now $\mathcal C$ at $p=(1:0:0)$,
  here is the Newton polygon w.r.t.~the line $z=0$:
$$\includegraphics[scale=1]{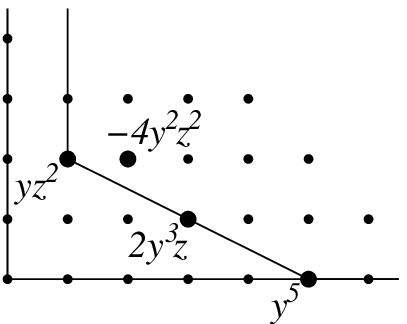}$$
It has one side with slope between $-1$ and $0$. The corresponding
germ will be a one-parameter subgroup with weights $(1,2)$:
$$\alpha(t)=\begin{pmatrix}
1 & 0 & 0\\
0 & t & 0\\
0 & 0 & t^2
\end{pmatrix}$$
yielding as limit the monomials of $F$ situated on the selected side:
$$x^2yz^2+2xy^3z+y^5=y(y^2+xz)^2\quad.$$
This is a double nonsingular conic, union a transversal line.

\item{\em Type V.\/} Finally, write formal branches for $\mathcal C$
  at $p$:
$$\left\{\aligned
y &=\phantom{-}0\\
z &=-y^2-2y^{5/2}-\dots\\
z &=-y^2+2y^{5/2}-\dots
\endaligned\right.$$
We find one characteristic $C=\frac 52$, corresponding to the second
and third branch. These branches truncate to $-y^2$; as will be
explained in \S\ref{quadritangent}, this information determines the germ
$$\alpha(t)=\begin{pmatrix}
1 & 0 & 0 \\
t^4 & t^5 & 0 \\
-t^8 & -2t^9 & t^{10}
\end{pmatrix}$$
yielding the limit
$$x(zx+y^2+2x^2)(zx+y^2-2x^2)$$
prescribed by the formula given in \S\ref{statement}. This is a pair
of quadritangent conics union a distinguished tangent
supported on the kernel line. 
\end{itemize}

\subsection{}\label{plan}
The rest of the section consists of the detailed analysis yielding the
list given in \S\ref{statement}. Our approach will be in the spirit of
Ghizzetti's paper, and indeed was inspired by reading it.
The general strategy consists of an elimination process: starting from
an arbitrary germ $\alpha(t)$, we determine the possible components
that arise unless $\alpha(t)$ is of some special kind, and keep
restricting the possibilities for $\alpha(t)$ until none is left.

Here is a guided tour of the successive reductions. First of all, we
will determine germs leading to type~I components (\S\ref{typeI});
this will account for all germs centered at a rank-2 matrix, so we
will then be able to assume that $\alpha(0)$ has rank 1. Next, we will
show (Proposition~\ref{faber},
\S\ref{faberlemma}--\S\ref{faberrefinement}) that every $\alpha$
centered at a rank-1 matrix is equivalent, in suitable coordinates, to
one in the form  
$$\alpha(t)=\begin{pmatrix}
1 & 0 & 0\\
q(t) & t^b & 0\\
r(t) & s(t)t^b & t^c
\end{pmatrix}$$
with $1\le b\le c$ and $q$, $r$, and $s$ polynomials satisfying
certain conditions. Considering the case $q\equiv r\equiv s\equiv 0$
(that is, when $\alpha$ is a `one-parameter subgroup') leads to
components of type~II, III, and IV; this is done in
\S\ref{1PS}--\S\ref{typeII}. The subtlest case, leading to type~V
components, takes the remaining
\S\ref{non1PS}--\S\ref{quadritangent}. The key step here consists of
showing that germs leading to new components are equivalent to germs
in the form
$$\begin{pmatrix}
1 & 0 & 0     \\
t^a & t^b & 0 \\
\underline{f(t^a)} & \underline{f'(t^a) t^b} & t^c
\end{pmatrix}$$
where $C=\frac ca$ is one of the characteristics considered in
\S\ref{statement}, $f$ is a corresponding branch, $b$ is determined by
the other data, and $\underline{\dots}$ stands for the truncation
modulo $t^c$. This key reduction is accomplished in
Proposition~\ref{standardform}, after substantial preparatory work.
A refinement of the reduction, given in Proposition~\ref{abc}, leads
to the definition of `characteristics' and to the description of
components of type~V given above (cf.~Proposition~\ref{typeV}).

As the process accounts for all possible germs, this will show that
the list of components given in \S\ref{statement} is exhaustive,
concluding our description of the PNC.\vskip 6pt

\subsection{}\label{lemmas}
Before starting on the path traced above, we discuss three results
that will be applied at several places in the discussion. The first
two are discussed in this subsection, and the third one in the
subsection which follows.

The first concerns a recurrent tool in establishing that a germ
$\alpha(t)$ does not contribute a component to the PNC.
As we discussed in
\S\ref{rank2obs}, this is the case for `rank-2 limits', that is,
limits that can also be obtained by germs entirely contained within
the locus of rank-2 transformations in $\Pbb^8$. Since we are looking
for germs determining components of $E$, we may ignore such `rank-2
limits' and the germs that lead to them.

\begin{lemma}\label{rank2lemma}
Assume that $\alpha(0)$ has rank~$1$.
If $\lim_{t\to 0}\mathcal C\circ\alpha(t)$ is a star with center on 
$\ker\alpha(0)$, then it is a rank-2 limit.
\end{lemma}

\begin{proof}
Assume $\mathcal X=\lim_{t\to 0}\mathcal C\circ\alpha(t)$ is a
star with center on $\ker\alpha(0)$. We may choose coordinates so that
$x=0$ is the kernel line and the generator for the ideal of
$\mathcal X$ is a polynomial in $x,y$ only. If
$$\alpha(t)=\begin{pmatrix}
a_{00}(t) & a_{01}(t) & a_{02}(t) \\
a_{10}(t) & a_{11}(t) & a_{12}(t) \\
a_{20}(t) & a_{21}(t) & a_{22}(t)
\end{pmatrix}\quad,$$
then $\mathcal X=\lim_{t\to 0}\mathcal C\circ\beta(t)$ for
$$\beta(t)=\begin{pmatrix}
a_{00}(t) & a_{01}(t) & 0 \\
a_{10}(t) & a_{11}(t) & 0 \\
a_{20}(t) & a_{21}(t) & 0
\end{pmatrix}\quad.$$
Since $\alpha(0)$ has rank~1 and kernel line $x=0$,
$$\alpha(0)=\begin{pmatrix}
a_{00}(0) & 0 & 0 \\
a_{10}(0) & 0 & 0 \\
a_{20}(0) & 0 & 0
\end{pmatrix}=\beta(0)\quad.$$
Now $\beta(t)$ is contained in the rank-2 locus, verifying the assertion.
\end{proof}

A limit $\lim_{t\to 0}\mathcal C\circ\alpha(t)$ as in the lemma will
be called a `kernel star'.

A second tool will be at the root of our reduction process: we will
replace a given germ $\alpha(t)$ with a different one in a more
manageable form, but leading to the same component of the PNC in a
very strong sense.

\begin{defin}\label{equivgerms} Two germs are `equivalent' with
respect to $\mathcal C$ if they are (possibly up to an invertible
change of parameter) fibers of a family of germs with constant center
and limit. More precisely, two germs $\alpha_0(t)$, $\alpha_1(t)$ are
equivalent if there exists a connected curve~$\mathcal H$, a regular
map $A:\mathcal H\times \Spec \Cbb[[t]] \to \Pbb^8$, two points $h_0$,
$h_1$ of $\mathcal H$, and a unit $\nu(t)\in\Cbb[[t]]$ such that:

\begin{itemize}
\item $A(h_0,t)=\alpha_0(t)$;
\item $A(h_1,t)=\alpha_1(t\nu(t))$;
\item $A(\_,0): \mathcal H \to \Pbb^8$ is constant;
\item If $\mathcal A:\mathcal H\times\Spec \Cbb[[t]] \to \Cbb^9$
is a lift of $A$, then $F\circ \mathcal A(h,t)\equiv\rho(h) G t^w 
\mod t^{w+1}$ for some $w$, 
with $\rho: \mathcal H\to\Cbb^*$ a nowhere vanishing function 
and $G$ a nonzero polynomial in $x,y,z$ independent of $h$.
\end{itemize}
\end{defin}

It is clear that Definition~\ref{equivgerms} gives an equivalence
relation on the set of germs: reflexivity and symmetry are immediate,
and transitivity is obtained by joining two families along a common
fiber (note that the limit $G$ and the weight $w$ are determined by
any of the fibers of a family, hence they are the same for any two
families extending a given germ).

Note the possible parameter change in the second condition: this
guarantees that a germ $\alpha(t)$ is equivalent to any of its
reparametrizations $\alpha(t\nu(t))$. This flexibility will 
play an important r\^ole in part of our discussion, especially
in \S\ref{inessential} and ff.

By Lemma~\ref{excdivsupp}, the third and fourth conditions amount to
the statement that all germs $\alpha_h(t)=A(h,t)$ lift to germs in
$\Til \Pbb^8$, the closure of the graph of $c$, centered at the same
point. If $\alpha_0$, $\alpha_1$ are centered at a point of $S$, then
$\alpha_0$ and $\alpha_1$ (and in fact all the intermediate
$\alpha_h$) determine the same point in the projective normal cone.
Thus, equivalent germs can be `continuously deformed' one into the
other while holding the center of their lift in $\Pbb^8\times\Pbb^N$
fixed.

Typically, the curve $\mathcal H$ will simply be a chain of affine
lines, minus some points.

Given an arbitrary germ $\alpha(t)$, we will want to produce an
equivalent and `simpler' germ. Our basic tool to produce an equivalent
germ will be the following.

\begin{lemma}\label{stequiv}
Assume $\alpha(t)\equiv\beta(t)\circ m(t)$, and that $M=m(0)$ is
invertible. Then $\alpha$ is equivalent to $\beta\circ M$ (w.r.t.~any
curve $\mathcal C$).
\end{lemma}

\begin{proof}
Write $m(t)=M+tm_1(t)$, and take $\mathcal H=\Abb^1$, $\nu(t)=1$. Let
$$\alpha_h(t)= A(h,t):=\beta(t)\circ(M+h t m_1(t))\quad.$$ 
Then
\begin{itemize}
\item $\alpha_0=\beta\circ M$,
\item $\alpha_1=\alpha$, and 
\item $\alpha_h(0)=\beta(0)\circ M$ does not depend on $h$.
\end{itemize}
For any given $F$,
$$F\circ\beta(t)=t^w G+ t^{w+1} G_1(t)$$
with $G=\lim_{t\to 0} F\circ \beta(t)$. Then
$$F\circ\alpha_h(t)=F\circ\beta(t)\circ (M+h t m_1(t))=t^w(G+t G_1(t))
\circ(M+h t m_1(t))=t^wG\circ M+\text{h.o.t.}$$
The term $G\circ M$ is not $0$, since $M$ is invertible, and does not
depend on $h$. 

Thus $\alpha=\alpha_1$ is equivalent to $\beta\circ M=\alpha_0$
according to Definition~\ref{equivgerms}, as needed.
\end{proof}

\subsection{}\label{formalbranches}
The third preliminary item concerns formal branches of $\mathcal C$ at
a point $p$, cf.~\cite{MR88a:14001} and \cite{MR1836037}, Chapter~6
and~7. Choose affine coordinates $(y,z)=(1:y:z)$ so that $p=(0,0)$,
and let $\Phi(y,z)=F(1:y:z)$ be the generator for the ideal of
$\mathcal C$ in these coordinates. Decompose $\Phi(y,z)$ in
$\Cbb[[y,z]]$:
$$\Phi(y,z)=\Phi_1(y,z)\cdot\cdots\cdot \Phi_r(y,z)$$
with $\Phi_i(y,z)$ irreducible power series. 
These define the {\em irreducible branches\/} of $\mathcal C$
at~$p$. Each $\Phi_i$ has a unique tangent line at $p$; if this
tangent line is {\em not\/} $y=0$, by the Weierstrass preparation
theorem we may write (up to a unit in $\Cbb[[y,z]]$) $\Phi_i$ as a
monic polynomial in $z$ with coefficients in $\Cbb[[y]]$, of degree
equal to the multiplicity $m_i$ of the branch at $p$
(cf.~for example \cite{MR1836037}, \S6.7).
If $\Phi_i$ {\em is\/} tangent to $y=0$, we may likewise write it as a
polynomial in $y$ with coefficients in $\Cbb[[z]]$; {\em mutatis
  mutandis,\/} the discussion which follows applies to this case as
well.

Concentrating on the first case, let
$$\Phi_i(y,z)\in \Cbb[[y]][z]$$
be a monic polynomial of degree $m_i$, defining an irreducible branch
of $\mathcal C$ at $p$, not tangent to $y=0$. Then $\Phi_i$ splits
(uniquely) as a product of linear factors over the ring $\Cbb[[y^*]]$
of power series with {\em rational nonnegative\/} exponents:
$$\Phi_i(y,z)=\prod_{j=1}^{m_i} \left(z-f_{ij}(y)\right)\quad,$$
with each $f_{ij}(y)$ in the form
$$f(y)=\sum_{k\ge 0} \gamma_{\lambda_k} y^{\lambda_k}$$
with $\lambda_k\in \Qbb$, $1\le \lambda_0<\lambda_1<\dots$, and
$\gamma_{\lambda_k}\ne 0$. We call each such $z=f(y)$ a {\em formal
  branch\/} of $\mathcal C$ at $p$. The branch is {\em tangent to
  $z=0$\/} if the dominating exponent $\lambda_0$ is~$>1$. The terms
$z-f_{ij}(y)$ in this decomposition are the Puiseux series for
$\mathcal C$ at $p$.

Summarizing: if $\mathcal C$ has multiplicity $m$ at $p$ then
$\mathcal C$ splits into $m$ formal branches at~$p$. In
\S\ref{non1PS} and ff.~we will want to determine $\lim_{t\to
  0}\mathcal C\circ\alpha(t)$ as a union of `limits' of the individual
formal branches at $p$. The difficulty here resides in the fact that
we cannot perform an arbitrary `change of variable' in a power series
with fractional exponents. In the case in which we will need to do
this, however, $\alpha(t)$ will have the following special form:
$$\alpha(t)=\begin{pmatrix}
1 & 0 & 0 \\
t^a & t^b & 0 \\
r(t) & s(t)t^b & t^c
\end{pmatrix}$$
with $a<b\le c$ positive integers and $r(t)$, $s(t)$ polynomials
(satisfying certain restrictions, which are immaterial here).
 We will circumvent
the difficulty we mentioned by the following {\em ad hoc\/} definition.

\begin{defin}\label{branchlimit} 
The {\em limit\/} of a formal branch $z=f(y)$, along a germ
$\alpha(t)$ as above, is defined by the dominant term in
$$(r(t)+s(t)t^by+t^cz)-f(t^a)-f'(t^a)t^by-f''(t^a)t^{2b}\frac{y^2}2
-\cdots$$
where $f'(y)=\sum \gamma_k\lambda_k y^{\lambda_k-1}$ etc.
\end{defin}

By `dominant term' we mean the coefficient of the lowest power of
$t$ after cancellations. This coefficient is a polynomial in $y$ and $z$,
giving the limit of the branch according to our definition.

Of course we need to verify that this definition behaves as expected,
that is, that the limit of $\mathcal C$ is the union of the limits of
its individual branches. We do so in the following lemma.

\begin{lemma}\label{limitlemma}
Let $\Phi(y,z)\in \Cbb[[y]][z]$ be a monic polynomial,
$$\Phi(y,z)=\prod_i \left(z-f_i(y)\right)$$
a decomposition over $\Cbb[[y^*]]$, and let $\alpha(t)$ be as
above. Then the dominant term in $\Phi\circ\alpha(t)$ is the product
of the limits of the branches $z=f_i(y)$ along $\alpha$, defined as in
Definition~\ref{branchlimit}.
\end{lemma}

\begin{proof}
We can `clear the denominators' in the exponents in $f_i$, by writing
$$\Phi(T^m,z)=\prod_i(z-\varphi_i(T))$$
where $\varphi_i(T)\in \Cbb[[T]]$ and $\varphi_i(T)=f_i(T^m)$. For an
integer $\ell$ such that $\ell a/m$ is integer, we may write
$$t=S^\ell\quad, \quad T^m=T(S)^m=S^{\ell a}+S^{\ell b}y=S^{\ell a}
(1+S^{\ell(b-a)}y)$$
with $T(S)\in \Cbb[[S,y]]$: explicitly, 
$$T(S)=S^{\frac{\ell a}m} \left(1+\frac 1m S^{\ell(b-a)}y + \frac 1m
\left(\frac 1m-1\right) S^{\ell 2(b-a)}\frac{y^2}2+ \cdots\right)$$
The dominant term (w.r.t.~$t$) in
$$\Phi\circ\alpha(t)=\Phi(t^a+t^by,r(t)+s(t)t^by+t^cz)$$
equals the dominant term (w.r.t.~$S$) in
\begin{multline*}
\Phi(S^{\ell a}+S^{\ell b}y,r(S^\ell)+s(S^\ell)S^{\ell b}y+S^{\ell
  c}z) =\Phi(T(S)^m,r(S^\ell)+s(S^\ell)S^{\ell b}y+S^{\ell c}z)\\
=\prod_i\left( (r(S^\ell)+s(S^\ell)S^{\ell b}y+S^{\ell c}z)
-\varphi_i(T(S))\right)
\end{multline*}
Thus the dominant term in $\Phi\circ\alpha(t)$ is the product of the
dominant terms in the factors
$$(r(S^\ell)+s(S^\ell)S^{\ell b}y+S^{\ell c}z)
-\varphi_i(T(S))$$
and we have to verify that the dominant term here agrees with the one
in Definition~\ref{branchlimit}.

For this, we use a `Taylor expansion' of $\varphi_i$. Write
$$\varphi_i(T(S))=\sum_{k\ge 0}
\frac{\partial^k\varphi_i(T(S))}{\partial y^k}|_{y=0}
\frac{y^k}{k!}\quad.$$
We claim that
$$\frac{\partial^k\varphi_i(T(S))}{\partial y^k}|_{y=0} =
f_i^{(k)}(S^{\ell a})S^{\ell b k}\quad:$$
indeed, this is immediately checked for $f_i(y)=y^\lambda$, hence
holds for any $f_i$.

Therefore we have
$$\varphi_i(T(S))=\sum_{k\ge 0} f_i^{(k)}(S^{\ell a}) \frac{S^{\ell b k}
y^k}{k!}$$
or, recalling $t=S^\ell$:
$$\varphi_i(T(S))=f_i(t^a)+f'_i(t^a) t^by+f''_i(t^a)t^{2b}\frac {y^2}2
+ \dots$$
This shows that 
$$(r(S^\ell)+s(S^\ell)S^{\ell b}y+S^{\ell c}z)-\varphi_i(T(S))$$
is in fact given by 
$$(r(t)+s(t) t^by+t^c z)-\left(f_i(t^a)+f'_i(t^a) t^by+f''_i(t^a)
t^{2b} \frac {y^2}2 + \dots\right)\quad,$$
and in particular the dominant terms in the two expressions must
match, as needed.
\end{proof}

The gist of this subsection is that we may use formal branches for
$\mathcal C$ at $p$ in order to compute the limit of $\mathcal C$
along germs $\alpha(t)$ of the type used above, provided that the
limit of a branch is computed by using the formal Taylor expansion
given in Definition~\ref{branchlimit}. This fact will be used several
times in \S\ref{non1PS} and ff.

\subsection{}\label{typeI}
We are finally ready to begin the discussion leading to the list of
components given in \S\ref{statement}.

Applying Lemma~\ref{excdivsupp} amounts to studying germs $\alpha(t)$
in $\Pbb^8$, centered at matrices $\alpha$ with image contained in
$\mathcal C$---these are precisely the matrices in the base locus $S$
of the rational map $c$ introduced in \S\ref{prelim}. The
corresponding component of $E$ is determined by the center of the
germ, and the limit.

As we may assume that $\alpha(0)$ is contained in $\mathcal C$, we may
assume that $\rk\alpha(0)=1$ or $2$. We first consider the case of
germs centered at a rank-2 matrix $\alpha$, hence with image equal to
a linear component of $\mathcal C$. We will show that any germ
centered at such a matrix leads to a point in the component of type~I
listed in \S\ref{statement}.

Write 
$$\alpha(t)=\alpha(0)+t\beta(t)\quad,$$
where $\alpha(0)$ has rank~2 and image defined by the linear
polynomial $L$; thus, we may write the generator of the ideal of
$\mathcal C$ as
$$F(x,y,z)=L(x,y,z)^m G(x,y,z)$$
with $L$ not a factor of $G$. The curve defined by $G$ is the
`residual of $L$ in $\mathcal C$'.

\begin{prop}\label{globallin}
The limit $\lim_{t\to 0} \mathcal C\circ\alpha(t)$ is a fan consisting
of an $m$-fold line $\ell$, supported on $\lim_{t\to
  0}L\circ\beta(t)$, and a star of lines through the point
$\ker\alpha$. This star reproduces projectively the tuple cut out on
$L$ by the residual of $L$ in $\mathcal C$.
\end{prop}

The terminology of stars and fans was introduced in
\cite{MR2002d:14084}, \S2.1. Here the $m$-fold line $\ell$ may contain
the point $\ker\alpha$, in which case the fan degenerates to a star.

\begin{proof}
Write
$$\alpha(t)=\alpha(0)+t\beta(t)\quad,$$
then the ideal of $\mathcal C\circ\alpha(t)$ is generated by
$$L(\alpha(t))^mG(\alpha(t))=t^m L(\beta(t))^mG(\alpha(0)+t\beta(t))
\quad,$$
since $L$ is linear and vanishes along the image of $\alpha(0)$. As
$t$ approaches 0, $L(\beta(t))^m$ converges to an $m$-fold line, while
the other factor converges to $G(\alpha(0))$, yielding the statement.
\end{proof}

A simple dimension count shows that the limits arising as in
Proposition~\ref{globallin} {\em do\/} produce components of the
projective normal cone. Indeed, matrices with image contained in a given line
form a $\Pbb^5$; for any given such matrix, the limits obtained
consist of a fixed star through the kernel, plus a (multiple) line
varying freely, accounting for 2 extra dimensions. These components
are the components of type I described in \S\ref{statement} (also
cf.~\cite{MR2001h:14068}, \S2, Fact~2(i)).

\subsection{}\label{faberlemma}
Having taken into account the case in which the center $\alpha(0)$ may
be a rank-2 matrix, we are reduced to considering germs $\alpha(t)$
with $\alpha(0)$ of rank~1. Proposition~\ref{faber} below will allow
us to further assume that $\alpha(t)$ has a particularly simple (and
polynomial) expression. 

We define the degree of the zero polynomial to be $-\infty$. We denote
by $v$ the `valuation' of a power series or polynomial, that is, its
order of vanishing at $0$; we define $v(0)$ to be $+\infty$.

\begin{prop}\label{faber}
With a suitable choice of coordinates, any germ $\alpha$ is
equivalent to a product
$$\begin{pmatrix}
1 & 0 & 0 \\
q & 1 & 0 \\
r & s & 1
\end{pmatrix} \cdot \begin{pmatrix}
t^a & 0 & 0 \\
0 & t^b & 0 \\
0 & 0 & t^c
\end{pmatrix}$$
with
\begin{itemize}
\item $a\le b\le c$ integers, $q,r,s$ polynomials;
\item $\deg(q)<b-a$, $\deg(r)<c-a$, $\deg(s)<c-b$;
\item $q(0)=r(0)=s(0)=0$.
\end{itemize}
If further $b=c$ and $q$, $r$ are not both zero, then we may assume that
$v(q)<v(r)$. 
\end{prop}

The proof of this proposition requires a few preliminary
considerations.

\subsection{}\label{sourcetarget}
The space $\Pbb^8$ of $3\times 3$ matrices considered above is, more
intrinsically, the projective space $\Pbb \Hom(V,W)$, where $V$ and
$W$ are vector spaces of dimension~3. The generator $F$ of the ideal
of a plane curve of degree $d$ is then an element of $\Sym^d
W^*$; for $\varphi\in \Hom(V,W)$, the composition $F\circ\varphi$ (if
nonzero) is the element of $\Sym^d V^*$ generating the ideal of
$\mathcal C\circ\varphi$. We denote by $\PGL(V,W)$ the Zariski open subset of
$\Pbb\Hom(V,W)$ consisting of invertible transformations,
and write $\PGL(V)$ for $\PGL(V,V)$ (which is a group under
composition).

Homomorphisms $\lambda$ of $\Cbb^*$ to $\PGL(V)$ will be called `1-PS'
(as in: `1-parameter subgroups'), as will be called their extensions
$\Cbb \to \Pbb\Hom(V,V)$. Recall that every 1-PS can be written as
$$t\mapsto \begin{pmatrix}
t^a & 0 & 0\\
0 & t^b & 0\\
0 & 0 & t^c
\end{pmatrix}$$
after a suitable choice of coordinates in $V$, where $a\le b\le c$ are
integers (and $a$ may in fact be chosen to equal~$0$). Thus we may
view a 1-PS as a $\Cbb((t))$-valued point of $\PGL(V)\subset
\Pbb\Hom(V,V)$. The following lemma shows that these are the basic
constituents of every germ $\alpha(t)$.

\begin{lemma}\label{GIT}
Every germ $\alpha(t)$ in $\Pbb\Hom(V,W)$ is equivalent to a germ
$$H\circ h_1\circ\lambda\quad,$$
where:
\begin{itemize}
\item $H$ is a constant invertible linear transformation $V\to W$; 
\item $h_1$ is a $\Cbb[[t]]$-valued point of $\PGL(V)$ with
  $h_1(0)=\text{Id}_V$; and
\item $\lambda$ is a 1-PS.
\end{itemize}
\end{lemma}

\begin{proof}
Every germ $\alpha$ can be written (cf.~\cite{MR86a:14006}, p.53) as a
composition:
$$\xymatrix{
V \ar@/_1pc/[rrr]_\alpha \ar[r]^k & U \ar[r]^{\lambda'} & U
\ar[r]^{h'} & W}$$
where $k$ and $h'$ are $\Cbb[[t]]$-valued points of $\PGL(V,U)$,
$\PGL(U,W)$ respectively and $\lambda'$ is a 1-PS. In particular,
$K=k(0)$ is an invertible linear transformation; by
Lemma~\ref{stequiv}, the composition is equivalent to the composition
$$\xymatrix{
V \ar[r]^K & U \ar[r]^{\lambda'} & U \ar[r]^{h'} & W
}\quad,$$
which can be written as
$$\xymatrix{
V \ar@/_1pc/[rrr]_\lambda \ar[r]^K & U \ar[r]^{\lambda'} & U
\ar[r]^{K^{-1}} & V \ar[r]^K & U \ar[r]^{h'} & W
}\quad.$$
Here $\lambda$ is again a 1-PS, as a conjugate of a 1-PS by a constant
transformation. The statement follows by writing $h'\circ K=H\circ
h_1$ as prescribed.
\end{proof}

Now we choose coordinates in $V$ so that $\lambda$ is diagonal:
$$\lambda=\begin{pmatrix} 
t^{a} & 0 & 0 \\
0 & t^{b} & 0 \\
0 & 0 & t^{c}
\end{pmatrix}$$
with $a\le b\le c$ integers; thus we may view $h_1$ and $\lambda$ as
matrices, and we are interested in putting $h_1$ in a `standard' form.

\begin{lemma}\label{MPI} Let
$$ h_1 =\left(\begin{matrix}
u_1 & b_1 & c_1 \\
a_2 & u_2 & c_2 \\
a_3 & b_3 & u_3
\end{matrix} \right) 
$$
be a $\Cbb[[t]]$-valued point of $\PGL(V)$, such that
$h_1(0)=I_3$. Then $h_1$ can be written as a product $h_1=h\cdot j$
with
$$
h=\begin{pmatrix}
1 & 0 & 0 \\
q & 1 & 0 \\
r & s & 1 \end{pmatrix}
\quad,\quad j=\left(\begin{matrix}
v_1 & e_1 & f_1 \\
d_2 & v_2 & f_2 \\
d_3 & e_3 & v_3 \end{matrix} \right)
$$
with $q$, $r$, $s$ {\em polynomials,\/} satisfying
\begin{enumerate}
\item $h(0)=j(0)=I_3$;
\item $\deg(q)<b-a$, $\deg(r)<c-a$, $\deg(s)<c-b$;
\item\label{refMPI} $d_2\equiv0\pmod{t^{b-a}}$, $d_3\equiv0\pmod{t^{c-a}}$, 
$e_3\equiv0\pmod{t^{c-b}}$.
\end{enumerate}
\end{lemma}

\begin{proof} Obviously $v_1=u_1, e_1=b_1$ and $f_1=c_1$. Use division 
with remainder to write 
$$ v_1^{-1}a_2=D_2t^{b-a}+q $$
with $\deg(q)<b-a$, and let $d_2=v_1D_2t^{b-a}$ (so that $qv_1+d_2=a_2$).
This defines $q$ and $d_2$, and uniquely determines $v_2$ and $f_2$.
(Note that $q(0)=d_2(0)=f_2(0)=0$ and that $v_2(0)=1$.)

Similarly, we let $r$ be the remainder of
$$(v_1v_2-e_1d_2)^{-1}(v_2a_3-d_2b_3)$$
under division by $t^{c-a}$; and $s$ be the remainder of
$$(v_1v_2-e_1d_2)^{-1}(v_1b_3-e_1a_3)$$
under division by $t^{c-b}$.

Then $\deg(r)<c-a$, $\deg(s)<c-b$ and $r(0)=s(0)=0$; moreover, we have
$$ v_1r+d_2s\equiv a_3\pmod{t^{c-a}},\qquad e_1r+v_2s\equiv b_3
\pmod{t^{c-b}},$$
so we take $d_3=a_3-v_1r-d_2s$, $e_3=b_3-e_1r-v_2s$. This defines $r$,
$s$, $d_3$ and $e_3$, and uniquely determines $v_3$.
\end{proof}

\subsection{}
We are now ready to prove Proposition~\ref{faber}. We have written a
germ equivalent to $\alpha$ as
$$H\cdot h\cdot j \cdot \lambda$$
with notations as above. Now, by (\ref{refMPI}) in Lemma~\ref{MPI} we
have $j\cdot \lambda=\lambda \cdot \ell$ for $\ell$ with entries in
$\Cbb[[t]]$, and $L=\ell(0)$ lower triangular, with 1's on the
diagonal. By Lemma~\ref{stequiv} this germ is equivalent to
$$H\cdot h\cdot \lambda\cdot L=(H\cdot L)\cdot L^{-1}\cdot(h\cdot
\lambda)\cdot L\quad.$$
We change coordinates in $V$ by $L^{-1}$, so that $L^{-1}\cdot(h\cdot
\lambda)\cdot L$ has matrix representation $h\cdot \lambda$. Finally,
we choose coordinates in $W$ so that $H\cdot L=I_3$, completing the
proof of the first part of Proposition~\ref{faber}.

If $b=c$, then the condition that $\deg s<c-b=0$ forces
$s=0$. Conjugating by
$$\begin{pmatrix}
1 & 0 & 0\\
0 & 0 & 1\\
0 & 1 & 0
\end{pmatrix}$$
interchanges $q$ and $r$; so we may assume $v(q)\le v(r)$ if $q$
and $r$ are not both~0. Conjugating by
$$\begin{pmatrix}
1 & 0 & 0\\
0 & 1 & 0\\
0 & u & 1
\end{pmatrix}$$
replaces $r$ by $uq+r$, allowing us to force $v(q)<v(r)$, and
completing the proof of Proposition~\ref{faber}.\hfill\qed

\subsection{}\label{faberrefinement}
By Proposition~\ref{faber}, and scaling the entries in the 1-PS so
that $a=0$, an arbitrary germ $\alpha$ is equivalent to one that, with
a suitable choice of coordinates, can be written as
$$\alpha(t)=\begin{pmatrix}
1 & 0 & 0\\
q(t) & t^b & 0\\
r(t) & s(t)t^b & t^c
\end{pmatrix}$$
with $0\le b\le c$ and $q$, $r$, and $s$ polynomials satisfying
certain conditions. We may in fact assume that $b>0$, since we are
already reduced to the case in which $\alpha(0)$ is a rank-1
matrix. If $b>0$, then the center $\alpha(0)$ is the matrix
$$\alpha(0)=\begin{pmatrix}
1 & 0 & 0\\
0 & 0 & 0\\
0 & 0 & 0
\end{pmatrix}$$
with image the point $(1:0:0)$ and kernel the line $x=0$.

Further, if $(1:0:0)$ is not a point of the curve $\mathcal C$ then
  $\lim_{t\to 0} \mathcal C\circ\alpha(t)$ is simply a multiple kernel
  line with ideal $(x^{\deg \mathcal C})$. Thus we may assume that
  $p=(1:0:0)$ is a point of $\mathcal C$. {\em In what follows, we
  will assume that $\alpha$ is a germ in the standard form given
  above, and all these conditions are satisfied.\/}

One last remark will be needed later in the section: if the polynomial
$q(t)$ is known to be nonzero, then Proposition~\ref{faber} admits the
following refinement.

\begin{lemma}\label{faberqnot0}
If $q\not\equiv 0$ in $\alpha(t)$, then $\alpha(t)$ is equivalent to a germ
$$\begin{pmatrix}
1 & 0 & 0     \\
t^a & t^b & 0 \\
r_1(t) & s_1(t)t^b & t^c\end{pmatrix}
\cdot\begin{pmatrix}
1 & 0 & 0 \\
0 & u & 0 \\
r & s & v 
\end{pmatrix}\quad,$$
with 
\begin{itemize}
\item $a<b\le c$ positive integers;
\item $r_1(t)$ and $s_1(t)$ polynomials of degree $<c$, $<(c-b)$
  respectively and vanishing at $t=0$; and
\item $u,r,s,v\in \Cbb$, with $uv\ne 0$.
\end{itemize}
If further $b=c$, then we may assume $a<v(r_1)$.
\end{lemma}

\begin{proof}
As $q\not\equiv 0$, we may write $q(t)=\tau(t)^a$ for $a=v(q)$ (so $0<a<b$)
and with $\tau(t)\in\Cbb[[t]]$, such that $\tau(t)/t$ is a unit in
$\Cbb[[t]]$. Expressing $t$ in terms of $\tau$, we can set
$$\beta(\tau)=\begin{pmatrix} 1 & 0 & 0\\
\tau^a & u(\tau) \tau^b & 0\\
\overline r(\tau) & \overline s(\tau) u(\tau) \tau^b & v(\tau) \tau^c
\end{pmatrix}$$
so that $\alpha(t)=\beta(\tau(t))$, for suitable $\overline r(\tau)$,
$\overline s(\tau)$, and invertible $u(\tau)$, $v(\tau)$ in
$\Cbb[[\tau]]$.

Since $\alpha(t)$ and $\beta(t)$ only differ by a change of parameter,
they are equivalent in the sense of Definition~\ref{equivgerms}.

Next, define $\rho(t),\sigma(t)\in\Cbb[[t]]$ so that
$$\overline r(t)=\underline{\overline r(t)}+\rho(t) t^c\quad,\quad 
\overline s(t) t^b=\underline{\overline s(t)t^b}+\sigma(t) t^c$$
with $\underline{\overline r(t)}$, $\underline{\overline s(t)t^b}$
polynomials of degree less than $c$,
and observe that then
$$\beta(t)=\begin{pmatrix} 
1 & 0 & 0\\
t^a & u(t) t^b & 0\\
\overline r(t) & \overline s(t) u(t) t^b & v(t) t^c
\end{pmatrix}=\begin{pmatrix}
1 & 0 & 0\\
t^a & t^b & 0\\
\underline{\overline r(t)} & \underline {\overline s(t)t^b} & t^c
\end{pmatrix} \cdot \begin{pmatrix}
1 & 0 & 0 \\
0 & u(t) & 0\\
\rho(t) & \sigma(t) u(t) & v(t)
\end{pmatrix}\quad.$$
The rightmost matrix is invertible at~0, so by Lemma~\ref{stequiv}
$\beta(t)$ (and hence $\alpha(t)$) is equivalent to
$$\begin{pmatrix}
1 & 0 & 0\\
t^a & t^b & 0\\
\underline{\overline r(t)} & \underline {\overline s(t)t^b} & t^c
\end{pmatrix} \cdot \begin{pmatrix}
1 & 0 & 0 \\
0 & u & 0\\
r & s & v
\end{pmatrix}$$
where $u=u(0)$, $r=\rho(0)$, $s=\sigma(0)u(0)$, and $v=v(0)$.
We have $uv\ne 0$ as both $u(t)$ and $v(t)$ are invertible. 

To obtain the stated form, let $r_1(t)=\underline{\overline r(t)}$ and
$s_1(t)$ so that $s_1(t)t^b=\underline {\overline s(t)t^b}$. Then
$r_1(t)$ and $s_1(t)$ are polynomials of degree $<c$, $<(c-b)$
respectively, and $r_1(0)=s_1(0)=0$ as an immediate consequence of
$r(0)=s(0)=0$.

Finally, note that $a=v(q)$ and $v(r_1)=v(r)$; if $b=c$, then
we may assume $v(q)<v(r)$ by Proposition~\ref{faber}, and hence
$a<v(r_1)$ as needed.
\end{proof}

The form obtained in Lemma~\ref{faberqnot0} will be needed in a key
reduction (Proposition~\ref{standardform}) later in the section.
The effect of the constant factor on the right in the germ appearing
in the statement of Lemma~\ref{faberqnot0} is simply to translate the
limit (by an invertible transformation fixing the flag consisting of
the line $x=0$ and the point $(0:0:1)$). Thus this factor will
essentially be immaterial in the considerations in this section. 

\subsection{}\label{affine}
In the following, it will be convenient to switch to affine
coordinates centered at the point $(1:0:0)$: we will denote by $(y,z)$
the point $(1:y:z)$; as we just argued, we may assume that the curve
$\mathcal C$ contains the origin $p=(0,0)$. We write
$$F(1:y:z)=F_m(y,z)+F_{m+1}(y,z)+\cdots +F_d(y,z)\quad,$$
with $d=\deg \mathcal C$, $F_i$ homogeneous of degree $i$, and $F_m\ne
0$. Thus, $F_m(y,z)$ generates the ideal of the {\em tangent cone\/}
of $\mathcal C$ at $p$.

\subsection{}\label{1PS}
In the next three
subsections we consider the case in which $q=r=s=0$,
that is, in which $\alpha(t)$ is itself a 1-PS:
$$\alpha(t)=\begin{pmatrix}
1 & 0 & 0\\
0 & t^b & 0\\
0 & 0 & t^c
\end{pmatrix}$$
with $1\le b \le c$. Also, we may assume that $b$ and $c$ are coprime:
this only amounts to a reparametrization of the germ by $t \mapsto
t^{1/d}$, with $d=gcd(b,c)$; the new germ is not equivalent to the old
one in terms of Definition~\ref{equivgerms}, but clearly achieves the
same limit.

Germs with $b=c(=1)$ lead to components of type~III:
\begin{prop}\label{tgcone}
If $q=r=s=0$ and $b=c$, then $\lim_{t\to 0} \mathcal C\circ\alpha(t)$
is a fan consisting of a star projectively equivalent to the tangent cone
to $\mathcal C$ at $p$, and of a residual $(d-m)$-fold line supported
on $\ker\alpha$.
\end{prop}

\begin{proof}
The composition $F\circ\alpha(t)$ is
$$F(x:t^by:t^bz)=t^{bm}x^{d-m}F_m(y,z)+t^{b(m+1)}x^{d-(m+1)}
F_{m+1}(y,z)+\cdots+ t^{dm}F_d(y,z)\quad.$$
By definition of limit, $\lim_{t\to 0}\mathcal C\circ\alpha(t)$ has ideal
$(x^{d-m}F_m(y,z))$, proving the assertion.
\end{proof}

A dimension count (analogous to the one in \S\ref{typeI}) shows that
the limits found in Proposition~\ref{tgcone} contribute a component to
the projective normal cone when the star is supported on three or more
lines. These are the components `of type III' in the terminology of
\S\ref{statement}; also cf.~\cite{MR2001h:14068}, \S2, Fact~4(i).

\subsection{}\label{Newton}
More components may arise due to 1-PS with $b<c$, but only if
$\mathcal C$ is in a particularly special position relative to
$\alpha$.

\begin{lemma}
If $q=r=s=0$ and $b<c$, and $z=0$ is not contained in the tangent cone
to $\mathcal C$ at $p$, then $\lim_{t\to 0} \mathcal C\circ\alpha(t)$
is supported on a pair of lines.
\end{lemma}

\begin{proof}
The condition regarding $z=0$ translates into $F_m(1,0)\ne 0$. 
Applying $\alpha(t)$ to~$F$, we find:
$$F(x:t^by:t^cz)=t^{bm}x^{d-m} F_m(y,t^{c-b}z)+t^{b(m+1)} x^{d-(m+1)}
F_{m+1}(y,t^{c-b}z)+\cdots$$
Since $F_m(1,0)\ne 0$, the dominant term on the right-hand-side is
$x^{d-m}y^m$, proving the assertion.
\end{proof}

By Lemma~\ref{rank2lemma}, these limits do not contribute components
to the projective normal cone.

Components that do arise due to 1-PS with $b<c$ may be described in
terms of the {\em Newton polygon\/} for $\mathcal C$ at $(0,0)$,
relative to the line $z=0$, which we may now assume (by the preceding
lemma) is part of the tangent cone to $\mathcal C$ at $p$. The Newton
polygon for $\mathcal C$ in the chosen coordinates is the boundary of
the convex hull of the union of the positive quadrants with origin at
the points $(j,k)$ for which the coefficient of $x^iy^jz^k$ in the
equation for $\mathcal C$ is nonzero (see \cite{MR88a:14001},
p.380). The part of the Newton polygon consisting of line segments
with slope strictly between $-1$ and $0$ does not depend on the choice
of coordinates fixing the flag $z=0$, $p=(0,0)$.

\begin{prop}\label{Newtonsides}
Assume $q=r=s=0$ and $b<c$.
\begin{itemize}
\item If $-b/c$ is not a slope of the Newton polygon for $\mathcal C$,
  then the limit $\lim_{t\to 0} \mathcal C\circ\alpha(t)$ is supported
  on (at most) three lines. Such limits do not contribute components to
  the projective normal cone.
\item If $-b/c$ is a slope of a side of the Newton polygon for
  $\mathcal C$, then the ideal of the limit $\lim_{t\to 0}\mathcal C
  \circ\alpha(t)$ is generated by the polynomial obtained by 
  setting to~$0$ the coefficients of
  the monomials in $F$ {\em not\/} on that side.
  Such polynomials are in the form
$$G=x^{\overline q}y^rz^q \prod_{j=1}^S(y^c+\rho_j x^{c-b}z^b)$$
\end{itemize}
\end{prop}

\begin{proof}
For the first assertion, simply note that under the stated hypotheses
only one monomial in $F$ is dominant in $F\circ\alpha(t)$; hence, the
limit is supported on the union of the coordinate axes. A simple
dimension count shows that such limits may only span a 6-dimensional
locus in $\Pbb^8\times\Pbb^N$, so they do not determine a component of
the projective normal cone.

The second assertion is analogous: the dominant terms in
$F\circ\alpha(t)$ are precisely those on the side of the Newton polygon
with slope equal to $-b/c$. It is immediate that the resulting
polynomial can be factored as stated.
\end{proof}

Limits arising as in the second part of
Proposition~\ref{Newtonsides} are the curves studied in
\cite{MR2002d:14083}, and appear as items (6) through (11) 
in the classification in \S1 of \cite{MR2002d:14084}.
The number $S$ of `cuspidal' factors in
$G$ is the number of segments cut out by the integer lattice on the
selected side of the Newton polygon.

Assume the point $p=(1:0:0)$ is a singular or an inflection point of
the support of~$\mathcal C$. If $b/c\ne 1/2$, then the corresponding limit will
contribute a component to the PNC:
indeed, the orbit of the corresponding limit curve has dimension~7. 
If $b/c=1/2$, then a dimension count shows that the corresponding limit 
will contribute a component to the PNC
unless it is supported on a conic union (possibly) the kernel line.

These are the components of type~IV in \S\ref{statement},
also cf.~\cite{MR2001h:14068}, \S2, Fact~4(ii).

\subsection{}\label{typeII}
If $p$ is a {\em nonsingular, non-inflectional\/} point of the support
of $\mathcal C$, then the Newton polygon consists of a single side
with slope $-1/2$, and the polynomial $G$ in the statement of
Proposition~\ref{Newtonsides} reduces to
$$x^{d-2S}(y^2+\rho x z)^S\quad,$$ 
that is, a (multiple) conic union a (multiple) tangent line supported
on $\ker\alpha$; here $S$ is the multiplicity of the corresponding
component of $\mathcal C$. The orbit of this limit curve has
dimension~6; but as there is one such limit at almost all points
of the support of every nonlinear component of $\mathcal C$, the
collection of these limits span one component of the projective normal
cone for each nonlinear component of $\mathcal C$. These components
are the components of type~II in \S\ref{statement}, also
cf.~\cite{MR2001h:14068}, Fact~2(ii).

\begin{example}
Consider the `double cubic' with (affine) ideal generated by 
$$(y^2+z^3+z)^2=y^4+2y^2z^3+2y^2z+z^6+2z^4+z^2\quad.$$
Its Newton polygon consists of one side, with slope $-1/2$:
$$\includegraphics[scale=1]{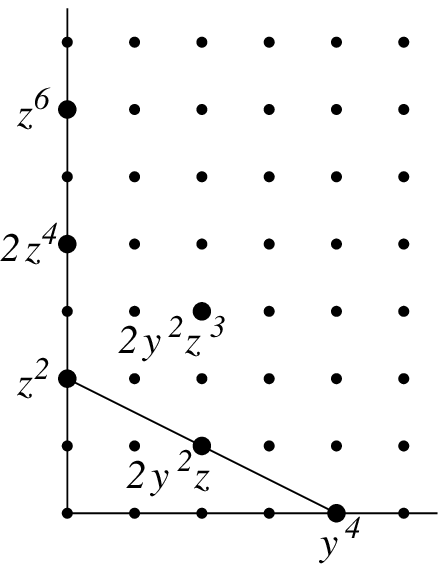}$$
The limit by the 1-PS
$$\begin{pmatrix}
1 & 0 & 0\\
0 & t & 0\\
0 & 0 & t^2
\end{pmatrix}$$
consists, according to the preceding discussion, of the part of the
above polynomial supported on the side. This may be checked directly:
$$(x(ty)^2+(t^2 z)^3+x^2(t^2z))^2 =x^2t^4(y^4+2xy^2z+x^2z^2)+
xt^8(2y^2z^3+2xz^4) +t^{12}z^6$$
the dominant terms in this expression are
$x^2y^4+2x^3y^2z+x^4z^2=x^2(y^2+xz)^2$. The support of the limit is a conic
union a tangent kernel line, as promised.
\end{example}

\subsection{}\label{non1PS}
Having dealt with the 1-PS case in the previous sections, we may now
assume that 
$$\alpha(t)=\begin{pmatrix}
1 & 0 & 0\\
q(t) & t^b & 0\\
r(t) & s(t)t^b & t^c
\end{pmatrix}$$
with the conditions listed in Proposition~\ref{faber}, and further
{\em such that $q,r$, and $s$ do not all vanish identically.\/} As
four of the five types of components listed in \S\ref{statement} have
been identified, we are left with the task of showing that the only
remaining components of the PNC to which such germs may lead are the
ones `of type~V'. This will take the rest of the section.

The key to the argument will be a further restriction on the germs we
need to consider. We are going to argue that the curve has a rank-2
limit unless $\alpha(t)$ and certain formal branches of the curve are
closely related. 

Work in affine coordinates $(y,z)=(1:y:z)$. If $\mathcal C$ has
multiplicity $m$ at $p=(0,0)$, then we can write the generator $F$ for
the ideal of $\mathcal C$ as a product of formal branches
(cf.~\S\ref{formalbranches})
$$F=f_1\cdot\cdots\cdot f_m$$
where each $f_i$ is expressed as a power series with fractional
exponents. Among these branches, we will especially focus on the ones
that are tangent to the line $z=0$, which may be written explicitly as
$$z=f(y)=\sum_{i\ge 0} \gamma_{\lambda_i} y^{\lambda_i}$$
with $\lambda_i\in \Qbb$, $1<\lambda_0<\lambda_1<\dots$, and
$\gamma_{\lambda_i}\ne 0$.

{\bf Notation.\/} For $C\in \Qbb$, we will denote by $f_{(C)}(y)$ the
finite sum (`truncation')
$$f_{(C)}(y)=\sum_{\lambda_i<C} \gamma_{\lambda_i}
y^{\lambda_i}\quad.$$
For $c\in \Zbb$, we will also write
$\underline{g(t)}$ for the truncation of $g(t)$ to $t^c$, so that
$\underline{f(t^a)}=f_{(C)}(t^a)$ when $C=\frac ca$.
Note that for all $b>a$ the
truncation $\underline{f'(t^a)t^b}$ is determined by $b$ and
$\underline{f(t^a)}$ (and hence by $f_{(C)}(y)$ and $a$, $b$).

\begin{prop}\label{standardform}
Let $\alpha(t)$ be as above, and assume that $\lim_{t\to 0}\mathcal
C\circ\alpha(t)$ is not a rank-2 limit. Then $\mathcal C$ has a formal
branch $z=f(y)$, tangent to $z=0$, such that $\alpha$ is equivalent to
a germ 
$$\begin{pmatrix}
1 & 0 & 0     \\
t^a & t^b & 0 \\
\underline{f(t^a)} & \underline{f'(t^a) t^b} & t^c\end{pmatrix}
\cdot\begin{pmatrix}
1 & 0 & 0 \\
0 & u & 0 \\
r & s & v 
\end{pmatrix}\quad,$$
with $a<b<c$ positive integers, $u,r,s,v\in \Cbb$, and $uv\ne 0$.

Further, it is necessary that $\frac ca\le \lambda_0+2(\frac ba-1)$.
\end{prop}

The proof of this key reduction requires the study of several distinct
cases. We will first show that under the hypothesis that $\lim_{t\to
  0} \mathcal C\circ\alpha(t)$ is not a rank-2 limit we may assume
that $q(t)\ne 0$, and this will allow us to replace it with a power of
$t$; next, we will deal with the $b=c$ case; and finally we will see
that if $b<c$ and $\alpha(t)$ is not in the stated form, then the
limit of {\em every\/} branch of $\mathcal C$ is a
$(0:0:1)$-star. This will imply that the limit of $\mathcal C$ is a
kernel star in this case, proving the assertion by Lemma~\ref{rank2lemma}.

\subsection{}
The first remark is that, under the assumptions that $q$, $r$, and $s$
do not vanish, we may in fact assume that $q(t)$ is not zero.

\begin{lemma}\label{qnot0}
If $\alpha(t)$ is as above, and $q\equiv 0$, then $\lim_{t\to 0}\mathcal
C\circ\alpha(t)$ is a rank-2 limit.
\end{lemma}

\begin{proof}
This is a case-by-case analysis. Assume $q\equiv 0$; thus $r$ and $s$
are not both zero. For $F$ the generator of the ideal of $\mathcal C$,
consider the fate of an individual monomial $x^Ay^Bz^C$ under
$\alpha(t)$:
$$m_{ABC}=x^A y^B (r(t) x+s(t) t^b y+t^c z)^C t^{bB}$$

If $r\equiv 0$ and $s\not\equiv 0$, note that $v(s)\le \deg s<c-b.$
Therefore, the dominating term in $m_{ABC}$ is
$$x^A y^{B+C} t^{bB+(b+v(s))C}\quad.$$
Since $v(s)>0$, 
the weights with which a fixed limit monomial $x^A y^{B+C}$ arises
are mutually distinct, hence the limit monomial with minimum weight
cannot be cancelled. Thus $\lim_{t\to 0}F\circ\alpha(t)$ is the
sum of all the limit monomials $x^A y^{B+C}$ with minimum weight. Thus
$\lim_{t\to 0}\mathcal C\circ\alpha(t)$ is a kernel star, and hence a
rank-2 limit by Lemma~\ref{rank2lemma}.

If  $r\not\equiv 0$ and $s\not\equiv 0$, but $v(r)> b+v(s)$, the same
discussion applies, with the same conclusion.

If  $r\not\equiv 0$, and $s\equiv 0$ or $v(r)< b+v(s)(<c)$, then the
dominating term in $m_{ABC}$ is
$$x^{A+C} y^B t^{bB+v(r)C}\quad;$$
as $v(r)>0$ these limit monomials again have different weights, so
$\lim_{t\to 0}\mathcal C \circ\alpha(t)$ is again a kernel star.

Finally, assume $r\not\equiv 0$, $s\not\equiv 0$, and $v(r)=b+v(s)$
(in particular, $v(r)\ne b$). The dominating terms are then those
$$x^A y^B (r_0 x+s_0 y)^C=r_0^C x^{A+C} y^B+\dots+s_0^C x^A y^{B+C}$$
with minimal $bB+v(r) C$, where $r_0$, $s_0$ are the leading
coefficients in $r(t)$, $s(t)$. These terms cannot all cancel: as
$b\ne v(r)$, there must be exactly one term with maximum $B+C$, and
the corresponding term $x^A y^{B+C}$ cannot be cancelled by other
terms. As the limit is again a kernel star, hence a rank-2 limit, the
assertion is proved.
\end{proof}

\subsection{}
By Lemma~\ref{qnot0} we may now assume that $q(t)\ne 0$. By
Lemma~\ref{faberqnot0} we may then replace $\alpha(t)$ with an
equivalent germ
$$\begin{pmatrix}
1 & 0 & 0     \\
t^a & t^b & 0 \\
r_1(t) & s_1(t)t^b & t^c\end{pmatrix}
\cdot\begin{pmatrix}
1 & 0 & 0 \\
0 & u & 0 \\
r & s & v 
\end{pmatrix}$$
with $a<b\le c$, $r_1(t)$, $s_1(t)$ polynomials, and an invertible
constant factor on the right. This constant factor is the factor
appearing in the statement of Proposition~\ref{standardform}. The
limit of any curve under this germ is a rank-2 limit if and only if
the limit by
$$\begin{pmatrix}
1 & 0 & 0\\
t^a & t^b & 0\\
r_1(t) & s_1(t)t^b & t^c
\end{pmatrix}$$
is a rank-2 limit, so we may ignore the constant matrix on the right
in the rest of the proof of Proposition~\ref{standardform}. 

Renaming $r_1(t)$, $s_1(t)$ by $r(t)$, $s(t)$ respectively, we are
reduced to studying germs
$$\alpha(t)=\begin{pmatrix}
1 & 0 & 0     \\
t^a & t^b & 0 \\
r(t) & s(t)t^b & t^c\end{pmatrix}$$
with $a<b\le c$ positive integers and $r(t)$, $s(t)$ polynomials of 
degree $<c$, $<(c-b)$ respectively and vanishing at $t=0$.

In order to complete the proof of Proposition~\ref{standardform}, we
have to show that if $\lim_{t\to 0}\mathcal C\circ\alpha(t)$ is not a
rank-2 limit then $b<c$ and $r(t)$, $s(t)$ are as stated.

\subsection{}\label{bequalc}
We first deal with the case $b=c$.
\begin{lemma}\label{b=c}
Let $\alpha(t)$ be as above. If $b=c$, then $\lim_{t\to 0} 
\mathcal C\circ\alpha(t)$ is a rank-2 limit.
\end{lemma}

\begin{proof}
If $b=c$, then $s=0$ necessarily:
$$\alpha(t)=\begin{pmatrix}
1 & 0 & 0 \\
t^a & t^b & 0 \\
r(t) & 0 & t^b
\end{pmatrix}\quad,$$
and further $a<v(r)$ (by Lemma~\ref{faberqnot0}). 

Decompose $F(1:y:z)$ in $\Cbb[[y,z]]$: $F(1:y:z)=G(y,z)\cdot H(y,z)$,
where $G(y,z)$ collects the branches that are {\em not\/} tangent to
$z=0$. Writing $G(y,z)$ as a sum of homogeneous terms in $y,z$:
$$G(y,z)=G_{m'}(y,z)+ \text{higher order terms}$$
with $G_{m'}(1,0)\ne 0$, and applying $\alpha(t)$ gives
$$G_{m'}(t^a x+ t^by,r(t) x+ t^bz) +\text{higher order terms}\quad.$$
As $a<v(r)$ and $a<b$, the dominant term in this expression is
$t^{m'a}x^{m'}$: that is, the limit of these branches is supported on
the kernel line $x=0$.

The (formal) branches collected in $H(y,z)$ are tangent to $z=0$, and
we can write such branches as power series with
fractional coefficients (cf.~\S\ref{formalbranches}):
$$z=f(y)=\sum_{i\ge 0} \gamma_{\lambda_i} y^{\lambda_i}$$
with $\lambda_i\in \Qbb$, $1<\lambda_0<\lambda_1<\dots$, and
$\gamma_{\lambda_i}\ne 0$. We will be done if we show that the
limit of such a branch (in the sense of Definition~\ref{branchlimit})
is given by an equation in $x$ and $z$: the limit $\lim_{t\to 0}
\mathcal C\circ\alpha(t)$ will then (cf.~Lemma~\ref{limitlemma}) be a
$(0:1:0)$-star, hence a rank-2 limit by Lemma~\ref{rank2lemma}.

The affine equation of the limit of $z=f(y)$ is given by the dominant
terms in 
$$r(t)+t^b z=f(t^a)+f'(t^a)t^b y+\cdots$$
we observe that $y$ appears on the right-hand-side with weight larger
than $b$ (as $\lambda_0>1$). On the other hand, $z$ only appears on
the left-hand-side, so it cannot be cancelled by other parts in the
expression. It follows that the weight of the dominant terms is $\le
b$, and in particular that $y$ does not appear in these dominant
terms. This shows that the equation of the limit does not depend on
$y$, and we are done. 
\end{proof}

\subsection{}\label{stars}
Next, assume that $\alpha$ is parametrized by
$$\alpha(t)=\begin{pmatrix}
1 & 0 & 0\\
t^a & t^b & 0\\
r(t) & s(t)t^b & t^c
\end{pmatrix}$$
with the usual conditions on $r(t)$ and $s(t)$, and further $b<c$.

We want to study the limits of individual branches of $\mathcal C$
under such a germ. We first deal with branches that are not tangent to
$z=0$:
\begin{lemma}\label{otherbranches}
Under these assumptions on $\alpha$, the limits of branches that
are not tangent to the line $z=0$ are necessarily $(0:0:1)$-stars.
Further, if $a<v(r)$ then the limit of such branches is the kernel
line $x=0$.
\end{lemma}

\begin{proof}
Formal branches that are not tangent to the line $y=0$ may be written
(cf.~\S\ref{formalbranches})
$$z=f(y)=\sum_{i\ge 0} \gamma_{\lambda_i} y^{\lambda_i}$$
with $\lambda_i\in \Qbb$, $1\le \lambda_0<\lambda_1<\dots$, and
$\gamma_{\lambda_i}\ne 0$, and have limit along $\alpha(t)$ given by
the dominant terms in
$$r(t)+s(t)t^by+t^cz=f(t^a)+f'(t^a)t^b y+\dots\quad.$$
Branches that are {\em not\/} tangent to $z=0$ have $\lambda_0=1$,
hence $f'(y)=\gamma_1+\dots$ with $\gamma_1\ne 0$. Hence, the term
$f'(t^a)t^b y$ on the right has weight $t^b$, and is not cancelled by
other terms in the expression (since $v(s)>0$). This implies that the
dominant weight is $\le b<c$, and in particular that the equation of
the limit does not involve~$z$. Hence the limit is a $(0:0:1)$ star,
as needed.

If $a<v(r)$ and $\lambda_0=1$, then the dominant weight is $a<b$,
hence the equation of the limit does not involve $y$ either, so the
limit is a kernel line, as claimed.

Analogous arguments can be used to treat formal branches that are
tangent to $y=0$.
\end{proof}

\subsection{}\label{tangentbranch}
Next, consider a formal branch that is tangent to $z=0$:
$$z=f(y)=\sum_{i\ge 0} \gamma_{\lambda_i} y^{\lambda_i}$$
with $1<\lambda_0<\lambda_1<\dots$

\begin{lemma}\label{tangentbranches}
Under the same assumptions on $\alpha$ as in
Lemma~\ref{otherbranches}, the limit of $z=f(y)$ by $\alpha$ is a
$(0:0:1)$-star unless
\begin{itemize}
\item $r(t)\equiv f(t^a)\pmod{t^c}$;
\item $s(t)\equiv f'(t^a)\pmod{t^{c-b}}$.
\end{itemize}
\end{lemma}

\begin{proof}
The limit of the branch is given by the dominant terms in
$$r(t)+s(t)t^by+t^cz=f(t^a)+f'(t^a)t^by+\dots$$
If $r(t)\not\equiv f(t^a)\pmod{t^c}$, then the weight of the branch
is necessarily $<c$, so the ideal of the limit is generated by a
polynomial in $x$ and $y$, as needed. 
The same reasoning applies if $s(t)\not\equiv f'(t^a)\pmod{t^{c-b}}$.
\end{proof}

\subsection{}\label{stand}
Proposition~\ref{standardform} is now essentially proved. We tie up
here the loose ends of the argument.

\begin{proof} Assume 
$$\alpha(t)=\begin{pmatrix}
1 & 0 & 0\\
q(t) & t^b & 0\\
r(t) & s(t)t^b & t^c
\end{pmatrix}$$
with the conditions listed in Proposition~\ref{faber}, and such that
$q,r$, and $s$ do not all vanish identically, and assume $\lim_{t\to
  0} \mathcal C\circ\alpha(t)$ is not a rank-2 limit. By
Lemma~\ref{qnot0} we may assume $q(t)\ne 0$; hence by
Lemma~\ref{faberqnot0} $\alpha(t)$ is equivalent to a germ
$$\begin{pmatrix}
1 & 0 & 0     \\
t^a & t^b & 0 \\
r_1(t) & s_1(t)t^b & t^c\end{pmatrix}
\cdot\begin{pmatrix}
1 & 0 & 0 \\
0 & u & 0 \\
r & s & v 
\end{pmatrix}\quad,$$
with $a<b\le c$, $r_1(t)$, $s_1(t)$ polynomials, and $u,r,s,v\in \Cbb$
with $uv\ne 0$. The limit along this germ is not a rank-2 limit if and
only if the limit along
$$\begin{pmatrix}
1 & 0 & 0     \\
t^a & t^b & 0 \\
r_1(t) & s_1(t)t^b & t^c\end{pmatrix}$$
is not a rank-2 limit. Assuming this is the case, necessarily $b<c$ by
Lemma~\ref{b=c}. Further, by Lemma~\ref{otherbranches} the limits of
all branches that are not tangent to $z=0$ are $(0:0:1)$-stars, hence
rank-2 limits (by Lemma~\ref{rank2lemma}); the same holds for all
formal branches $z=f(y)$ tangent to $z=0$ unless
$r_1(t)=\underline{f(t^a)}$ and $s_1(t)t^b=\underline{f'(t^a)t^b}$, by
Lemma~\ref{tangentbranches}. Hence, if $\lim_{t\to 0}\mathcal C
\circ\alpha(t)$ is not a rank-2 limit then $\alpha(t)$ must be
equivalent to one in the form given in the statement of the
proposition, for some formal branch $z=f(y)$ tangent to $z=0$.

Finally, to see that the stated condition on $\frac ca$ must hold,
look again at the limit of the formal branch $z=f(y)$, that is, the
dominant term in
$$r(t)+s(t) t^b y+t^c z=f(t^a)+f'(t^a)t^b y+\frac{f''(t^a)t^{2b}y^2}2
+\dots:$$
the dominant weight will be less than $c$ (causing the limit to be a
$(0:0:1)$-star) if $c>2b+v(f''(t^a))=2b+a(\lambda_0-2)$. The stated
condition follows at once.\end{proof}

The effect of the constant factor on the right in the germ appearing
in the statement of Proposition~\ref{standardform} is simply to
translate the limit (by an invertible transformation fixing the flag
consisting of the line $x=0$ and the point $(0:0:1)$). Hence, for the
remaining considerations in this section we may and will ignore this
factor.

Also, we will replace $t$ by $t^{1/d}$ in the germ obtained in
Proposition~\ref{standardform} to ensure that the exponents appearing
in its expression are relatively prime; the resulting germ determines
the same component of the PNC. 

\subsection{}
The next reduction concerns the possible triples $a<b<c$ determining
limits contributing to components of the PNC. This is best expressed
in terms of $B=\frac ba$ and $C=\frac ca$. Let
$$z=f(y)=\sum_{i\ge 0} \gamma_{\lambda_i} y^{\lambda_i}$$
with $\lambda_i\in \Qbb$, $1<\lambda_0<\lambda_1<\dots$, and
$\gamma_{\lambda_i}\ne 0$, be a formal branch tangent to $z=0$.
Every choice of such a branch and of a rational number $C=\frac ca>1$
determines a truncation 
$$f_{(C)}(y)=\sum_{\lambda_i<C} \gamma_{\lambda_i} y^{\lambda_i}\quad.$$
With this notation, the truncation $\underline{f(t^a)}$ equals
$f_{(C)}(t^a)$.

The choice of a rational number $B=\frac ba$ satisfying $1<B<C$
determines now a germ as prescribed by Proposition~\ref{standardform}:
$$\alpha(t)=\begin{pmatrix}
1 & 0 & 0     \\
t^a & t^b & 0 \\
\underline{f(t^a)} & \underline{f'(t^a) t^b} & t^c\end{pmatrix}$$
(choosing the smallest positive integer $a$ for which the entries of 
this matrix have integer exponents). Observe that the truncation
$\underline{f(t^a)}=f_{(C)}(t^a)$ is identically~0 if and only if
$C\le\lambda_0$. Also observe that $\underline{f'(t^a)t^b}$ is
determined by $f_{(C)}(t^a)$, as it equals the truncation to $t^c$ 
of~${(f_{(C)})}'(t^a)t^b$.

\begin{prop}\label{abc}
If $C\le\lambda_0$ or $B\ne \frac{C-\lambda_0}2+1$, then $\lim_{t\to
  0}\mathcal C\circ\alpha(t)$ is a rank-2 limit.
\end{prop}

We deal with the different cases separately.

\begin{lemma}\label{Clelambda0}
If $C\le\lambda_0$, then $\lim_{t\to 0}\mathcal C\circ\alpha(t)$ is a
$(0:1:0)$-star.
\end{lemma}

\begin{proof}
If $C=\frac ca\le\lambda_0$, then $f_{(C)}(y)=0$, so
$$\alpha(t)=\begin{pmatrix}
1 & 0 & 0     \\
t^a & t^b & 0 \\
0 & 0 & t^c\end{pmatrix}\quad.$$
Collect the branches that are not tangent to $z=0$ into
$\Phi(y,z)\in\Cbb[[y,z]]$, with initial form $\Phi_m(y,z)$. Applying
$\alpha(t)$ to these branches gives 
$$\Phi_m(t^a x+t^b y,t^c z)+\dots$$
with limit a kernel line since $a<c$ and $\Phi_m(1:0)\ne 0$.

As for the branches that are tangent to $z=0$, let $z=f(y)$ be such a
formal branch. The limit along $\alpha(t)$ is given by the dominant
terms in
$$t^c z=f(t^a)+f'(t^a)t^by+\dots$$
All terms on the right except the first one have weight larger than 
$a\lambda_0\ge a C=c$, hence the dominant term 
does not involve $y$, concluding the proof.
\end{proof}

By Lemma~\ref{rank2lemma}, the limits obtained in
Lemma~\ref{Clelambda0} are rank-2 limits, so the first part of
Proposition~\ref{abc} is proved. As for the second part, if $B<
\frac{C-\lambda_0}2+1$ then $C>\lambda_0+2(B-1)$, and the limit is a
rank-2 limit by the last assertion in Proposition~\ref{standardform}.

For $B\ge \frac{C-\lambda_0}2+1$, the limit of a branch tangent to
$z=0$ depends on whether the branch truncates to $f_{(C)}(y)$ or not.
These cases are studied in the next two lemmas.

\begin{lemma}\label{nottrunc}
Assume $C>\lambda_0$ and $B\ge \frac{C-\lambda_0}2+1$, and let
$z=g(y)$ be a formal branch tangent to $z=0$, such that $g_{(C)}(y)\ne
f_{(C)}(y)$. Then the limit of the branch is supported on a kernel line.
\end{lemma}

\begin{proof}
The limit of the branch is determined by the dominant terms in
$$\underline{f(t^a)}+\underline{f'(t^a)t^b}y+t^cz=g(t^a)+g'(t^a)
t^by+\dots $$
Assume the truncations $g_{(C)}$ and $f_{(C)}$ do not agree. If the
first term at which they disagree has weight lower than
$B+\lambda_0-1$, then the dominant terms in the expansion have weight
lower than the weight of $\underline{f'(t^a)t^b}y$, and it follows
that the limit is supported on $x=0$. So we may assume that 
$$g_{(C)}(y)=f_{(C)}(y)+y^{B+\lambda_0-1}\rho(y)$$
for some $\rho(y)$. We claim that then the terms
$\underline{f'(t^a)t^b}y$ and $g'(t^a)t^by$ agree modulo $t^c$: this
implies that the dominant term is independent of $y$. As the dominant
term is also independent of $z$ (since the truncations $g_{(C)}$ and
$f_{(C)}$ do not agree), the statement will follow from our claim.

In order to prove the claim, observe that
$$(g_{(C)})'(y)=(f_{(C)})'(y)+(B+\lambda_0-1)y^{B+\lambda_0-2}\rho(y)
+y^{B+\lambda_0-1}\rho'(y)\quad:$$
thus, $(g_{(C)})'(y)y^B$ and $(f_{(C)})'(y)y^B$ must agree modulo
$y^{2B+\lambda_0-2}$. Since $B\ge\frac{C-\lambda_0}2+1$, we have
$$2B+\lambda_0-2\ge(C-\lambda_0)+2+\lambda_0-2=C\quad;$$
hence $(g_{(C)})'(t^a)t^b$ and $(f_{(C)})'(t^a)t^b$ must agree modulo
$t^{aC}=t^c$. It follows that $\underline{f'(t^a)t^b}y$ and
$g'(t^a)t^by$ agree modulo $t^c$, and we are done.
\end{proof}

\begin{lemma}\label{dominant}
Assume $C>\lambda_0$ and $B\ge \frac{C-\lambda_0}2+1$, and let
$z=g(y)$ be a formal branch tangent to $z=0$, such that $g_{(C)}(y)=
f_{(C)}(y)$. Denote by $\gamma_C^{(g)}$ the coefficient of $y^C$ in
$g(y)$. 
\begin{itemize}
\item If $B> \frac{C-\lambda_0}2+1$, then the limit of the branch
  $z=g(y)$ by $\alpha(t)$ is the line
$$z=(C-B+1)\gamma_{C-B+1} y+\gamma_C^{(g)}\quad.$$

\item If $B= \frac{C-\lambda_0}2+1$, then the limit of the branch
  $z=g(y)$ by $\alpha(t)$ is the conic
$$z=\frac {\lambda_0(\lambda_0-1)}2\gamma_{\lambda_0}y^2+\frac
{\lambda_0+C}2\gamma_{\frac{\lambda_0+C}2}y+\gamma_C^{(g)}\quad.$$
\end{itemize}
\end{lemma}

\begin{proof}
Rewrite the expansion whose dominant terms give the limit of the
branch as:
$$t^c z=(g(t^a)-\underline{f(t^a)})+(g'(t^a)t^b-\underline{f'(t^a)t^b})y
+\frac{g''(t^a)}2 t^{2b} y^2+\dots$$
The dominant term has weight $c=Ca$ by our choices; if
$B> \frac{C-\lambda_0}2+1$ then the weight of the coefficient of $y^2$
exceeds $c$, so it does not survive the limiting process, and the limit
is a line. If $B= \frac{C-\lambda_0}2+1$, the term in $y^2$ is
dominant, and the limit is a conic. 

The explicit expressions given in the statement are obtained by
reading the coefficients of the dominant terms.
\end{proof}

We can now complete the proof of Proposition~\ref{abc}:
\begin{lemma}\label{completion}
If $B>\frac{C-\lambda_0}2+1$, then the limit $\lim_{t\to 0} \mathcal C
\circ\alpha(t)$ is a rank-2 limit.
\end{lemma}

\begin{proof}
We will show that the limit is necessarily a kernel star, which
gives the statement by Lemma~\ref{rank2lemma}.

As $B>1$, the coefficient $\gamma_{C-B+1}$ is determined by the
truncation $f_{(C)}$, and in particular it is the same for
all formal branches with that truncation. If
$B>\frac{C-\lambda_0}2+1$, by Lemma~\ref{dominant} the branches
contributes a line through the fixed point
$(0:1:(C-B+1)\gamma_{C-B+1})$. We are done if we check that all other
branches contribute a kernel line $x=0$: and this is implied by
Lemma~\ref{otherbranches} for branches that are not tangent to $z=0$
(note $a<v(r)$ for the germs we are considering), and by
Lemma~\ref{nottrunc} for formal branches $z=g(y)$ tangent to $z=0$ but
whose truncation $g_{(C)}$ does not agree with $f_{(C)}$.
\end{proof}

\subsection{}\label{quadritangent}
Finally we are ready to complete the description of the components
given in \S\ref{statement}. By Propositions~\ref{standardform} and
\ref{abc}, germs leading to components of the PNC that have not yet
been accounted for must be (up to a constant translation, and up to
replacing $t$ by $t^{1/d}$) in the form
$$\alpha(t)=\begin{pmatrix}
1 & 0 & 0 \\
t^a & t^b & 0 \\
\underline{f(t^a)} & \underline{f'(t^a)t^b} & t^c
\end{pmatrix}$$
for some branch $z=f(y)=\gamma_{\lambda_0}y^{\lambda_0}+\dots$ of 
$\mathcal C$ tangent to $z=0$ at $p=(0,0)$, and further satisfying
$C>\lambda_0$ and $B=\frac{C-\lambda_0}2+1$ for $B=\frac ba$, $C=\frac
ca$. New components of the PNC will arise depending on the limit
$\lim_{t\to 0}\mathcal C\circ\alpha(t)$, which we now determine.

\begin{lemma}\label{quadconics}
If $C>\lambda_0$ and $B=\frac{C-\lambda_0}2+1$, then the limit
$\lim_{t\to 0}\mathcal C\circ\alpha(t)$ consists of a union of
quadritangent conics, with distinguished tangent equal to the kernel
line $x=0$, and of a multiple of the distinguished tangent line.
\end{lemma}

\begin{proof}
Both $\gamma_{\lambda_0}$ and $\gamma_{\frac{\lambda_0+C}2}$
are determined by the truncation $f_{(C)}$ (since $C>\lambda_0$);
hence the equations of the conics 
$$z=\frac {\lambda_0(\lambda_0-1)}2\gamma_{\lambda_0}y^2+\frac
{\lambda_0+C}2\gamma_{\frac{\lambda_0+C}2}y+\gamma_C$$
contributed (according to Lemma~\ref{dominant}) by different branches
with truncation $f_{(C)}$ may only differ in the $\gamma_C$
coefficient.

It is immediately verified that all such conics are tangent to the
kernel line $x=0$, at the point $(0:0:1)$, and that any two such
conics meet only at the point $(0:0:1)$; thus they are necessarily
quadritangent.

Finally, the branches that do not truncate to $f_{(C)}(y)$ must
contribute kernel lines, by Lemmas~\ref{otherbranches} and \ref{nottrunc}.
\end{proof}

The type of curves arising as the limits described in
Lemma~\ref{quadconics} are studied in \cite{MR2002d:14083},~\S4.1; also
see item~(12) in the classification of curves with small orbit in
\cite{MR2002d:14084},~\S1. The degenerate case in which only one conic
appears does not lead to a component of the projective normal cone, by
the usual dimension considerations. A component is present as
soon as there are two or more conics, that is, as soon as two branches
contribute distinct conics to the limit.

This leads to the description given in \S\ref{statement}. We say that
a rational number $C$ is `characteristic' for $\mathcal C$ (with
respect to $z=0$) if at least two formal branches of $\mathcal C$
(tangent to $z=0$) have the same nonzero truncation, but different
coefficients for $y^C$.

\begin{prop}\label{typeV}
The set of characteristic rationals is finite.

The limit $\lim_{t\to 0}\mathcal C \circ\alpha(t)$ obtained in
Lemma~\ref{quadconics} determines a component of the projective normal
cone precisely when $C$ is characteristic.
\end{prop}

\begin{proof}
If $C\gg 0$, then branches with the same truncation must in fact be
identical, hence they cannot differ at $y^C$, hence $C$ is not
characteristic. Since the set of exponents of any branch is discrete,
the first assertion follows.

The second assertion follows from Lemma~\ref{quadconics}: if
$C>\lambda_0$ and $B=\frac{C-\lambda_0}2+1$, then the limit is a union
of a multiple kernel line and conics with equation
$$z=\frac {\lambda_0(\lambda_0-1)}2\gamma_{\lambda_0}y^2+\frac
{\lambda_0+C}2\gamma_{\frac{\lambda_0+C}2}y+\gamma_C\quad:$$
these conics are different precisely when the coefficients $\gamma_C$
are different, and the statement follows.
\end{proof}

Proposition~\ref{typeV} leads to the procedure giving components of
type~V explained in \S\ref{statement} (also cf.~\cite{MR2001h:14068},
\S2, Fact~5), concluding the set-theoretic description of the
projective normal cone given there.

\subsection{}
Rather than reproducing from \S\ref{statement} the procedure leading
to components of type~V, we offer another explicit example by applying
it to the curve considered in ~\S\ref{7ic}.

\begin{example}\label{7icb}
We obtain the components of type~V due to the singularity at the point
$p=(1:0:0)$ on the curve
$$x^3 z^4-2 x^2 y^3 z^2+x y^6-4 x y^5 z-y^7=0$$
The ideal of the tangent cone at $p$ is $(z^4)$, so that this curve
has four formal branches, all tangent to the line $z=0$. These can be
computed with ease:
$$\left\{\aligned
z&=\phantom{-}y^{3/2}-\phantom{i}y^{7/4}\\
z&=\phantom{-}y^{3/2}+\phantom{i}y^{7/4}\\
z&=-y^{3/2}+iy^{7/4}\\
z&=-y^{3/2}-iy^{7/4}
\endaligned\right.\quad.$$
We find a single characteristic $C=\frac 74$, and two truncations
$y^{3/2}$, $-y^{3/2}$. In both cases $\lambda_0=\frac 32$, so
$B=\frac{\frac 74-\frac 32}2+1=\frac 98$. The lowest integer $a$
clearing all denominators is $8$, and we find the two germs
$$\alpha_1(t)=\begin{pmatrix}
1 & 0 & 0\\
t^8 & t^9 & 0\\
t^{12} & \frac 32 t^{13} & t^{14}
\end{pmatrix}\quad,\quad
\alpha_2(t)=\begin{pmatrix}
1 & 0 & 0\\
t^8 & t^9 & 0\\
-t^{12} & -\frac 32 t^{13} & t^{14}
\end{pmatrix}\quad.$$
As in each case there are two branches (with different
$y^{7/4}$-coefficients), the limits must in both cases be pairs of
quadritangent conics, union a triple kernel line. The limit along
$\alpha_1$ is listed in \S\ref{7ic}; the limit along $\alpha_2$ must
be, according to the formula given above,{\small
$$x^3\left(zx-\frac{\frac 32\cdot \frac 12}2(-1)y^2-\frac{\frac
  32+\frac 74}2\cdot 0\cdot yx-ix^2\right)
\left(zx-\frac{\frac 32\cdot \frac 12}2(-1)y^2-\frac{\frac
  32+\frac 74}2\cdot 0\cdot yx+ix^2\right)$$}
that is (up to a constant factor),
$$x^3\left(64x^2z^2+48xy^2z+9y^4+64x^4\right)\quad,$$
as may also be checked directly.
\end{example}

\subsection{}\label{Ghizz}
In this subsection we briefly describe the contents of Aldo Ghizzetti's
first paper \cite{ghizz}. Following this paper, Ghizzetti turned
to analysis under the mentorship of Guido Fubini and Mauro Picone.
His 50 year career was crowned by the election into the 
{\em Accademia Nazionale dei Lincei\/}. See~\cite{MR1284903}. 

In \cite{ghizz}, Ghizzetti reports the results of his 1930 thesis under
Alessandro Terracini. He determines here the limits of one-parameter
families of `homographic' plane curves (that is, of curves in the same
orbit under the $\PGL(3)$ action).

A main feature of loc.~cit.~is the classification of one-dimensional
systems of `homographies' approaching a degenerate homography. Let
$\Omega_t=(a_{ik}(t))$ be such a system, where the coefficients
$a_{ik}(t)$ are power series in $t$ in a neighborhood of $t=0$; we
assume that $\det(a_{ik}(t))$ vanishes at $t=0$ but is not identically
zero. The $\Omega_t$ are viewed as transformations from a projective
plane $\pi$ to another plane $\pi'$:
$$x_i'=\sum_{k=1}^3 a_{ik}(t)x_k\quad,\qquad(i=1,2,3).$$
A plane curve $C'$ of degree $d$ 
in $\pi'$ with equation $F(x_1',x_2',x_3')=0$ is
given; transforming it by $\Omega_t$ gives a curve $C_t$ in $\pi$ and the
goal is to determine the limiting curve $C_0$. The coordinates $x_k$
may be modified by a system of non-degenerate homographies of $\pi$
to reduce $\Omega_t$ to a simpler form.

When $\Omega_0$ has rank~$2$, it defines a kernel point $S$ in $\pi$
and an image line $s'$ in $\pi'$. When $C'$ does not contain $s'$, it is
easy to see that $C_0$ consists of $d$ lines through $S$, the $d$-tuple
being projectively equivalent to the $d$ points that $C'$ cuts out on $s'$.
Let $s_0$ be the limiting position of the line $s_t$ that results from
transforming $s'$ by $\Omega_t$. 
When $C'$ contains $s'$ with multiplicity $m$, the limiting curve $C_0$
consists of a $(d-m)$-tuple of lines through $S$ (projectively equivalent to
the $(d-m)$-tuple of points that the residual of $C'$ cuts out on $s'$) and
the line $s_0$ with multiplicity $m$. These are the limit curves of type~I.

When $\Omega_0$ has rank $1$, it defines a kernel line $s$ in $\pi$
and an image point $S'$ in $\pi'$. When $C'$ doesn't pass through $S'$,
the limiting curve $C_0$ consists of the line $s$ with multiplicity $d$.
Assume from now on that $C'$ passes through $S'$. For $t\neq0$, let
$S_t$ be the inverse image point under $\Omega_t$ of $S'$ and let $S_0$
be the limiting point. The author considers two cases: $S_0$ is not
contained in $s$ (case~I) or it is (case~II). 
(We reserve the word `cases' for Ghizzetti's classification and continue
to use `types' for ours.) After a change of coordinates
in $\pi$ it may be assumed that $S_t$ is fixed and coincides with $S_0$.
Then $\Omega_t$ induces a linear map $\omega_t$ on the pencils through
$S_0$ and $S'$ and one obtains a limiting transformation $\omega_0$,
which is either non-degenerate (cases~${\rm I}_1$ and~${\rm II}_1$)
or degenerate (cases~${\rm I}_2$ and~${\rm II}_2$).

Case~${\rm I}_1$ yields one-parameter subgroups with two equal weights;
this gives the limit curves of type~III, the fans of 
Proposition~\ref{tgcone}. In case~${\rm II}_1$ the kernel line $s$
contains $S_0$ and the limit curves become kernel stars. Ghizzetti
explicitly shows that they are rank-2 limits.
\newcommand{\os}{\overline{s}}

In cases~${\rm I}_2$ and~${\rm II}_2$, the degenerate map $\omega_0$
defines a kernel line $\os$ through $S_0$ and an image line $\os'$
through $S'$. Let $\os_0$ be the limiting position of the line $\os_t$
that results from transforming $\os'$ by $\omega_t$. The author
distinguishes two subcases: $\os$ and $\os_0$ are distinct
(cases~${\rm I}_{21}$ and~${\rm II}_{21}$) or they coincide
(cases~${\rm I}_{22}$ and~${\rm II}_{22}$). Moreover, in case~${\rm II}_{22}$
it is necessary to distinguish whether $s$ differs from $\os=\os_0$
(case~${\rm II}_{22}'$) or coincides with it (case~${\rm II}_{22}''$).

Case~${\rm I}_{21}$ yields one-parameter subgroups with three distinct
weights: in the notation of \S\ref{1PS} (with $b<c$), $s$ is $x=0$,
$\os$ is $y=0$, and $\os_0$ is $z=0$. This gives the limit curves
of types~II and~IV, cf.~\S\S\ref{Newton} and \ref{typeII}.

Each of the cases~${\rm I}_{22}$, ${\rm II}_{21}$, and~${\rm II}_{22}'$
leads to limit curves that are kernel stars; again, Ghizzetti shows
explicitly that they are rank-2 limits.

The case~${\rm II}_{22}''$ remains: $s=\os=\os_0$. Assume that $C'$
is tangent to $\os'$ in $S'$; if not, $C_0$ consists of the line $s$ with 
multiplicity $d$. After a change of coordinates in $\pi$ and a change
of parameter (similar to the one in Lemma~\ref{faberqnot0}), the matrix
of $\Omega_t$ has three zero entries, one entry equal to $1$, and one
entry a positive power of $t$; each of the remaining entries vanishes at
$t=0$ but is not identically zero. (The matrix is {\em not\/} in triangular
form.) The vanishing orders of two of the entries ($m$ and $m+n$ in 
Ghizzetti's notation) arose naturally in the classification leading
to the present case. The other three vanishing orders are called $m+p$,
$q$, and $r$, and it is necessary to analyze the various possibilities
for the triple $p$, $q$, and $r$. (In the last section of loc.~cit.~the author
remarks that $m$ and $n$ are the degree and class of the curve of image
points under $\Omega_t$ of a fixed point in $\pi$ and that $p$, $q$, and $r$
are similarly related to the curve of inverse image points of a fixed
point in $\pi'$.)

Ghizzetti distinguishes five cases in his analysis of the triple
$p$, $q$, and $r$. After simplifying coordinate changes, the $\Omega_t$
are applied to the branches of $C'$ at $S'$ that are tangent to $\os'$.
The first four cases lead to limit curves that are kernel stars.
In the study of the fifth case, it is necessary to distinguish five subcases.
The first three lead to limit curves that are kernel stars. The fourth case
is not worked out in detail. The limit curve consists of lines; in general, 
not all of these lines belong to the same pencil, according to Ghizzetti,
but in fact they form a kernel star (cf.~Lemma~\ref{completion}).
The fifth and final case leads to a union of quadritangent conics and
a multiple of the kernel line, as in Lemma~\ref{quadconics}.
These are the limit curves of type~V; with them, Ghizzetti concludes
his analysis of the possible limit curves arising from a system
of non-degenerate homographies approaching a degenerate homography.

One of the consequences of Ghizzetti's work is that the irreducible
components of limit curves are very special from the projective
standpoint; for example, they are necessarily isomorphic to their own
dual. It is clear that such curves have `small' linear orbit;
Ghizzetti's characterization was extended to curves with small orbit
in projective spaces of arbitrary dimension in one of Ciro Ciliberto's
first papers, \cite{MR58:22057}. Plane curves with small linear orbits
are classified in \cite{MR2002d:14084}.

It will be obvious from the above that Ghizzetti's approach and ours
are quite similar. Let us then conclude this section by indicating
some of the differences. From a technical viewpoint,
Proposition~\ref{faber} and its refinement Lemma~\ref{faberqnot0} lead
us to the key reduction of  Proposition~\ref{standardform}; this result
allows us to restrict our attention to germs that are essentially
determined by a branch of $\mathcal C$. Although this reduction is
perhaps not as geometrically meaningful as Ghizzetti's approach, it
appears to lead to a considerable simplification. Also, the fact that
the limit curves necessarily have infinite stabilizer plays a less
prominent role in~\cite{ghizz}. Most importantly, Ghizzetti's goal was
essentially the set-theoretic description of the boundary components
of the linear orbit of a curve, while our enumerative applications in
\cite{MR2001h:14068} require the more refined information carried by
the projective normal cone dominating the boundary. This forces us to
be more explicit concerning equivalence of germs, and leads us to a
rather different classification than the one considered by
Ghizzetti. In fact, the set-theoretic description alone of (even) the
projective normal cone does not suffice for our broader goals, so that
we need to refine our analysis considerably in order to determine the
projective normal cone {\em as a cycle,\/} in the next section. For
this, the notion of equivalence introduced in
Definition~\ref{equivgerms} will play a crucial role,
cf.~Proposition~\ref{equivgermsagain} and Lemma~\ref{strongeq}.

In any case, Ghizzetti's contribution remains outstanding for its
technical prowess, and it is an excellent example of concreteness and
concision in the exposition of very challenging material.
The fact that this is his first paper makes our admiration of it
only stronger.


\section{The projective normal cone as a cycle}\label{multi}

\subsection{}
As mentioned in \S\ref{intro} and \S\ref{prelim}, the enumerative
computations in \cite{MR2001h:14068} require the knowledge of the {\em
cycle\/} supported on the projective normal cone; that is, we need to
compute the {\em multiplicities\/} with which the components
identified in \S\ref{setth} appear in the PNC. This is what we do in
this section.

Our general strategy will be the following. By normalizing the graph
of the basic rational map $\Pbb^8 \dashrightarrow \Pbb^N$, we will
distinguish different `ways' in which a component may arise, and
compute a contribution to the multiplicity due to each way. 
This contribution will be obtained by carefully evaluating different
ingredients: the order of contact of certain germs with the base
scheme of the rational map, the number of components in the
normalization dominating a given component of the PNC, and the degree
of the restriction of the normalization map to these components.

\subsection{}\label{multstat}
Here is the result. The `multiplicities' in the following list are
contributions to the multiplicity of each individual component from
the possibly different ways to obtain it.
This list should be compared with~\cite{MR2001h:14068}, \S2, Facts~1
through~5.

\begin{itemize}

\item{\em Type I.\/} The multiplicity of the component determined by a
  line $\ell\subset \mathcal C$ equals the multiplicity of $\ell$ in
  $\mathcal C$.

\item{\em Type II.\/} The multiplicity of the component determined by
  a nonlinear component $\mathcal C'$ of $\mathcal C$ equals $2m$,
  where $m$ is the multiplicity of $\mathcal C'$ in $\mathcal C$.

\item{\em Type III.\/} The multiplicity of the component determined by
  a singular point $p$ of $\mathcal C$ such that the tangent cone
  $\lambda$ to $\mathcal C$ at $p$ is supported on three or more lines
  equals $mA$, where $m$ is the multiplicity of $\mathcal C$ at $p$
  and $A$ equals the number of automorphisms of $\lambda$ as a tuple
  in the pencil of lines through $p$.

\item{\em Type IV.\/} The multiplicity of the component determined by
  one side of a Newton polygon for $\mathcal C$, with vertices
  $(j_0,k_0)$, $(j_1,k_1)$ (where $j_0<j_1$) and limit
$$x^{\overline q} y^r z^q \prod_{j=1}^S\left(y^c+\rho_j
x^{c-b} z^b\right)\quad,$$
(with $b$ and $c$ relatively prime) equals
$$\frac{j_1k_0-j_0k_1}SA\quad,$$
where $A$ is the number of automorphisms $\Abb^1 \to \Abb^1$, $\rho
\mapsto u\rho$ (with $u$ a root of unity) preserving the $S$-tuple
$\{\rho_1,\dots, \rho_S\}$.\footnote{Note: the number $A$ given here
  is denoted $A/\delta$ in \cite{MR2001h:14068}.}

\item{\em Type V.\/} The multiplicity of the component corresponding
  to the choice of a characteristic $C$ and a truncation $f_{(C)}(y)$
  at a point $p$, with limit
$$x^{d-2S}\prod_{i=1}^S\left(zx-\frac {\lambda_0(\lambda_0-1)}2
\gamma_{\lambda_0}y^2 -\frac{\lambda_0+C}2
\gamma_{\frac{\lambda_0+C}2}yx-\gamma_C^{(i)}x^2\right)$$
is $\ell WA$, where:
\begin{itemize}
\item $\ell$ is the least positive integer $\mu$ such that
  $f_{(C)}(y^\mu)$ has integer exponents.
\item $W$ is defined as follows. For each formal branch $\beta$ of
  $\mathcal C$ at $p$, let $v_\beta$ be the first exponent at which
  $\beta$ and $f_{(C)}(y)$ differ, and let $w_\beta$ be the minimum of
  $C$ and $v_\beta$. Then $W$ is the sum $\sum w_\beta$.
\item $A$ is twice the number of automorphisms $\gamma\to u\gamma+v$
preserving the $S$-tuple $\{\gamma_C^{(1)},\dots,\gamma_C^{(S)}\}$.
\end{itemize}

\end{itemize}

Concerning the `different ways' in which a component may be obtained
(each producing a multiplicity computed by the above recipe), no
subtleties are involved for components of type I, II, or III: there is
only one contribution for each of the specified data---that is,
exactly one contribution of type~I from each line contained in
$\mathcal C$, one contribution of type~II from each nonlinear
component of $\mathcal C$, and one of type~III from each singular
point of $\mathcal C$ at which the tangent cone is supported on three
or more distinct lines.

As usual, the situation is a little more complex for components of
type~IV and~V.

{\em Components of type~IV\/} correspond to sides of Newton polygons;
one polygon is obtained for each line in the tangent cone at a fixed
singular point $p$ of $\mathcal C$, and each of these polygons
provides a set of sides (with slope strictly between $-1$ and $0$).
Exactly one contribution has to be counted for each side obtained in
this fashion. Note that sides of different Newton polygons may lead to
the same limits, hence to the same component of the PNC.

\begin{example}
The curve 
$$(y-z^3)(z-y^3)=0$$
has a node at the origin $p$: $(y,z)=(0,0)$; the two lines in the
tangent cone both give the Newton polygon
$$\includegraphics[scale=1]{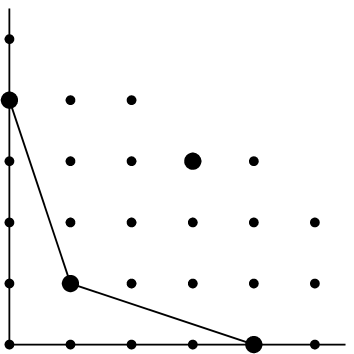}$$
each with one side with slope $-1/3$. The corresponding limits,
obtained (as prescribed in \S\ref{Newton})
via the germs
$$\begin{pmatrix}
1 & 0 & 0 \\
0 & t & 0 \\
0 & 0 & t^3
\end{pmatrix}
\quad,\quad
\begin{pmatrix}
1 & 0 & 0 \\
0 & t^3 & 0 \\
0 & 0 & t
\end{pmatrix}$$
are respectively
$$y(z-y^3)\quad,\quad z(y-z^3)\quad.$$
These limits belong to the same $\PGL(3)$ orbit; thus the two sides
determine the same component of the PNC. According to the result given
above, the multiplicity of this component receives a contribution of $4$
from each of the sides, so the component appears with multiplicity~8
in the PNC.
\end{example}

{\em Components of type~V\/} are determined by a choice of a singular
point $p$ of $\mathcal C$, a line $L$ in the tangent cone to $\mathcal
C$ at $p$, a characteristic $C$ and a truncation $f_{(C)}(y)$ of a
formal branch of $\mathcal C$
tangent to $L$. Recall that this data determines a triple of positive
integers $a<b<c$ with $C=c/a$: the number $C$ and the truncation
$f_{(C)}(y)$ determine $B$ as in the beginning of
\S\ref{quadritangent}, and $a$ is the smallest integer clearing
denominators of all exponents in the
corresponding germ $\alpha(t)$. Again, different choices may
lead to the same component of the PNC, and we have to specify when
choices should be counted as giving separate contributions. Of course
different points $p$ or different lines in the tangent cone at $p$
give separate contributions; the question is when two sets of data
$(C,f_{(C)}(y))$ for the same point, with respect to the same tangent
line, should be counted separately.

To state the result, we say that $(C,f_{(C)}(y))$, $(C',g_{(C')}(y))$
(or the truncations $f_{(C)}$, $g_{(C')}$ for short) are {\em
  sibling\/} data if $C=C'$ and
$f_{(C)}(t^a)=g_{(C)}((\xi t)^a)$ for an $a$-th root $\xi$ of $1$.

\begin{example}\label{sibl} 
If $z=f(y)$, $z=g(y)$ are formal branches belonging to
  the same irreducible branch of $\mathcal C$ at $p$, then the
  corresponding truncations $f_{(C)}(y)$, $g_{(C)}(y)$ are siblings for
  all $C$.

Indeed, if the branch has multiplicity $m$ at $p$ then
$f(\tau^m)=\varphi(\tau)$ and $g(\tau^m)=\psi(\tau)$, with
$\psi(\tau)=\varphi(\zeta\tau)$
for an $m$-th root $\zeta$ of $1$ (\cite{MR1836037}, \S7.10). That is,
$$\text{if}\quad 
f(y) =\sum \gamma_{\lambda_i} y^{\lambda_i}\quad,\quad\text{then}\quad
g(y) =\sum \zeta^{m\lambda_i}\gamma_{\lambda_i} y^{\lambda_i}$$
for $\zeta$ an $m$-th root of $1$. Note that the positive integer $a$
determined by $C$ is such that $a\lambda_i$ is integer for all the
exponents $\lambda_i$ lower than $C$.

Now let $\rho$ be an $(am)$-th root of $1$ such that $\rho^a=\zeta$,
and set $\xi=\rho^m$; since the exponents $a\lambda_i$ in the
truncations are integers, as well as all exponents $m\lambda_i$, we
have
$$\zeta^{m\lambda_i}=\rho^{ma\lambda_i}=\xi^{a\lambda_i}$$
for all exponents $\lambda_i<C$, and this shows that the truncations
are siblings.
\end{example}

With this notion, we can state precisely when two truncations
$(C,f_{(C)}(y))$ at the same point, with respect to the same tangent
line, yield separate contributions: they do if and only if they are
{\em not\/} siblings.\vskip 6pt

The proof that the formulas presented above hold occupies the rest of 
this paper. Rather direct arguments can be given for the `global'
components, of type~I and type~II, cf.~\S\S\ref{multtypeI} and
\ref{multtypeII}. The `local' types III, IV, and~V require some
preliminaries, covered in the next several sections.

\subsection{}\label{beginprel}
The main character in the story will be the {\em normalization\/}
$\nG$ of the closure $\Til\Pbb^8$ of the graph of the basic rational
map $\Pbb^8 \dashrightarrow \Pbb^N$ introduced in \S\ref{prelim}. We
denote by
$$\nor: \nG \rightarrow \Til\Pbb^8$$ 
the normalization map, and by $\onu$ the composition $\nG \to
\Til\Pbb^8 \to \Pbb^8$.

In \S\ref{prelim}, Lemma~\ref{excdivsupp} we have realized the PNC as
a subset of $\Til\Pbb^8\subset \Pbb^8\times \Pbb^N$. Recall that the
PNC is in fact the exceptional divisor $E$ in $\Til\Pbb^8$, where the
latter is viewed as the blow-up of $\Pbb^8$ along the base scheme $S$
of the basic rational map. If $F\in\Cbb[x,y,z]$ generates the ideal of
$\mathcal C$ in $\Pbb^2$, then the ideal of $S$ in $\Pbb^8$ is
generated by all
$$F(\varphi(x_0,y_0,z_0))$$
viewed as polynomials in $\varphi\in\Pbb^8$, as $(x_0,y_0,z_0)$
ranges over $\Pbb^2$.
We denote by $E_i$ the supports of the
components of $E$, and by $m_i$ the multiplicity of $E_i$ in $E$.

We also denote by $\oE$ the Cartier divisor
$\nor^{-1}(E)=\onu^{-1}(S)$ in $\nG$, and by 
$$\oE_{i1},\dots,\oE_{ir_i}$$
the supports of the components of $\oE$ lying above a given component
$E_i$ of $E$. Finally, we let $m_{ij}$ be the multiplicity of
$\oE_{ij}$ in $\oE$. That is:
$$[E]=\sum m_i [E_i]\quad,\quad [\oE]=\sum m_{ij} [\oE_{ij}]\quad.$$

\begin{prop}\label{multcount}
We have
$$m_i=\sum_{j=1}^{r_i} e_{ij} m_{ij}$$
where $e_{ij}$ is the degree of $\nor|_{\oE_{ij}}:\oE_{ij} \to E_i$.
\end{prop}

\begin{proof}
This follows from the projection formula and $(\nor|_{\oE_{ij}})_*
[\oE_{ij}]=e_{ij}[E_i]$.
\end{proof}

\subsection{}
In order to apply Proposition~\ref{multcount}, we have to develop
tools to evaluate the multiplicities $m_{ij}$ and the degrees
$e_{ij}$. For `local' components, we will obtain this information by
describing $\oE_{ij}$ in terms of lifts of certain germs from $\Pbb^8$.
Until the end of \S\ref{weights} we focus on components of type
III,~IV, and~V.

Every germ $\alpha(t)$ in $\Pbb^8$, whose general element is
invertible, and such that $\alpha(0)\in S$, lifts to a unique germ in
$\Til\Pbb^8$ centered at a point $(\alpha(0),\mathcal X)$ of the
PNC. The germ lifts to a unique germ in $\nG$, centered at a point of
$\oE$.

\begin{defin}\label{markers}
We denote by $\oalpha$ the center of the lift of $\alpha(t)$ to $\nG$.
We will say that $\alpha(t)$ is a `marker' germ if
\begin{itemize}
\item $\oalpha$ belongs to exactly one component of $\oE$, and the
lift of $\alpha(t)$ is transversal to (the support of) $\oE$ at
$\oalpha$; 
\item $\nG$ is nonsingular at $\oalpha$; 
\item $\lim \mathcal C\circ\alpha(t)$ is a curve of the type described
in \S\ref{statement}.
\end{itemize}
\end{defin}

Thus, marker germs `mark' one component of $\oE$,
and $\nG$ is particularly well-behaved around marker germs. Note that
since $\nG$ is normal, it is nonsingular along a dense open set
in each component of $\oE$: so the second requirement in the
definition of marker germ is satisfied for a general $\oalpha$ on
every component of $\oE$. It follows that every component of $\oE$
admits marker germs. Our next result is that, for
marker germs, our notion of `equivalence' (Definition~\ref{equivgerms}) 
translates nicely in terms of the lifts to $\oE$.

\begin{prop}\label{equivgermsagain}
Let $\alpha_0(t)$, $\alpha_1(t)$ be germs, and assume that
$\alpha_1(t)$ is a marker germ. Then $\alpha_0(t)$ is equivalent to
$\alpha_1(t)$ if and only if $\oalpha_0=\oalpha_1$ and the lift of
$\alpha_0(t)$ is transversal to $\oE$.
In particular, $\alpha_0(t)$ is then a marker germ as well.
\end{prop}

\begin{proof}
Assume first that the germs $\alpha_0(t)$ and $\alpha_1(t)$ are
equivalent. Then $\alpha_0$, $\alpha_1$ are fibers of a family
$\alpha_h(t)=A(h,t)$ where $A:\mathcal H\times \Spec \Cbb[[t]] \to
\Pbb^8$ is as specified in Definition~\ref{equivgerms}. In particular,
the map $\mathcal H \to \Til\Pbb^8$, 
$h\mapsto (\alpha_h(0),\lim_{t\to 0}\mathcal C\circ\alpha_h(t))$ 
is constant. This map factors
$$h \mapsto \oalpha_h \stackrel{\nor}\mapsto
(\alpha_h(0),\lim_{t\to 0}\,\mathcal C\circ\alpha_h(t))$$
and $\nor$ is finite over $(\alpha_h(0),\lim_{t\to 0}
\mathcal C\circ\alpha(t))$, so $h\mapsto
\oalpha_h$ is constant; in particular, $\oalpha_0=\oalpha_1$. 

Further, the weight of $F(\alpha_h(t))$ is constant as $h$ varies;
this implies that the intersection numbers of all lifts of the germs
$\alpha_h(t)$ with $\oE$ are equal. In particular, the lift of
$\alpha_0(t)$ is transversal to $\oE$ if and only if the lift of
$\alpha_1(t)$ is. Since $\alpha_1(t)$ is a marker germ, it follows
that the lift of $\alpha_0(t)$ is transversal to $\oE$, as needed.

For the other implication, $\nG$ is nonsingular at
$\oalpha:=\oalpha_0=\oalpha_1$ since $\alpha_1(t)$ is a marker germ.
Let $(\overline z)=(z_1,\dots,z_8)$ be a system of local
parameters for $\nG$ centered at $\oalpha$, and write the lifts
$\oalpha_i(t)$ as 
$$t\to (\overline z^{(i)}(t))=(z^{(i)}_1(t),\dots,
z^{(i)}_8(t))\quad.$$
Now
$$A(h,t):=\onu\left((1-h)\overline z^{(0)}(t)+h \overline
z^{(1)}(t)\right)$$
defines a map $\Abb^1\times \Spec \Cbb[[t]] \to \Pbb^8$, so that
$\alpha_h(t):=A(h,t)$ interpolates between $\alpha_0(t)$ and
$\alpha_1(t)$. The map $A(\_,0)$ is constant since all lifts
$\oalpha_h(t)$ meet $\oE$ at $\oalpha$. By the same token, writing
$$F(\alpha_h(t))=G(h) t^w+\text{higher order terms}\quad,$$
necessarily $G(h)=\rho(h) G$ for $\rho(h)\in\Cbb$ and $G$ a polynomial
in $x,y,z$ independent of $h$. Since both $\oalpha_0(t)$ and
$\oalpha_1(t)$ are transversal to $\oE$, we have that $\rho(0)$ and
$\rho(1)$ are both nonzero.
Taking $\mathcal H$ to be the complement of the zero-set of $\rho$
in $\Abb^1$ and restricting $A$ to $\mathcal H\times\Spec \Cbb[[t]]$
we obtain a map as prescribed in Definition~\ref{equivgerms}, showing
that $\alpha_0(t)$ and $\alpha_1(t)$ are equivalent.
\end{proof}

\begin{corol}\label{markergerms}
The germs $\alpha(t)$ obtained in section \ref{setth} in order to identify
components of the PNC are marker germs.
\end{corol}

\begin{proof}
Indeed, we have shown in \S\ref{setth} that {\em every\/}
germ determining a component of the PNC is equivalent to one such
germ $\alpha(t^d)$, with $d$ a positive integer; in particular this
holds with $d=1$ for marker germs, showing that such $\alpha(t)$ lift
to germs that are transversal to $\oE$, and meet it at points at which
$\nG$ is nonsingular.
\end{proof}

\subsection{}\label{action}
Another important ingredient is the $\PGL(3)$ action on $\nG$.

The $\PGL(3)$ action on $\Pbb^8$ given by multiplication on the right
makes the basic rational map
$$\Pbb^8 \dashrightarrow \Pbb^N$$
equivariant, and hence induces a right $\PGL(3)$ action on
$\Til\Pbb^8$ and $\nG$, fixing each component of $\oE$. Explicitly,
the action is realized by setting $\oalpha\cdot N$ to be the center of
the lift of
the germ $\alpha(t)\cdot N$, for $N\in \PGL(3)$. We record the
following trivial but useful remark:

\begin{lemma}\label{actionlift}
If $\beta(t)=\alpha(t)\cdot N$ for $N\in\PGL(3)$, then $\oalpha$ and
$\obeta$ belong to the same component of $\oE$.
\end{lemma}

\begin{proof}
Indeed, then $\obeta$ belongs to the $\PGL(3)$ orbit of $\oalpha$.
\end{proof}

\begin{lemma}\label{formal}
Let $\oD$ be a component of $\oE$ of type III,~IV,~or V.  The orbit of 
a general $\oalpha$ in $\oD$ is dense in $\oD$.
\end{lemma}

\begin{proof}
This follows immediately from the description of the general elements in the 
components of the PNC, given in \S\ref{statement}.
\end{proof}

Combining this observation with Lemma~\ref{actionlift} and
Proposition~\ref{equivgermsagain} gives a precise description of
the fibers of $\oE$ over $S$:

\begin{corol}\label{fiberdesc}
Let $\alpha(t)$, $\beta(t)$ be marker germs such that
$\alpha(0)=\beta(0)$. Then $\oalpha$ and $\obeta$ belong to the
same component of $\oE$ if and only if $\beta(t)$ is equivalent to
$\alpha(t)\cdot N$ for some constant invertible matrix $N$ such that
$\alpha(0)\cdot N=\alpha(0)$.
\end{corol}

\begin{proof}
Let $\alpha(t)$, $\beta(t)$ be marker germs such that $\alpha(0)=\beta(0)$. 
If $\oalpha$ and $\obeta$ belong to the same component of $\oE$, 
then by Lemma~\ref{formal} there exists
$N\in \PGL(3)$ such that $\obeta=\oalpha\cdot N$; by
Proposition~\ref{equivgermsagain}, $\beta(t)$ is equivalent to
$\alpha(t)\cdot N$. Further, $\alpha(0)=\beta(0)=\alpha(0)\cdot N$ (by
definition of equivalent germs), hence the stated condition on $N$
must hold.

The other implication is immediate from Lemma~\ref{actionlift}.
\end{proof}

This description yields our main tool for computing the degrees
$e_{ij}$, Proposition~\ref{degreetool} below. 

First, equivalence of marker germs can be recast in the following
apparently stronger form.

\begin{lemma}\label{strongeq}
Two marker germs $\alpha_0(t)$ and $\alpha_1(t)$ are equivalent
(w.r.t.~$\mathcal C$) if and only if there exists a unit $\nu(t)\in
\Cbb[[t]]$ and a $\Cbb[h][[t]]$-valued point $N(h,t)$ of $\PGL(3)$, 
such that
\begin{itemize}
\item $N(0,t)$ is the identity;
\item $N(h,0)$ is the identity; and
\item $\alpha_1(t\nu(t))=\alpha_0(t)\cdot N(1,t)$.
\end{itemize}
\end{lemma}

\begin{proof}
If a matrix $N(h,t)$ exists as in the statement, 
define $A: \Abb^1 \times \Spec\Cbb[[t]] \to \Pbb^8$ by setting
$A(h,t):=\alpha_0(t)\cdot N(h,t)$. Then the conditions prescribed by
Definition~\ref{equivgerms} are satisfied with $h_0=0$, $h_1=1$,
showing that $\alpha_0(t)$ is equivalent to $\alpha_1(t)$.

For the converse: since $\alpha_0(t)$ and $\alpha_1(t)$ are
equivalent, $\oalpha_0=\oalpha_1$ by
Proposition~\ref{equivgermsagain}, and $(\alpha_0(0),\lim_{t\to 0}
\mathcal C\circ\alpha_0(t))=(\alpha_1(0),\lim_{t\to 0}\mathcal C\circ
\alpha_1(t)).$ Let $(\alpha,\mathcal X)$ be this point of
$\Til\Pbb^8$, and let $D$ be the (unique) component of $E$ which contains it.
For the remainder of the argument, we use the standing assumption that
$D$ is a component of type III,~IV, or~V.

Under this assumption, the stabilizer of $(\alpha,\mathcal X)$ has
dimension~1; consider an $\Abb^7$ transversal to the stabilizer at the
identity $I$, let $U=\Abb^7\cap\PGL(3)$, and consider the action map 
$U\times \Spec\Cbb[[t]]\to \nG$:
$$(\varphi,t) \mapsto \oalpha_0(t)\circ\varphi$$
where $\oalpha_0(t)$ is the lift of $\alpha_0(t)$ to $\nG$.

Note that $\alpha_0(t)$ factors through this map:
$$\Spec\Cbb[[t]] \rightarrow U\times \Spec\Cbb[[t]] \rightarrow
\nG \rightarrow \Pbb^8$$
by
$$t \mapsto (I,t) \mapsto \oalpha_0(t) \mapsto \alpha_0(t)\quad.$$
Lifting $\alpha_1(t)$ we likewise get a factorization
$$t \mapsto (M(t),z(t)) \mapsto \oalpha_0(z(t))\circ M(t)=\oalpha_1(t) 
\mapsto \alpha_1(t)$$
for suitable $M(t)$, $z(t)$.
We may assume that the center $(M(0),z(0))$ of the lift of $\alpha_1(t)$ 
equals the center $(I,0)$ of the lift of $\alpha_0(t)$; also, $z(t)$ 
vanishes to order~1 at $t=0$, since the lift of $\alpha_1(t)$
is transversal to $\oE$.
Hence there exists a unit $\nu(t)$ such that $z(t\nu(t))=t$, and we
can effect the parameter change
$$\alpha_1(t\nu(t))=\alpha_0(t)\circ M(t\nu(t))=\alpha_0(t)\circ
N(t)\quad,$$
where we have set $N(t)=M(t\nu(t))$, a $\Cbb[[t]]$-valued point of 
$\PGL(3)$.

Finally, interpolate between $(I,t)$ and $(N(t),t)$ in 
$U\times\Spec\Cbb[[t]]\subset \Abb^7\times \Spec\Cbb[[t]]$ (much as we did 
already in the proof of Proposition~\ref{equivgermsagain}), by
$$(h,t) \mapsto ((1-h)I+h N(t),t)\quad.$$
Setting $N(h,t):=(1-h)I+h N(t)$, we obtain the sought
$\Cbb[h][[t]]$-valued point of 
$\PGL(3)$.
Indeed: $N(0,t)=I$ for all $t$; 
$N(h,0)=I$ for all $h$;
and $\alpha_0(t)\cdot N(1,t)=\alpha_0(t)\cdot N(t)= \alpha_1(t\nu(t))$.
\end{proof}

\subsection{}\label{inessential}
The dependence of this observation on a parameter change prompts the
following definition. Let $\alpha(t)$ be a marker germ for a component
$\oD$, and consider the $\Cbb((t))$-valued points of $\PGL(3)$ obtained
as products 
$$M_\nu(t):=\alpha(t)^{-1}\cdot \alpha(t \nu(t) )$$
as $\nu(t)$ ranges over all units in $\Cbb[[t]]$. Among all the
$M_\nu(t)$, consider those that are in fact $\Cbb[[t]]$-valued points
of $\PGL(3)$, and in that case let
$$M_\nu:=M_\nu(0)\quad.$$

Let $\alpha=\alpha(0)$, and $\mathcal X=\lim_{t\to 0}\mathcal
C\circ\alpha(t)$. 

\begin{lemma}\label{trivialrem}
The set of all $M_\nu$ so obtained is a subgroup of the stabilizer of
$(\alpha,\mathcal X)$.
\end{lemma}

\begin{proof}
Let $M_\nu$ be as above; that is, $M_\nu=M_\nu(0)$, where
$M_\nu(t)=\alpha(t)^{-1}\cdot \alpha(t\nu(t))$ is a $\Cbb[[t]]$-valued
point of $\PGL(3)$. Since
$$\alpha(t\nu(t))=\alpha(t)\cdot M_\nu(t)\quad,$$
we have
$$\alpha=\alpha(0\nu(0))=\alpha(0)\cdot M_\nu(0)=\alpha\cdot
M_\nu\quad,$$
showing that $M_\nu$ stabilizes $\alpha$. Further
$$\mathcal C\circ \alpha(t\nu(t))=\mathcal C\circ \alpha(t)\circ
M_\nu(t)\quad,$$
and taking the limit
$$\lim_{t\to 0} C\circ\alpha(t\nu(t))=(\lim_{t\to 0}
C\circ\alpha(t))\circ M_\nu(0)\quad,$$
that is,
$$\mathcal X=\mathcal X\circ M_\nu$$
since the limits along $\alpha(t)$ and $\alpha(t\nu(t))$ must
agree as the two germs only differ by a change of parameter. Thus
$M_\nu$ stabilizes $\mathcal X$ as well, as needed.

To verify that the set $\{M_\nu\}_\nu$ forms a group, let
$\overline\nu(t)$ be the unit such that
$$\overline\nu(t\nu(t))=\nu(t)^{-1}\quad.$$
Then we find
$$\aligned
M_{\nu}(t)\cdot M_{\overline\nu}(t\nu(t)) 
&= \alpha(t)^{-1} \alpha(t\nu(t))\,\alpha(t\nu(t))^{-1}
\alpha(t\nu(t) \overline\nu(t\nu(t)))\\
&=\alpha(t)^{-1} \alpha(t\nu(t) \nu(t)^{-1})\\ 
&=I\quad,
\endaligned$$
the identity matrix. That is, $M_{\overline\nu}(t\nu(t))
=M_{\nu}(t)^{-1}$, and hence
$$M_{\overline\nu}=M_{\overline\nu}(0)= M_{\nu}(0)^{-1}
=M_{\nu}^{-1}\quad.$$
Similarly, for $\nu_1(t)$, $\nu_2(t)$ units, let
$$\nu_3(t)=\nu_1(t) \nu_2(t\nu_1(t))\quad.$$
Then we find
$$\aligned
M_{\nu_3}(t) &=\alpha(t)^{-1} \alpha(t\nu_3(t))\\
&=\alpha(t)^{-1}\alpha(t\nu_1(t)) \alpha(t\nu_1(t))^{-1}
\alpha(t\nu_1(t)\nu_2(t\nu_1(t)))\\
&= M_{\nu_1}(t)M_{\nu_2}(t\nu_1(t))\quad,
\endaligned$$
and hence $M_{\nu_3}=M_{\nu_3}(0)=M_{\nu_1}(0) M_{\nu_2}(0)=M_{\nu_1}
M_{\nu_2}$ as needed.
\end{proof}

We call {\em inessential\/} the components of the stabilizer of
$(\alpha, \mathcal X)$ containing elements of the subgroup identified
in Lemma~\ref{trivialrem}. These components form the {\em inessential
  subgroup\/} of the stabilizer of $(\alpha, \mathcal X)$,
corresponding to the choice of $\alpha(t)$.

This apparently elusive notion is crucial for our tool to
compute the degrees~$e_{ij}$.

\begin{prop}\label{degreetool}
Let $D$ be a component of $E$, let $\oD$ be any component of $\oE$
dominating $D$. Then the degree of $\oD$ over $D$ is the index of the
inessential subgroup in the stabilizer of a general point of $D$.
\end{prop}

\begin{proof}
Let $(\alpha,\mathcal X)$ be a general point of $D$, and let
$\oalpha_1,\dots,\oalpha_r$ be its preimages in $\oD$; so the degree
of $\oD$ over $D$ equals $r$. We may assume that $\mathcal X$ is a
limit curve of the type described in \S\ref{statement}, and that $\nG$
is nonsingular at all $\oalpha_i$,
hence for $i=1,\dots,r$ we may choose a {\em marker\/} germ
$\alpha_i(t)$ whose lift is centered at $\oalpha_i$. Let
$\alpha(t)=\alpha_1(t)$.
By Corollary~\ref{fiberdesc}, every $\alpha_i(t)$ is equivalent to
$\alpha(t)\cdot N_i$ for some $N_i\in \PGL(3)$ stabilizing
$(\alpha,\mathcal X)$;
conversely, if $N\in PGL(3)$ fixes $(\alpha,\mathcal X)$ then
$\alpha(t)\cdot N$ is equivalent to one of the $\alpha_i(t)$.
Thus, the action $N\mapsto \oalpha\cdot N$ maps the stabilizer of
$(\alpha,\mathcal X)$ onto the fiber of $\oD$ over $(\alpha,\mathcal X)$.

Therefore we simply have to check that two elements of the stabilizer map
to the same point in $\oE$ if and only if they are in the same coset
w.r.t.~the inessential subgroup; that is, it is enough to check that
$\oalpha=\oalpha\cdot N$ if and only if $N$ is in the inessential
subgroup.

First assume that $N$ is in the inessential subgroup, that is, $N$ is
in a component containing an element $M_\nu$ as above. Since the fiber
of $\oE$ over $(\alpha,\mathcal X)$ is finite,
then $\oalpha\cdot N=\oalpha\cdot M_\nu$. Now $\alpha(t)\cdot M_\nu$
and $\alpha(t)\cdot M_\nu(t)$ are equivalent by Lemma~\ref{stequiv};
since $\alpha(t)\cdot M_\nu(t)=\alpha(t\nu(t))$, we have that
$\oalpha\cdot M_\nu=\oalpha$. Thus $\oalpha\cdot N=\oalpha$ as needed.

For the converse, assume $\oalpha=\oalpha\cdot N$. By
Proposition~\ref{equivgermsagain} and Lemma~\ref{strongeq}, if
$\oalpha=\oalpha\cdot N$ then there is a $\Cbb[h][[t]]$-valued point
$N(h,t)$ of $\PGL(3)$ and a unit $\nu(t)$ in $\Cbb[[t]]$ such that
$N(0,t)=N(h,0)=I$ and
$$\alpha(t)\cdot N(1,t)=\alpha(t\nu(t))\cdot N.$$

Now $M_\nu(t):=\alpha(t)^{-1}\alpha(t\nu(t))=N(1,t)N^{-1}$ is a
$\Cbb[[t]]$-valued point of $\PGL(3)$, thus 
$M_\nu=N(1,0)N^{-1}=N^{-1}$ is in the inessential subgroup. 
This shows that $N$ is in the inessential subgroup, completing the proof.
\end{proof}

\subsection{}\label{weights}
Our work in \S\ref{setth} has produced a list of marker germs; these
can be used to evaluate the multiplicities $m_{ij}$.

\begin{defin}
The {\em weight\/} of a germ $\alpha(t)$ is the order of vanishing in
$t$ of $F\circ\alpha(t)$.
\end{defin}

Note that the weight of $\alpha(t)$ is the `order of contact' of
$\alpha(t)$ with $S$: indeed, it is the minimum intersection
multiplicity of $\alpha(t)$ and generators
$F\circ\varphi(x_0,y_0,z_0)$ of the ideal of $S$, at $\alpha(0)$.

\begin{lemma}\label{weightsmult}
The multiplicity $m_{ij}$ is the minimum weight of a germ $\alpha(t)$
such that $\oalpha\in \oE_{ij}$. This weight is achieved by a marker
germ for $\oE_{ij}$.
\end{lemma}

\begin{proof}
Let $\oalpha$ be a general point of $\oE_{ij}$. Since $\nG$ is normal,
we may assume that it is nonsingular at $\oalpha$. Let $(\overline
z)=(z_1,\dots,z_8)$ be a system of local parameters for $\nG$ centered
at $\oalpha$, and such that the ideal of $\oE_{ij}$ is $(z_1)$ near
$\oalpha$; thus the ideal of $\oE$ is $(z_1^{m_{ij}})$ near $\oalpha$.
Consider the germ $\oalpha(t)$ in $\nG$ defined by
$$\oalpha(t)=(t,0,\dots,0)\quad,$$
and its push-forward $\alpha(t)=\onu(\oalpha(t))$ in $\Pbb^8$.

The weight of $\alpha(t)$ is the order of contact of $\alpha(t)$ with
$S$; hence it equals the order of contact of $\oalpha(t)$ with
$\onu^{-1}(S)=\oE$; pulling back the ideal of $\oE$ to $\oalpha(t)$,
we see that this equals $m_{ij}$.

In fact, this argument shows that the weight of any germ in $\Pbb^8$
lifting to a germ in $\nG$ meeting the support of $\oE$ {\em
transversally\/} at a general point of $\oE_{ij}$ is $m_{ij}$. It
follows that the weight of {\em any\/} germ lifting to one meeting
$\oE_{ij}$ must be $\ge m_{ij}$, completing the proof of the first
assertion.

The second assertion is immediate, as the germ $\alpha(t)$ constructed
above is a marker germ for $\oE_{ij}$.
\end{proof}

Applying Proposition~\ref{multcount} requires the list of the
components $\oE_{ij}$ of $\oE$ dominating a given component $E_i$ of
$E$, and for each $\oE_{ij}$ the two numbers $e_{ij}$ and
$m_{ij}$. These last two elements of information will be obtained by
applying Proposition~\ref{degreetool} and
Lemma~\ref{weightsmult}. Obtaining the list $\oE_{ij}$ and a local
description of $\nG$ requires a case-by-case analysis. 

\subsection{}\label{multtypeI}
We start with components of type~I in this subsection. As it happens,
$\nG \to \Til \Pbb^8$ is an isomorphism near the general point of such
a component, and we can perform the multiplicity computation directly
on $\Til\Pbb^8$.

\begin{prop}\label{typeImul}
Assume $\mathcal C$ contains a line $\ell$ with multiplicity $m$, and
let $(\alpha, \mathcal X)$ be a general point of the corresponding
component $D$ of $E$. Then $\Til\Pbb^8$ is nonsingular near
$(\alpha,\mathcal X)$, and $D$ appears with multiplicity $m$ in $E$. 
\end{prop}

\begin{proof}
We are going to show that, in a neighborhood of $(\alpha, \mathcal
X)$, $\Til\Pbb^8$ is isomorphic to the blow-up of $\Pbb^8$ along
the $\Pbb^5$ of matrices whose image is contained in $\ell$. 
The nonsingularity of $\Til\Pbb^8$ near $(\alpha, \mathcal X)$ follows
from this.

Choose coordinates so that $\ell$ is the line $z=0$, and the affine
open set $U$ in $\Pbb^8$ with coordinates
$$\begin{pmatrix}
1 & p_1 & p_2 \\
p_3 & p_4 & p_5 \\
p_6 & p_7 & p_8
\end{pmatrix}$$
contains $\alpha$. The $\Pbb^5$ of matrices with image contained in
$\ell$ intersects this open set along $p_6=p_7=p_8=0$, so we can
choose coordinates $q_1,\dots,q_8$ in an affine open subset $V$ of the
blow-up of $\Pbb^8$ along $\Pbb^5$ so that the blow-up map is given by
$$\left\{
\aligned
p_i &=q_i\quad i=1,\dots,5 \\
p_6 &=q_6 \\
p_7 &=q_6 q_7\\
p_8 &=q_6 q_8
\endaligned\right.$$
(the part of the blow-up over $U$ is covered by three such open sets;
it should be clear from the argument that the choice made here is
immaterial).

Under the hypotheses of the statement, the ideal of $\mathcal C$ is
generated by $z^m G(x,y,z)$, where $z$ does not divide $G$; that is,
$G(x,y,0)\ne 0$. The rational map $\Pbb^8 \dashrightarrow \Pbb^N$ acts
on $U$ by sending $(p_1,\dots,p_8)$ to the curve with ideal generated
by
$$(p_6 x+p_7 y+p_8 z)^m G(x+p_1 y+p_2 z,p_3 x+p_4 y+p_5 z, p_6 x+p_7
y+p_8 z)\quad.$$ 
Composing with the blow-up map:
$$V \longrightarrow U \dashrightarrow \Pbb^N$$
we find that $(q_1,\dots,q_8)$ is mapped to the curve with ideal
generated by
$$(x+q_7 y+q_8 z)^m G(x+q_1 y+q_2 z,q_3 x+q_4 y+q_5 z, q_6 x+q_6
q_7 y+q_6 q_8 z)\quad,$$ 
where a factor of $q_6^m$ has been eliminated. Note that no other
factor of $q_6$ can be extracted, by the hypothesis on $G$.

A coordinate verification shows that the induced map
$$V \longrightarrow U\times \Pbb^N \subset \Pbb^8\times \Pbb^N\quad,$$
which clearly maps $V$ to $\Til\Pbb^8\subset \Pbb^8\times\Pbb^N$ and
the exceptional divisor $q_6=0$ to $D$, is an isomorphism onto the
image in a neighborhood of a general point of the exceptional divisor,
proving that $\Til\Pbb^8$ is nonsingular in a neighborhood of the
general $(\alpha, \mathcal X)$ in $D$. 

For the last assertion in the statement, pull-back the generators of
the ideal of $S$ to $V$:
$$q_6^m (x+q_7 y+q_8 z)^m G(x+q_1 y+q_2 z,q_3 x+q_4 y+q_5 z, q_6 x+q_6
q_7 y+q_6 q_8 z)\quad,$$ 
as $(x:y:z)$ ranges over $\Pbb^2$. For $(x:y:z)=(1:0:0)$ this gives
the generator
$$q_6^m G(1,q_3,q_6)\quad.$$
At a general $(\alpha, \mathcal X)$ in $D$ we may assume that
$G(1,q_3,0)\ne 0$ (again by the hypothesis on $G$), and we find that
the ideal of $E$ is $(q_6^m)$ near $(\alpha, \mathcal X)$. This shows
that the multiplicity of the component is $m$, as stated.
\end{proof}

Proposition~\ref{typeImul} yields the multiplicity statement
concerning type~I components in \S\ref{multstat}; also cf.~Fact~2~(i)
in \S2 of \cite{MR2001h:14068}.

\subsection{}\label{multtypeII} 
Components of type~II can be analyzed almost as explicitly as
components of type~I, without employing the tools developed in
\S\S\ref{beginprel}--\ref{weights}.

Recall from \S\ref{typeII} that every nonlinear component $\mathcal
C'$ of a curve $\mathcal C$ determines a component $D$ of the
exceptional divisor $E$.

\begin{prop}\label{typeIImul}
Assume $\mathcal C$ contains a nonlinear component $\mathcal C'$, with
multiplicity $m$, and let $D$ be the corresponding component of $E$.
Then $D$ appears with multiplicity $2m$ in $E$.
\end{prop}

This can be proved by using the blow-ups described in
\cite{MR94e:14032}, which resolve the indeterminacies of the basic
rational map $\Pbb^8 \dashrightarrow \Pbb^N$ over nonsingular 
non-inflectional points
of $\mathcal C$. We sketch the argument here, leaving detailed
verifications to the reader.

\begin{proof}
In \cite{MR94e:14032} it is shown that two blow-ups at smooth centers
suffice over nonsingular, non-inflectional points of $\mathcal
C$. While the curve was assumed to be reduced and irreducible in
loc.~cit., the reader may check that the same blow-ups resolve the
indeterminacies over a possibly multiple component $\mathcal C'$, near
nonsingular, non-inflectional points of the support of $\mathcal C'$.
Let $V$ be the variety obtained after these two blow-ups.

Since the basic rational map is resolved by $V$ over a general point
of $\mathcal C'$, the inverse image of the base scheme $S$ is locally
principal in $V$ over such points. By the universal property of
blow-ups, the map $V \to \Pbb^8$ factors through $\Til\Pbb^8$ over a
neighborhood of a general point of $\mathcal C'$.
It may then be checked that the second exceptional divisor obtained in
the sequence maps birationally onto~$D$, and appears with a
multiplicity of $2m$.

The statement follows.
\end{proof}

Proposition~\ref{typeIImul} yields the multiplicity statement
concerning type~II components in \S\ref{multstat}; also cf.~Fact~2~(ii)
in \S2 of \cite{MR2001h:14068}.

\subsection{}\label{multtypeIII}
Next, we consider components of type~III. Recall that there is one such
component for every singular point $p$ of $\mathcal C$ at which the
tangent cone to $\mathcal C$ consists of at least three lines, and
that we have shown (cf.~Propositions~\ref{faber} and~\ref{tgcone})
that every marker germ $\alpha(t)$ leading to one of these components is
equivalent to one which, in suitable coordinates $(x:y:z)$, may be
written as 
$$\begin{pmatrix}
1 & 0 & 0\\
0 & t & 0\\
0 & 0 & t
\end{pmatrix}\quad,$$
where $p$ has coordinates $(1:0:0)$ and the kernel line has equation
$x=0$. Call $D$ the component of $E$ corresponding to one such point;
note that $D$ dominates the subset of $S$ consisting of matrices whose
image is the point $(1:0:0)$.

\begin{prop}\label{typeIIImul}
There is exactly one component $\oD$ of $\oE$ dominating $D$, and the
degree of the map $\oD \to D$ equals the number of linear
automorphisms of the tuple determined by the tangent cone to $\mathcal
C$ at $p$. The minimal weight of a germ leading to $D$ equals the
multiplicity of $p$ on $\mathcal C$.
\end{prop}

In this statement, the {\em linear automorphisms\/} of a tuple of
points in $\Pbb^1$ are the elements of its $\PGL(2)$-stabilizer; this
is a finite set if and only if the tuple is supported on at least
three points.

\begin{proof}
First we will verify that if $\alpha(t)$, $\beta(t)$ are marker germs
leading to $D$, then $\oalpha$ and $\obeta$ belong to the same
component of $\oE$. This will show that there is only one component of
$\oE$ over $D$. By Proposition~\ref{equivgermsagain}, we may replace
$\alpha(t)$ and $\beta(t)$ by equivalent germs; as recalled above, in
suitable coordinates we may then assume
$$\alpha(t)=\begin{pmatrix}
1 & 0 & 0 \\
0 & t & 0 \\
0 & 0 & t
\end{pmatrix}\quad,$$
and $\beta(t)$ may be assumed to be in the same form after a change of
coordinates. That is, we may assume
$$\beta(t)=M\cdot \alpha(t)\cdot N$$
for constant matrices $M,N$.
As $D$ dominates the subset of $S$ consisting of matrices with image 
$(1:0:0)$, the image of $\alpha(0)$ and $\beta(0)$ is necessarily
$(1:0:0)$; this implies that $M$ is in the form
$$\begin{pmatrix}
A & B & C \\
0 & E & F \\
0 & H & I
\end{pmatrix}$$
with $A(EI-FH)\ne 0$.

Now we claim that $M\cdot\alpha(t)\cdot N$ is equivalent to
$\alpha(t)\cdot N'$ for another constant matrix $N'$. Indeed, note that 
$$M\cdot \alpha(t)\cdot N=\alpha(t)\cdot\begin{pmatrix}
A & Bt & Ct \\ 0 & E & F \\ 0 & H & I \end{pmatrix}\cdot N
=: \alpha(t)\cdot N'(t)\quad.$$
Writing $N'=N'(0)$ and applying Lemma~\ref{stequiv}, we establish the claim.

By Proposition~\ref{equivgermsagain}, we may thus assume that
$\beta(t)=\alpha(t)\cdot N'$. It follows that $\obeta$ and $\oalpha$
are on the same component of $\oE$, by Lemma~\ref{actionlift}.

The degree of $\oD \to D$ is evaluated by using
Proposition~\ref{degreetool}. By the description in \S\ref{statement},
the limit of $\mathcal C$ along a marker germ $\alpha(t)$ as above
consists of a fan $\mathcal X$ whose star reproduces the tangent cone
to $\mathcal C$ at $p$, and whose free line is supported on the kernel
line $x=0$. It is easily checked that the stabilizer of
$(\alpha(0),\mathcal X)$ has one component for each element
of $\PGL(2)$ fixing the tuple determined by the tangent cone to
$\mathcal C$ at $p$ and that the inessential subgroup equals the
identity component.

Finally, the foregoing considerations show that the minimal weight of
a germ leading to $D$ is achieved by $\alpha(t)$, and this is
immediately computed to be the multiplicity of $\mathcal C$ at $p$.
\end{proof}

By Proposition~\ref{multcount}, Proposition~\ref{typeIIImul} implies
the multiplicity statement for type~III components in
\S\ref{multstat}; also cf.~Fact~4~(i) in \S2 of \cite{MR2001h:14068}.

\subsection{}\label{multtypeIV}
Recall from \S\ref{setth} that components of type~IV arise from
certain sides of the Newton polygon determined by the choice of a
point $p$ on $\mathcal C$ and of a line in the tangent cone to
$\mathcal C$ at $p$. If coordinates $(x:y:z)$ are chosen so that
$p=(1:0:0)$, and the tangent line is the line $z=0$, then the Newton
polygon (see \S\ref{Newton}) consists of the convex hull of the union
of the positive quadrants with origin at the points $(j,k)$ for which
the coefficient of $x^iy^jz^k$ in the equation for $\mathcal C$ is
nonzero. The part of the Newton polygon consisting of line segments
with slope strictly between $-1$ and $0$ does not depend on the choice
of coordinates fixing the flag $z=0$, $p=(0,0)$. We have found (see
Proposition~\ref{Newtonsides} and ff.) that if $-b/c$ is a slope
of the Newton polygon, with $b$, $c$ relatively prime, then
$$\alpha(t)=\begin{pmatrix}
1 & 0 & 0 \\
0 & t^b & 0 \\
0 & 0 & t^c
\end{pmatrix}$$
is a marker germ for a component of type~IV if $p$ is a singular or
inflection point of the support of $\mathcal C$, and the limit
$$x^{\overline q}y^rz^q \prod_{j=1}^S(y^c+\rho_j x^{c-b}z^b)$$
is not supported on a conic union (possibly) the kernel line.

It is clear that germs arising from sides of Newton polygons
corresponding to different lines in the tangent cone at $p$ cannot be
equivalent, so we concentrate on one side of one polygon. 

\begin{lemma}\label{typeIVmul1}
Let $D$ be the component of $E$ of type~IV corresponding to one side
of the Newton polygon of slope $-b/c$, as above. Then there is exactly
one component $\oD$ of $\oE$ over~$D$ corresponding to this side, 
and the degree of the map $\oD \to D$ equals the number of
automorphisms $\Abb^1 \to \Abb^1$, $\rho\mapsto u\rho$ (with $u$ a
root of unity) preserving the $S$-tuple $\{\rho_1,\dots,\rho_S\}$.
\end{lemma}

\begin{proof}
As in the proof of Proposition~\ref{typeIIImul}, we begin by verifying
that if $\alpha(t)$, $\beta(t)$ are marker germs leading to $D$, then
$\oalpha$ and $\obeta$ belong to the same component, by essentially
the same strategy.

By Proposition~\ref{equivgermsagain}, we may replace $\alpha(t)$ and
$\beta(t)$ with equivalent germs; thus we may choose
$$\alpha(t)=\begin{pmatrix}
1 & 0 & 0 \\
0 & t^b & 0 \\
0 & 0 & t^c
\end{pmatrix}$$
and $\beta(t)=M\cdot \alpha(t) \cdot N$ for constant invertible
matrices $M$ and $N$. As $\beta(t)$ leads to $D$, the matrix $M$ must
preserve the flag consisting of $p=(1:0:0)$ and the line $z=0$ (as
this is the data which determines the Newton polygon).
This implies that
$$M=\begin{pmatrix}
A & B & C \\
0 & E & F \\
0 & 0 & I
\end{pmatrix}$$
with $AEI\ne 0$. 
Note that 
$$M\cdot \alpha(t)=\alpha(t)\cdot\begin{pmatrix}
A & t^b & t^c \\ 0 & E & t^{c-b} \\ 0 & 0 & I \end{pmatrix}\quad;$$
applying Lemma~\ref{stequiv}, we find
that $M\cdot \alpha(t)\cdot N$ is equivalent to
$$\alpha(t)\cdot \begin{pmatrix}
A & 0 & 0 \\
0 & E & 0 \\
0 & 0 & I
\end{pmatrix}
\cdot N\quad,$$
since $c>b>0$.

We conclude that $\beta(t)$ is equivalent to $\alpha(t)\cdot N'$ for
some constant invertible matrix $N'$, and Lemma~\ref{actionlift} then
implies that $\oalpha$, $\obeta$ belong to the same component of
$\oD$, as needed.

Next, we claim that the inessential subgroup of the stabilizer
consists of the component of the identity; by
Proposition~\ref{degreetool}, it follows that the degree of the map
$\oD \to D$ equals the number of components of the stabilizer of
$(\alpha, \mathcal X):= (\alpha(0),\lim_{t\to 0}\mathcal C\circ\alpha(t))$. 

To verify our claim, note that for all units $\nu(t)$
$$\alpha(t)^{-1}\cdot \alpha(t\nu(t))=
\begin{pmatrix}
1 & 0 & 0 \\
0 & t^b & 0 \\
0 & 0 & t^c
\end{pmatrix}^{-1}\cdot\begin{pmatrix}
1 & 0 & 0 \\
0 & t^b\nu(t)^b & 0 \\
0 & 0 & t^c\nu(t)^c
\end{pmatrix}=
\begin{pmatrix}
1 & 0 & 0 \\
0 & \nu(t)^b & 0\\
0 & 0 & \nu(t)^c
\end{pmatrix}$$
is a $\Cbb[[t]]$-valued point of $\PGL(3)$. Thus the inessential
components of the stabilizer are those containing elements
$$\begin{pmatrix}
1 & 0 & 0 \\
0 & \nu^b & 0\\
0 & 0 & \nu^c
\end{pmatrix}\quad,$$
for $\nu\in\Cbb$, $\nu\ne 0$. This is the component containing the
identity, as claimed.

The number of components of the stabilizer of $(\alpha,\mathcal X)$ 
is determined as follows. The limit is
$$x^{\overline q}y^rz^q \prod_{j=1}^S(y^c+\rho_j x^{c-b}z^b)\quad,$$
and the orbit of $(\alpha,\mathcal X)$ has dimension~7.
The degree of $\oD \to D$ equals the number of components of the
stabilizer of $(\alpha,\mathcal X)$, that is, the subset of the
stabilizer of $\mathcal X$ fixing the kernel line $x=0$. 
If the orbit of $\mathcal X$ has dimension~7, this number
equals the number of components of the stabilizer of $\mathcal X$, or
the same number divided by~2, according to whether the kernel line is
identified by $\mathcal X$ or not. The latter eventuality occurs
precisely when $c=2$ and $q=\overline q$; the stated conclusion
follows then from Lemma~3.1 in \cite{MR2002d:14083}.
Analogous arguments apply when the orbit of $\mathcal X$ has dimension
less than~7.
\end{proof}

In order to complete the proof of the multiplicity statement for
type~IV components we just need a weight computation.

\begin{lemma}\label{typeIVmul2}
The minimum weight of a marker germ for the component corresponding to
a side of the Newton polygon with vertices $(j_0,k_0)$, $(j_1,k_1)$,
$j_0<j_1$, is
$$\frac{j_1k_0-j_0k_1}S\quad.$$
\end{lemma}

\begin{proof}
With notations as above, we have seen that every marker germ leading to the
component is equivalent to $\alpha(t)\cdot N$, where 
$$\alpha(t)=\begin{pmatrix}
1 & 0 & 0 \\
0 & t^b & 0 \\
0 & 0 & t^c
\end{pmatrix}\quad;$$
thus the minimum weight is achieved by this germ.

The limit
$$x^{\overline q}y^rz^q \prod_{j=1}^S(y^c+\rho_j x^{c-b}z^b)$$
appears with weight $br+cq+Sbc$, so we just have to show that
$$Sbc+br+cq=\frac{j_1k_0-j_0k_1}S\quad.$$
This is immediate, as $(j_0,k_0)=(r,q+Sb)$, $(j_1,k_1)=(r+Sc,q)$.
\end{proof}

The prescription for type~IV components now follows from
Lemmas~\ref{typeIVmul1} and~\ref{typeIVmul2},
Proposition~\ref{multcount}, and Lemma~\ref{weightsmult}.  

We note that the same prescription yields the correct multiplicity for
type~II limits as well: indeed, the side of the Newton polygon
corresponding to type~II limits (as in \S\ref{typeII}) has vertices
$(0,m)$ and $(2m,0)$, where $m$ is the multiplicity of the
corresponding nonlinear component of $\mathcal C$; so $S=m$ and
$(j_1k_0-j_0k_1)/S=2m^2/m=2m$, in agreement with
Proposition~\ref{typeIImul}.

Fact~4(ii) in \cite{MR2001h:14068}, \S2, reproduces the result proved
here; the reader should note that the number denoted $A$ here is
denoted $A/\delta$ in loc.~cit.

The weight computed in Lemma~\ref{typeIVmul2}:
$$\frac{j_1k_0-j_0k_1}S\quad,$$
happens to equal $2/S$ times the area of the triangle with vertices
$(0,0)$, $(j_0,k_0)$, and $(j_1,k_1)$. 

An example illustrating this for $b/c=1/3$:
$$\includegraphics[scale=.66]{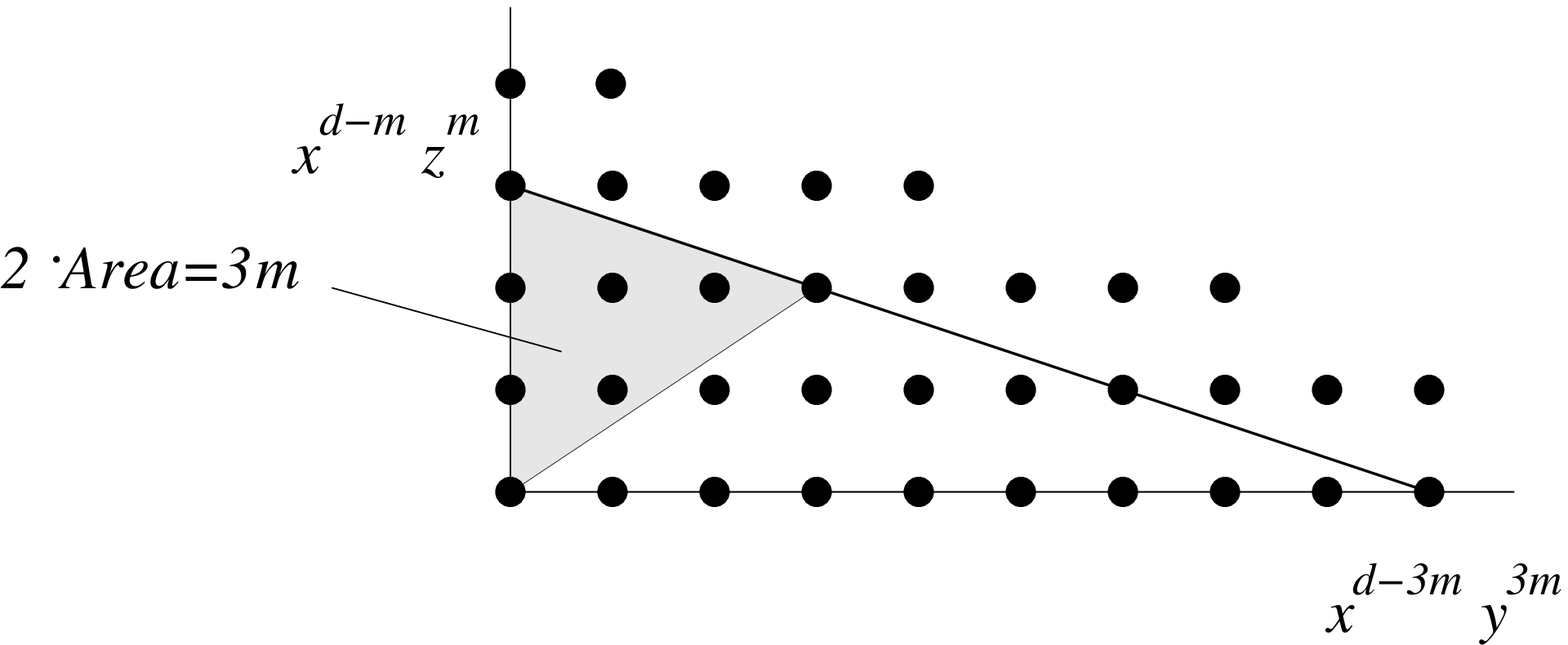}$$
We do not have a conceptual explanation for this observation.

\subsection{}\label{multtypeV}
We are left with components of type~V, whose analysis is predictably
subtler. In \S\ref{setth} we have found that such components arise
from suitable truncations of the Puiseux expansion of the branches of
$\mathcal C$ at a singular point $p$ of its support. Choosing
coordinates $(x:y:z)$ so that $p=(1:0:0)$ and the branch has tangent
cone supported on the line $z=0$, we have found a marker germ
$$\alpha(t)=\begin{pmatrix}
1 & 0 & 0 \\
t^a & t^b & 0 \\
\underline{f(t^a)} & \underline{f'(t^a)t^b} & t^c
\end{pmatrix}$$
where 
$$z=f(y)=\sum_{i\ge 0} \gamma_{\lambda_i} y^{\lambda_i}$$
is the Puiseux expansion of a corresponding formal branch, $a<b<c$ are
positive integers, $C=\frac ca$ is `characteristic' in the sense
explained in \ref{statement}, underlining denotes truncation to $t^c$,
and $a$ is the smallest positive integer for which all entries 
in $\alpha(t)$ are polynomials. The limit obtained along
this germ is:
$$x^{d-2S}\prod_{i=1}^S\left(zx-\frac {\lambda_0(\lambda_0-1)}2
\gamma_{\lambda_0}y^2 -\frac{\lambda_0+C}2
\gamma_{\frac{\lambda_0+C}2}yx-\gamma_C^{(i)}x^2\right)\quad,$$
where $\gamma_C^{(i)}$ are the coefficients of $y^C$ for all formal
branches sharing the truncation
$f_{(C)}(y)=\sum_{\lambda_i<C}\gamma_{\lambda_i}y^{\lambda_i}$.

As the situation is more complex than for other components, we proceed
through the proof of the multiplicity statement given in
\S\ref{multstat} one step at the time.

\subsection{}
Through the procedure recalled above, the choice of a characteristic
$C$ and of a truncation $f_{(C)}(y)$ of a formal branch determines a
germ, and hence (by lifting to the normalization) a component of $\oE$
over a fixed type~V component $D$ of $E$. In fact, by
Proposition~\ref{standardform} and Lemma~\ref{actionlift}, every
component $\oD$ over $D$ is marked in this fashion.

Clearly different points or different lines in the tangent cone yield
different components $\oD$ over $D$. 
As stated in \S\ref{multstat}, for a fixed point and line there are
different contributions for truncations that are not `siblings'. In
other words, we must show that there is a bijection between the set
of components $\oD$ over $D$ corresponding to a given point and line,
and set of data $(C,f_{(C)}(y))$ as above, modulo the sibling
relation. We will now recall this notion, and prove this fact in
Proposition~\ref{typeVmultsiblings} below.

We say that $(C,f_{(C)}(y))$, $(C',g_{(C')}(y))$ (or the truncations
$f_{(C)}$, $g_{(C')}$ for short) are {\em sibling\/} data if the
corresponding integers $a<b<c$, $a'<b'<c'$ are the same (so in
particular $C=C'$) and further 
$$g_{(C)}(y)=\sum_{\lambda_i<C} \xi^{a\lambda_i} \gamma_{\lambda_i}
y^{\lambda_i}$$
(that is, $f_{(C)}(t^a)=g_{(C)}((\xi t)^a)$) for an $a$-th root $\xi$
of $1$.

Both Proposition~\ref{typeVmultsiblings} and the determination of the
inessential subgroup rely on the following technical lemma.

\begin{lemma}\label{technical}
Let 
$$\alpha(t)=\begin{pmatrix}
1 & 0 & 0 \\
t^a & t^b & 0 \\
\underline{f(t^a)} & \underline{f'(t^a)t^b} & t^c 
\end{pmatrix}\quad,\quad
\beta(t)=\begin{pmatrix}
1 & 0 & 0 \\
t^{a'} & t^{b'} & 0 \\
\underline{g(t^{a'})} & \underline{g'(t^{a'})t^{b'}} & t^{c'} 
\end{pmatrix}$$
be two marker germs of the type considered above, and assume that
$\alpha(t)^{-1}\beta(\tau(t))$ is a $\Cbb[[t]]$-valued point of
$\PGL(3)$, for a change of parameter $\tau(t)=t\nu(t)$ with
$\nu(t)\in\Cbb[[t]]$ a unit. Then $a'=a$, $b'=b$, $c'=c$, and
$\nu(t)=\xi(1+t^{b-a}\mu(t))$, where $\xi$ is an $a$-th root of~$1$
and $\mu(t)\in\Cbb[[t]]$; further, $\underline{g((\xi
t)^a)}=\underline{f(t^a)}$.
\end{lemma}

\begin{proof}
Write $\varphi(t)=\underline{f(t^a)}$ and
$\psi(t)t^b=\underline{f'(t^a)t^b}$. The hypothesis is that
$$\alpha(t)^{-1}\cdot\beta(\tau)=\begin{pmatrix}
1 & 0 & 0 \\
\frac{\tau^{a'}-t^a}{t^b} & \frac{\tau^{b'}}{t^b} & 0 \\
\frac{\underline{g(\tau^{a'})}-\varphi(t)-(\tau^{a'}-t^a)\psi(t)}
     {t^c} &
\frac{\underline{g'(\tau^{a'})\tau^{b'}}-\psi(t)\tau^{b'}}{t^c} &
\frac{\tau^{c'}}{t^c}
\end{pmatrix}$$
has entries in $\Cbb[[t]]$, and its determinant is a unit in
$\Cbb[[t]]$. The latter condition implies $b'=b$ and $c'=c$. As 
$$\frac{\tau^{a'}-t^a}{t^b}\in \Cbb[[t]]\quad,$$
necessarily $a'=a$ and $t^a(\nu(t)^a-1)=(\tau^a-t^a)\equiv
0\mod{t^b}$. Since $b>a$, this implies 
$$\nu(t)=\xi(1+t^{b-a}\mu(t))$$ 
for $\xi$ an $a$-th root of $1$ and $\mu(t)\in\Cbb[[t]]$. 
Also note that since the triples $(a,b,c)$ and $(a',b',c')$ coincide,
necessarily the dominant term in $g(y)$ has the same exponent
$\lambda_0$ as in $f(y)$,
since $a\lambda_0=2a-2b+c$.

Now we claim that
$$g(\tau^a)-(\tau^a-t^a)\psi(t)\equiv g((\xi t)^a)\mod t^c\quad.$$
Granting this for a moment, it follows that
$$\underline{g(\tau^{a'})}-\varphi(t)-(\tau^{a'}-t^a)\psi(t) \equiv
g((\xi t)^a)-f(t^a)\mod t^c\quad;$$
hence, the fact that the $(3,1)$ entry is in $\Cbb[[t]]$ implies that
$$g((\xi t)^a)\equiv f(t^a)\mod t^c\quad,$$
which is what we need to show in order to complete the proof.

Since the $(3,2)$ entry is in $\Cbb[[t]]$, necessarily
$$g'(\tau^a)\equiv \psi(t)\mod {t^{c-b}}\quad;$$
so our claim is equivalent to the assertion that
$$g(\tau^a)-(\tau^a-t^a)g'(\tau^a)\equiv g((\xi t)^a)\mod t^c\quad.$$
By linearity, in order to prove this it is enough to verify the
stated congruence for $g(y)=y^\lambda$, with $\lambda\ge
\lambda_0$. That is, we have to verify that if $\lambda\ge \lambda_0$
then
$$\tau^{a\lambda}-(\tau^a-t^a)\lambda \tau^{a\lambda-a}\equiv (\xi
t)^{a\lambda} \mod t^c\quad.$$
For this, observe
$$\tau^{a\lambda}=(\xi t)^{a\lambda}(1+t^{b-a}\mu(t))^{a\lambda}
\equiv(\xi t)^{a\lambda}(1+a\lambda t^{b-a}\mu(t)) \mod
      {t^{a\lambda+2(b-a)}}$$
and similarly
$$\tau^{a\lambda-a}=(\xi t)^{a\lambda-a}(1+t^{b-a}\mu(t))^{a\lambda-a}
\equiv t^{-a}(\xi t)^{a\lambda}\mod {t^{a\lambda-a+(b-a)}}\quad,$$
$$(\tau^a-t^a)=(\xi t)^a(1+t^{b-a}\mu(t))^a-t^a\equiv at^b\mu(t)
\mod {t^{a+2(b-a)}}\quad.$$
Thus
$$(\tau^a-t^a)\lambda \tau^{a\lambda-a}\equiv (\xi t)^{a\lambda}
a\lambda t^{b-a}\mu(t)\mod {t^{a\lambda+2(b-a)}}$$
and
$$\tau^{a\lambda}-(\tau^a-t^a)\lambda \tau^{a\lambda-a}\equiv (\xi
t)^{a\lambda} \mod {t^{a\lambda+2(b-a)}}\quad.$$
Since
$$a\lambda+2(b-a)\ge a\lambda_0+2b-2a=c\quad,$$
our claim follows.
\end{proof}

\subsection{}
The first use of this observation is in the following result.

\begin{prop}\label{typeVmultsiblings}
Two truncations $f_{(C)}(y)$, $g_{(C')}(y)$ determine the same
component $\oD$ over $D$ if and only if they are siblings.
\end{prop}

\begin{proof}
Assume that $f_{(C)}(y)$, $g_{(C')}(y)$ are siblings. Then $C=C'$, and
for an $a$-th root $\xi$ of $1$ the corresponding germs
$$\alpha(t)=\begin{pmatrix}
1 & 0 & 0 \\
t^a & t^b & 0 \\
\underline{f(t^a)} & \underline{f'(t^a)t^b} & t^c 
\end{pmatrix}
\quad,\quad
\beta(t)=\begin{pmatrix}
1 & 0 & 0 \\
t^a & t^b & 0 \\
\underline{g(t^a)} & \underline{g'(t^a)t^b} & t^c 
\end{pmatrix}$$
satisfy 
$$\aligned
\alpha(\xi t) &=\begin{pmatrix}
1 & 0 & 0 \\
(\xi t)^a & (\xi t)^b & 0 \\
\underline{f((\xi t)^a)} & \underline{f'((\xi t)^a)(\xi t)^b} & 
(\xi t)^c 
\end{pmatrix} = \begin{pmatrix}
1 & 0 & 0 \\
t^a & t^b\xi^b  & 0 \\
\underline{g(t^a)} & \underline{g'(t^a)t^b}\xi^b  & t^c \xi^c 
\end{pmatrix}\\
& = \beta(t) \cdot \begin{pmatrix}
1 & 0 & 0 \\
0 & \xi^b & 0 \\
0 & 0 & \xi^c
\end{pmatrix}\quad.
\endaligned$$
By Lemma~\ref{actionlift}, the lifts of $\alpha(\xi t)$ and $\beta(t)$
belong to the same component of $\oE$. As $\alpha(\xi t)$ only differs
from $\alpha(t)$ by a reparametrization, this shows that $\oalpha$ and
$\obeta$ belong to the same component of $\oE$, as needed.

For the converse, assume
$$\alpha(t)=\begin{pmatrix}
1 & 0 & 0 \\
t^a & t^b & 0 \\
\underline{f(t^a)} & \underline{f'(t^a)t^b} & t^c 
\end{pmatrix}
\quad,\quad
\beta(t)=\begin{pmatrix}
1 & 0 & 0 \\
t^{a'} & t^{b'} & 0 \\
\underline{g(t^{a'})} & \underline{g'(t^{a'})t^{b'}} & t^{c'} 
\end{pmatrix}$$
mark the same component $\oD$. Since the limit of $\mathcal C$ along
$\beta(t)$ has 7-dimensional orbit, the action of $\PGL(3)$ on $\oD$
is transitive on a dense open set.
Therefore, we have that $\alpha(t)\cdot N$ and $\beta(t)$
are equivalent for some $N\in\PGL(3)$. By Lemma~\ref{strongeq}, there
is a $\Cbb[h][[t]]$-valued point $N(h,t)$ of $\PGL(3)$ such that
$\beta(t\nu(t)) =\alpha(t)\cdot N(1,t)$, for a unit $\nu(t)$. That is,
$$N(1,t)=\alpha(t)^{-1}\cdot\beta(\tau(t))$$
is a $\Cbb[[t]]$-valued point of $\PGL(3)$, for $\tau(t)=t\nu(t)$. By 
Lemma~\ref{technical}, this implies $a'=a$, $b'=b$, $c'=c$, and 
$\underline{g((\xi t)^a)}=\underline{f(t^a)}$,
for an $a$-th root $\xi$ of~1, showing that the truncations are siblings.
\end{proof}

\subsection{}
By Proposition~\ref{multcount}, the multiplicity of $D$ is a sum over
the distinct sibling classes of truncations producing a given
limit. Evaluating the degree of the map $\oD \to D$ requires the
determination of the inessential subgroup of the stabilizer, which
also makes crucial use of Lemma~\ref{technical}.

\begin{lemma}\label{typeVines}
Let
$$\alpha(t)=\begin{pmatrix}
1 & 0 & 0 \\
t^a & t^b & 0 \\
\underline{f(t^a)} & \underline{f'(t^a)t^b} & t^c 
\end{pmatrix}$$
be the germ determined by $C$ and the truncation
$f_{(C)}(y)=\sum_{\lambda_i<C} \gamma_{\lambda_i}y^{\lambda_i}$ as
above. Then the inessential subgroup corresponding to $\alpha(t)$
consists of the components of the stabilizer of $(\alpha(0),\lim_{t\to
  0} C\circ\alpha(t))$ containing matrices
$$\begin{pmatrix}
1 & 0 & 0 \\
0 & \eta^b & 0\\
0 & 0 & \eta^c
\end{pmatrix}$$
with $\eta$ an $h$-th root of~$1$, where $h$ is the greatest common
divisor of $a$ and all $a\lambda_i$ ($\lambda_i<C$).
\end{lemma}

\begin{proof}
For every $h$-th root $\eta$ of~$1$, any component of the stabilizer
containing a diagonal matrix of the given form is in the inessential
subgroup: indeed, such a diagonal matrix can be realized as
$\alpha(t)^{-1}\cdot \alpha(\eta t)$.

To see that, conversely, every component of the inessential subgroup
is as stated, apply Lemma~\ref{technical} with
$\beta(t)=\alpha(t)$. We find that if $\alpha(t)^{-1}\cdot
\alpha(t\nu(t))$ is a $\Cbb[[t]]$-valued point of $\PGL(3)$, then
$\nu(t)=\eta(1+t^{b-a}\mu(t))$, with $\eta$ an $a$-th root of~$1$, and
further 
$$\underline{f(t^a)}=\underline{f((\eta t)^a)}\quad,$$
that is,
$$\sum_{\lambda_i<C} \gamma_{\lambda_i} y^{\lambda_i}
=\sum_{\lambda_i<C} \eta^{a\lambda_i}\gamma_{\lambda_i}
y^{\lambda_i}\quad.$$ 
Therefore $\eta^{a\lambda_i}=1$ for all $i$ such that $\lambda_i<C$,
and $\eta$ is an $h$-th root of~$1$. 

For $\nu(t)=\eta(1+t^{b-a}\mu(t))$, the matrix $\alpha(t)^{-1}\cdot
\alpha(t\nu(t))|_{t=0}$ is lower triangular, of the form 
$$\begin{pmatrix}
1 & 0 & 0 \\
a\mu_0 & \eta^b & 0\\
 \gamma_{\lambda_0}\binom{\lambda_0}2 (a\mu_0)^2 
+ \gamma_{\frac{\lambda_0+C}2} \frac{\lambda_0+C}2 (a\mu_0) &
2\gamma_{\lambda_0}\binom{\lambda_0}2 (a\mu_0) \eta^b & \eta^c
\end{pmatrix}$$
where $\mu_0=\mu(0)$.
If a component of the stabilizer contains this matrix for some
$\nu(t)$, then it must contain all such matrices for all $\mu_0$, and
in particular that component must contain the diagonal matrix 
$$\begin{pmatrix}
1 & 0 & 0 \\
0 & \eta^b & 0\\
0 & 0 & \eta^c
\end{pmatrix}\quad;$$
the statement follows.
\end{proof}

Note that $\eta^c=(\eta^b)^2$ since $c-2b=a\lambda_0-2a$ is divisible by $h$;
this is in fact a necessary condition for the diagonal matrix above to belong
to the stabilizer. Moreover, if $\gamma_{\frac{\lambda_0+C}2}\neq0$, then
necessarily $\eta^b=1$; as the proof of the following proposition shows,
this implies $h=1$.

\begin{prop}\label{typeVmuldeg}
For the component $\oD$ determined by the truncation $f_{(C)}(y)$ as
above, 
let $A$ be the number of components of the stabilizer of the limit
$$x^{d-2S}\prod_{i=1}^S\left(zx-\frac {\lambda_0(\lambda_0-1)}2
\gamma_{\lambda_0}y^2 -\frac{\lambda_0+C}2
\gamma_{\frac{\lambda_0+C}2}yx-\gamma_C^{(i)}x^2\right)$$
(that is, by \cite{MR2002d:14083}, \S4.1, twice the number of
automorphisms $\gamma\to u\gamma+v$ preserving the $S$-tuple
$\{\gamma_C^{(1)},\dots,\gamma_C^{(S)}\}$).
Then the degree of the map $\oD \to D$ equals $\frac Ah$, where $h$ is
the number determined in Lemma~\ref{typeVines}.
\end{prop}

\begin{proof}
As the kernel line $\alpha$ must be supported on the distinguished
tangent of the limit $\mathcal X$, the stabilizer of $(\alpha,
\mathcal X)$ equals the stabilizer of $\mathcal X$, and in particular
it consists of $A$ components.

Next, observe that for $\eta_1\ne \eta_2$ two $h$-th roots of~$1$, the
two matrices
$$\begin{pmatrix}
1 & 0 & 0 \\
0 & \eta_1^b & 0\\
0 & 0 & \eta_1^c
\end{pmatrix}\quad,\quad
\begin{pmatrix}
1 & 0 & 0 \\
0 & \eta_2^b & 0\\
0 & 0 & \eta_2^c
\end{pmatrix}$$
are distinct: indeed, if $\eta^b=\eta^c=1$, then the order of $\eta$
divides every exponent of every entry of $\alpha(t)$, hence it equals~$1$
by the minimality of~$a$. Further, the components of the stabilizer
containing these two matrices must be distinct: indeed, the
description of the identity component of the stabilizer of a curve
consisting of quadritangent conics given in \cite{MR2002d:14084}, \S1,
shows that the only diagonal matrix in the component of the identity
is in fact the identity itself. 

Hence the index of the inessential subgroup equals $A/h$, and the
statement follows then from Proposition~\ref{degreetool}.
\end{proof}

\subsection{}
All that is left is the computation of the {\em weight\/} of the
marker germ
$$\alpha(t)=\begin{pmatrix}
1 & 0 & 0 \\
t^a & t^b & 0 \\
\underline{f(t^a)} & \underline{f'(t^a)t^b} & t^c 
\end{pmatrix}\quad.$$
For every formal branch $\beta$ of $\mathcal C$ at $p$, define an
integer $w_\beta$ as follows:
\begin{itemize}
\item if the branch is not tangent to the line $z=0$, then
  $w_\beta=1$;
\item if the branch is tangent to the line $z=0$, but does not
  truncate to $f_{(C)}$, then $w_\beta=$ the first exponent at which
  $\beta$ and the truncation differ;
\item if the branch truncates to $f_{(C)}$, then $w_\beta=C$.
\end{itemize}

\begin{lemma}\label{typeVmulweight}
The weight of $\alpha(t)$ equals $aW$, with $W=\sum w_\beta$.
\end{lemma}

\begin{proof}
It is immediately checked that $aw_\beta$ equals the order of vanishing
in $t$ of the composition of each formal branch with $\alpha(t)$, so
$a\sum w_\beta$ equals the order of vanishing of $F\circ\alpha(t)$,
which is the claim.
\end{proof}

We are finally ready to conclude the verification of the multiplicity
statement for components of type~V given in \S\ref{multstat}.

The upshot of the foregoing discussion is that the multiplicity equals
the sum of contributions from each sibling class of truncations
$f_{(C)}(y)$.

\begin{prop}\label{typeVlast}
Let $D$ be the component of type~V determined by the choice of $C$ and
of the truncation $f_{(C)}(y)$, and let $\ell$ be the minimum among the
positive integers $\mu$ such that $f_{(C)}(y^\mu)$ has integer
exponents. Then, with notations as above, the contribution of the
sibling class of $f_{(C)}(y)$ to the multiplicity of $D$ is $\ell W A$.
\end{prop}

\begin{proof} By Proposition~\ref{typeVmuldeg} and
Lemma~\ref{typeVmulweight}, the sibling class of $f_{(C)}(y)$
contributes $a W\frac Ah$. So all we have to prove is that $\ell=\frac
ah$, with $\ell$ as in the statement.

For this, let $\lambda_i$, $i=1,\dots,r$ be the exponents appearing in
$f_{(C)}(y)$. If $h'$ is any divisor of $a$ and all $a\lambda_i$, then
as $\frac a{h'}\lambda_i$ are integers, necessarily $\frac a{h'}$ is a
multiple of $\ell$. That is, $h'$ divides $\frac a\ell$. On the other
hand, $\frac a\ell$ is a divisor of $a$ and all $a\lambda_i$. Hence
$\frac a\ell$ equals the greatest common divisor of $a$ and all
$a\lambda_i$, which is the claim.
\end{proof}

Proposition~\ref{typeVlast} completes the proof of the multiplicity
statement for type~V components; also cf.~\cite{MR2001h:14068}, \S2,
Fact~5.

This concludes the proof of the multiplicity statement given in
\S\ref{multstat}.


\section{Examples}\label{exs}

\subsection{}
In this final section we collect several explicit examples of limits
of translates of plane curves, obtained by applying the results
presented in this paper. We will describe the limits corresponding
to the different components of the PNC for the curves we will
consider, and marker germs for these components. We will generally
pass in silence degenerate limits such as multiple lines (obtained for
example as $\lim\mathcal C\circ\alpha(t)$, for $\alpha(0)$ a rank~1
matrix with image not contained in $\mathcal C$), 
or rank-2 limits.
Limits will often be described in terms of the geometry of the curve,
and representative pictures will be superimposed on the curve to
emphasize this relation; of course such pictures should not be taken
too literally.

We will also compute the degrees of the orbit closures of the curves
we will consider, as an illustration of the formulas in
\cite{MR2001h:14068} that we obtained as main application of the
results presented here.
Several other enumerative examples can be found in
\cite{MR2001h:14068}. Concerning these enumerative computations, we
will freely use the terminology introduced in loc.~cit., and in
particular the notions of {\em predegree\/} and {\em adjusted predegree
polynomial (a.p.p.)\/} (cf.~\S1 of \cite{MR2001h:14068}).

\subsection{}
Let $d_1,d_2,m_1,m_2$ be positive integers. Consider a curve
$\mathcal C$ consisting of the union of two general curves $\mathcal
C_1$, $\mathcal C_2$ of degrees $d_1\le d_2$, in general position and
appearing with multiplicity $m_1$, $m_2$ respectively.

We distinguish three cases:
\begin{enumerate}
\item $1=d_1=d_2$;
\item $1=d_1<d_2$;
\item $1<d_1\le d_2$.
\end{enumerate}

In case (1) the curve is the union of two distinct lines with
multiplicity $m_1$, $m_2$.
$$\includegraphics[scale=.5]{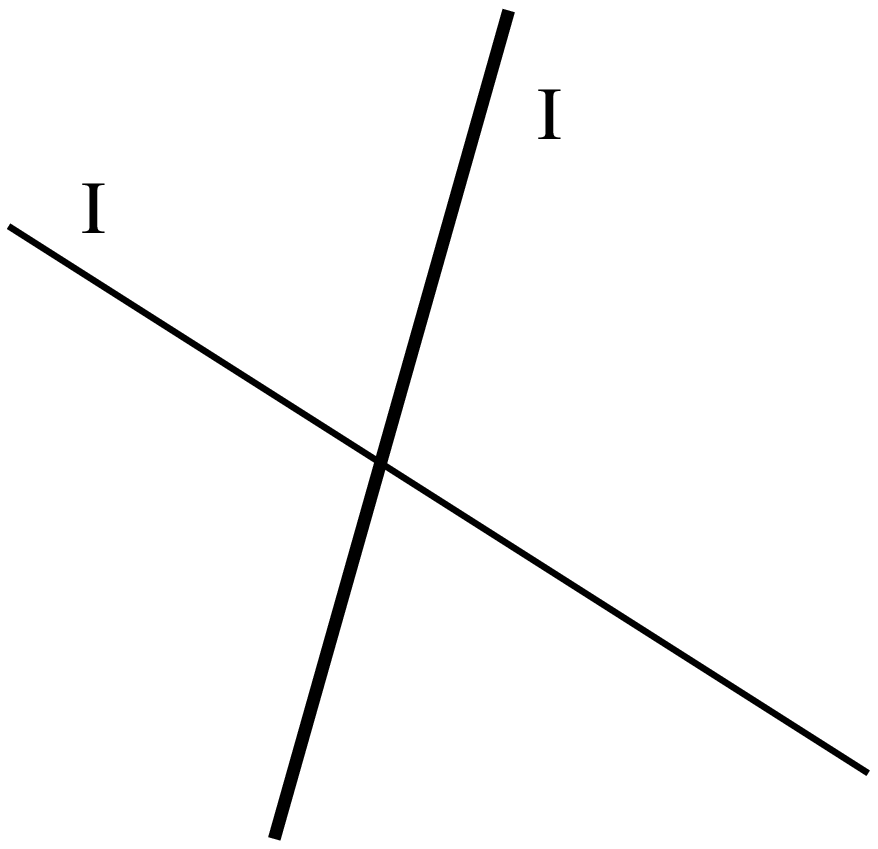}$$
According to \S\ref{statement} and \S\ref{multstat}, the PNC consists
of two components of type~I, appearing with multiplicity $m_1$, $m_2$
(also cf.~\S\ref{typeI} and \S\ref{multtypeI}); and the limits
attained by $\mathcal C$ are either translates of $\mathcal C$, or
multiple lines.

Since the line $\mathcal C_1$ meets the rest of $\mathcal C$ at one
point with multiplicity $m_2$, the contribution of $\mathcal C_1$ to
the a.p.p.~of $\mathcal C$ is the antiderivative w.r.t.~$H$ of
$$-\frac{m_1^3}2 \exp(-(m_1+m_2)H) H^2
\left(1+m_2H+\frac{m_2^2H^2}2\right)$$ 
(Proposition 3.1 in \cite{MR2001h:14068}); and similarly for $\mathcal
C_2$. The contribution of both together is therefore
$$-\frac{(m_1^3+m_2^3)H^3}6+\frac{(m_1^4+m_2^4)H^4}8
-\frac{(m_1^5+m_2^5)H^5}{20}+\frac{(m_1^3+m_2^3)^2H^6}{72}-\dots\quad,$$
so the a.p.p.~of $\mathcal C$ is (Proposition~1.1 in \cite{MR2001h:14068})
$$\exp((m_1+m_2)H)\left(1-\frac{(m_1^3+m_2^3)H^3}6+\frac{(m_1^4+m_2^4)H^4}8
-\frac{(m_1^5+m_2^5)H^5}{20}+\dots\right)$$
$$=1+(m_1+m_2)H+\frac{(m_1+m_2)^2}2H^2+\frac{m_1m_2(m_1+m_2)}2H^3
+\frac{m_1^2m_2^2}4H^4$$
$$=\left(1+m_1 H+\frac{m_1^2H^2}2\right)\left(1+m_2
H+\frac{m_2^2H^2}2\right)\quad.$$
It follows (\S1 in \cite{MR2001h:14068}) that the orbit closure of
$\mathcal C$ has dimension~4, and predegree $4!\frac{m_1^2m_2^2}4$;
that is, degree
$$\left\{\aligned
6m_1^2m_2^2\quad &\text{if $m_1\ne m_2$}\\
3m^4\quad &\text{if $m_1=m_2=m$}
\endaligned\right.$$
This of course agrees with the naive dimension count, and multiplicity
and combinatorial considerations.

In case (2) the curve is the transversal union of a line, with
multiplicity $m_1$, and a general nonsingular curve of degree $d_2$,
with multiplicity $m_2$. 
$$\includegraphics[scale=.7]{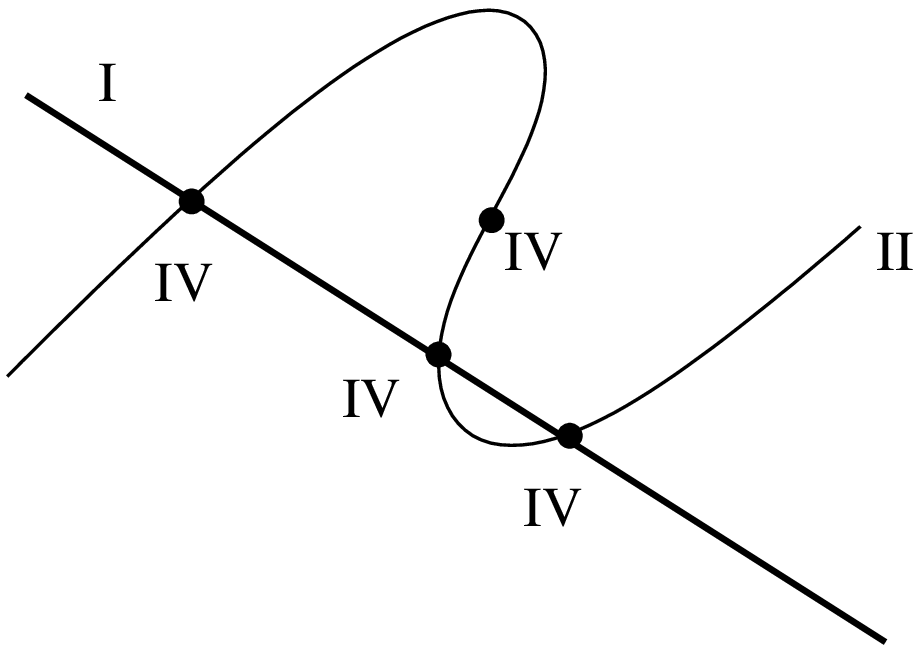}$$
According to \S\ref{statement} and \S\ref{multstat} the PNC has one
component of type~I, with multiplicity~$m_1$, one component of
type~II, with multiplicity~$2m_2$, and several `local components' of
type~IV: one for each of the $3d_2(d_2-2)$ ordinary flexes of
$\mathcal C_2$, and one for each of the $d_2$ points of intersection
of $\mathcal C_1$ and $\mathcal C_2$. 

Limits corresponding to the type~I component are obtained as
$\lim\mathcal C\circ \alpha(t)$ for $\alpha(0)$ a rank-2 matrix whose
image is the line $\mathcal C_1$. Such limits are fans determined by
the intersection of $\mathcal C_1$ with the rest of the curve:
$$\includegraphics[scale=.7]{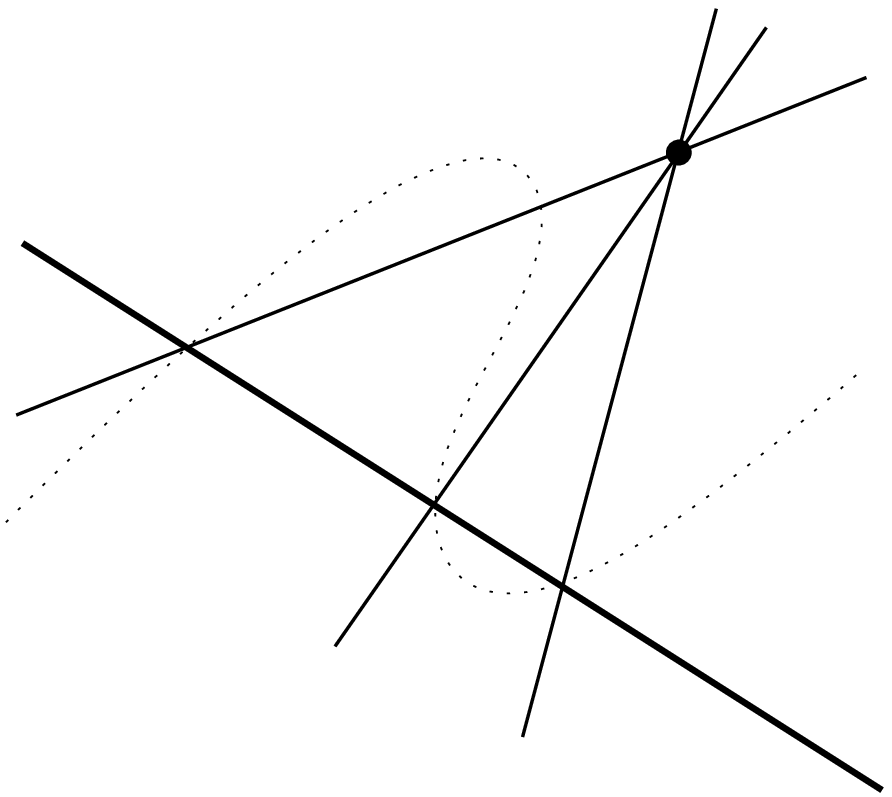}$$
The multiplicity of the non-concurrent line in the fan is $m_1$, and
the multiplicities of the star lines all equal $m_2$.

The type~II component can be marked by a 1-PS with weights $(1,2)$:
$$\alpha(t)=\begin{pmatrix}
1 & 0 & 0 \\
0 & t & 0 \\
0 & 0 & t^2
\end{pmatrix}$$
if $p=(1:0:0)$ is a general point of $\mathcal C_2$, and the highest
weight line $z=0$ is the tangent line to $\mathcal C_2$ at $p$
(cf.~\S\ref{typeII}). The corresponding limit consists of a conic,
with multiplicity $m_2$, and tangent to the kernel line, union the
kernel line with multiplicity $m_1+m_2(d_2-2)$:
$$\includegraphics[scale=.7]{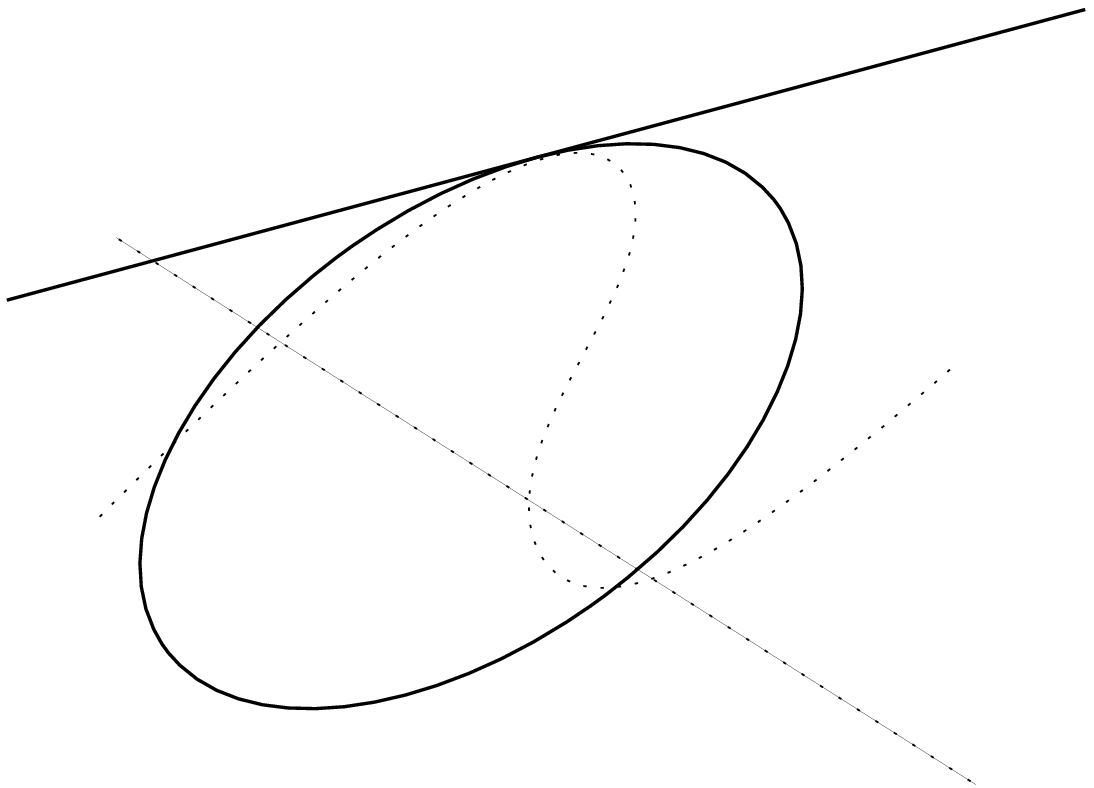}$$

At each inflection point the relevant side of the Newton polygon is
$$\includegraphics[scale=.6]{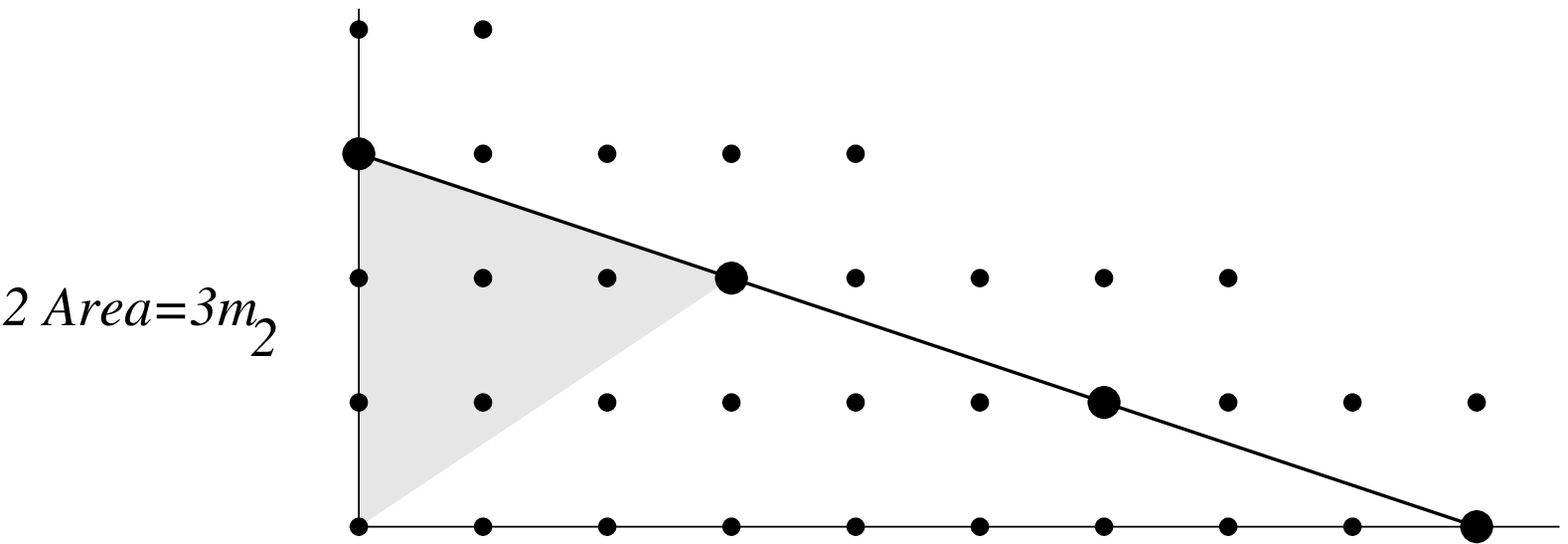}$$
joining $(0,m_2)$ and $(3m_2,0)$; according to \S\ref{multstat} such
type~IV components appear with multiplicity~$3m_2$ (cf.~\S\ref{multtypeIV}).
Marker germs for one of these components can be chosen to be 1-PS
$\alpha(t)$ with weights $(1,3)$, image of $\alpha(0)$ equal to the flex,
and highest weight line on the inflectional tangent (\S\ref{Newton}). The
limits consist of a cuspidal cubic with multiplicity $m_2$ and
cuspidal tangent on the kernel line, union the kernel line with
multiplicity $m_1+m_2(d_2-3)$:
$$\includegraphics[scale=.7]{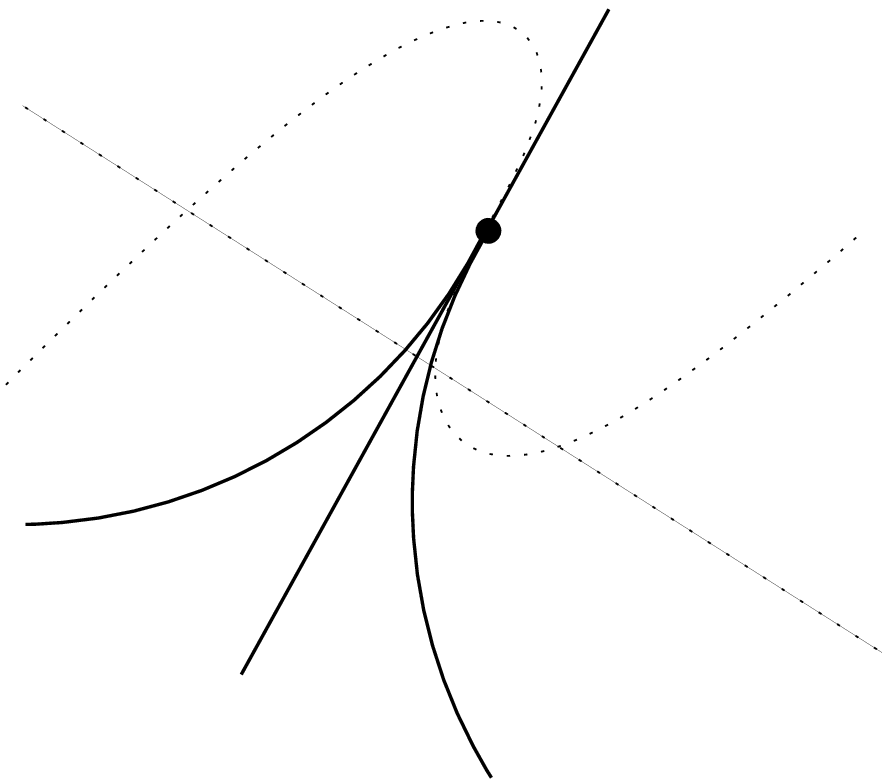}$$

The $d_2$ points of intersection are nodes, with multiple branches,
and one of whose branches is supported on a line. The Newton polygons
corresponding to the two tangent directions are
$$\includegraphics[scale=.6]{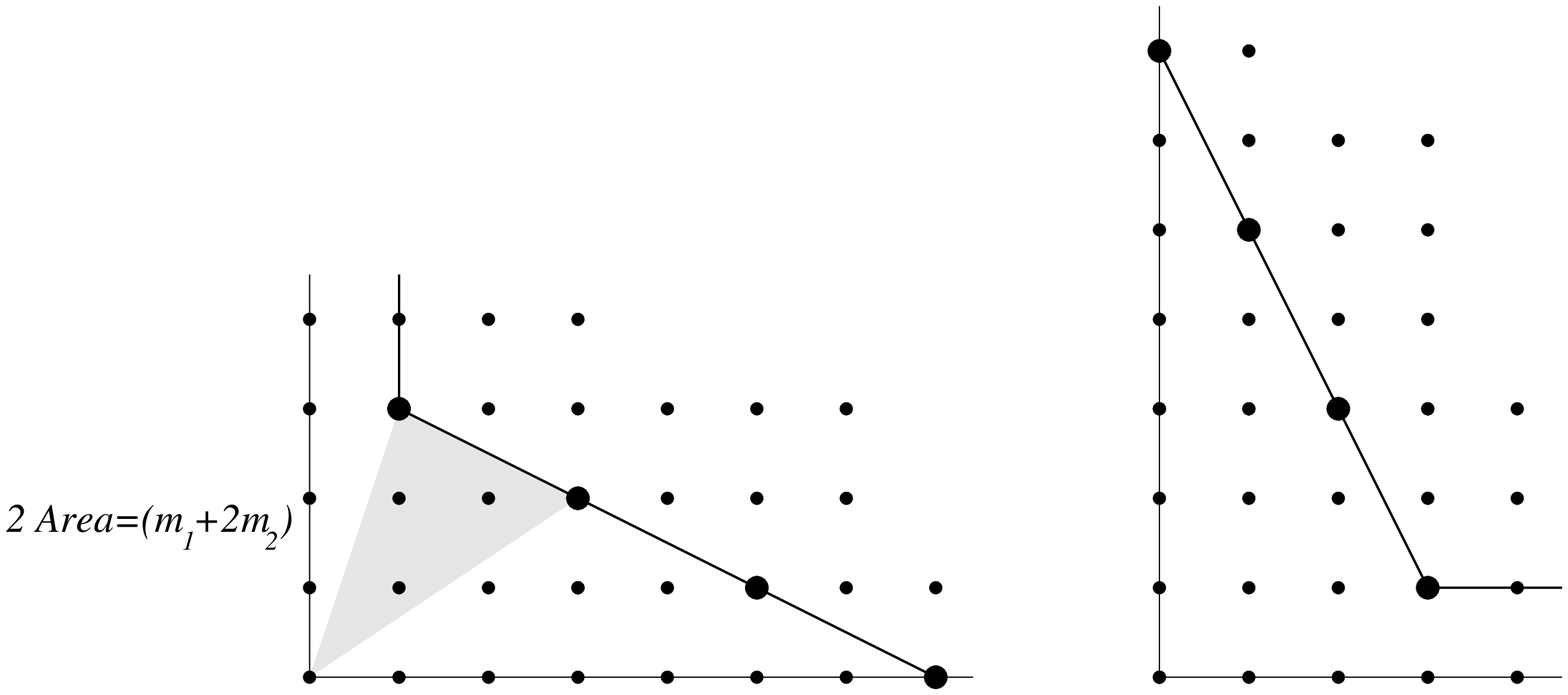}$$
and only one side (joining $(m_1,m_2)$ and $(m_1+2m_2,0)$) has slope
strictly between $-1$ and $0$. So each of these points contributes one
component of type~IV, appearing with multiplicity~$m_1+2m_2$.
This component is marked by a 1-PS $\alpha(t)$ with weights $(1,2)$
with $\im\alpha(0)$ an intersection point, and highest weight line tangent to
the non-linear branch. Limits consist of a conic with multiplicity
$m_2$ and tangent to the kernel line, the kernel line with
multiplicity $m_2(d_2-2)$, and a transversal line through the point of
intersection, with multiplicity $m_1$:
$$\includegraphics[scale=.7]{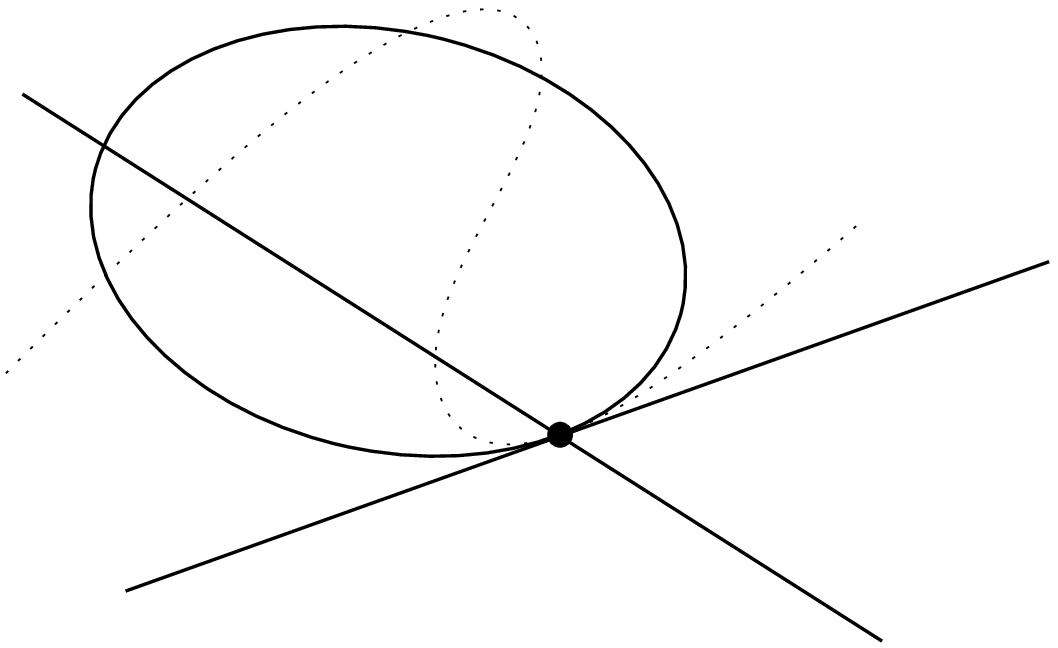}$$

No components of type~III appear, since the tangent cone at each point
of type $\mathcal C$ is supported on~$\le 2$ lines. The curve
$\mathcal C$ has no characteristics at its singularities, so the PNC has
no components of type~V in this case, cf.~again \S\ref{statement}.

Proposition~3.1 in \cite{MR2001h:14068} gives the contribution due to
the type~I component as the antiderivative of
$$-\frac{m_1^3}2 \exp(-(m_1+m_2 d_2)H) H^2
\left(1+m_2H+\frac{m_2^2H^2}2\right)^{d_2}\quad:$$ 
indeed, the line meets $\mathcal C_2$ at $d_2$ points, each with
multiplicity $m_2$. Explicitly:{\small
\begin{multline*}
-\frac{m_1^3H^3}6+\frac{m_1^4H^4}8-\frac{m_1^5H^5}{20}
+\frac{m_1^3(m_1^3+m_2^3 d_2)H^6}{72}\\
-\frac{m_1^3(m_1^4+4m_1m_2^3d_2+3m_2^4 d_2)H^7}{336}
+\frac{m_1^3(m_1^5+10m_1^2m_2^3d_2+15m_1m_2^4d_2+6m_2^5d_2)H^8}{1920}
\end{multline*}
}

As for the type~II component, its contribution is{\small
\begin{multline*}
-2 m_2^5d_2\left(\frac{H^5}{20}-\frac{(5(m_1+m_2d_2)+18m_2)H^6}{360}
+\frac{(9(m_1+m_2 d_2)+8m_2)m_2H^7}{420}\right.\\
\left.-\frac{(m_1+m_2 d_2)m_2^2H^8}{60}\right)
\end{multline*}
}
according to Proposition~3.2 in \cite{MR2001h:14068}.

`Local contributions' from type~IV components are evaluated using
Proposition~3.4 in \cite{MR2001h:14068}. The data needed in order to
apply this formula consists of the vertices of the corresponding side
of the Newton polygon, and the multiplicities of the curvilinear
components in the limit. We obtain a contribution of 
$$-\frac{m_2^6 H^6}{48}+\frac{3m_2^7 H^7}{70}-\frac{197 m_2^8
  H^8}{4480}$$
from each of the $3d_2(d_2-2)$ inflection points, and of{\small
\begin{multline*}
m_1 m_2^3 (m_1+2 m_2) \left(
-\frac{(m_1+m_2)H^6}{72}
+\frac{(20m_1^2+45m_1m_2+36m_2^2)H^7}{1680}\right.\\
\left.
-\frac{(10m_1^3+35m_1^2m_2+48m_1m_2^2+32m_2^3)H^8}{1920}
\right)
\end{multline*}
}
from each of the $d_2$ nodes. Combining all these contributions and
applying Proposition~1.1 in \cite{MR2001h:14068} yields the a.p.p.~for
$\mathcal C$. The coefficient of $H^8$ in this polynomial, multiplied
by $8!$, gives the predegree of $\mathcal C$ for $d_2>2$:
\begin{multline*}
(d_2-2)d_2 m_2^6 
\left(28(d_2^4 + 2d_2^3 + 4d_2^2 - 22d_2 -33)m_1^2\right.\\
+8(d_2^5 + 2d_2^4 + 4d_2^3 + 8d_2^2 - 411d_2 + 744)m_1 m_2\\
\left.+(d_2^6 + 2d_2^5 + 4d_2^4 + 8d_2^3 - 1356d_2^2 + 5280d_2
-5319)m_2^2\right)
\end{multline*}
This expression vanishes for $d_2=2$ because in that case the orbit
closure has dimension $<8$; the predegree can then be computed from
the coefficient of $H^7$ in the a.p.p., giving $84 m_1^2
m_2^5$. Accounting for the stabilizer, this gives $21 m_1^2 m_2^5$ for
the {\em degree\/} of the orbit closure in this case, again agreeing
with naive combinatorial considerations. 

For $d_2>2$ and $m_1=m_2=1$ the expression given above counts the
number of configurations containing 8~points in general position. The
individual sub-expressions in this formula can also be given a
concrete enumerative interpretation. For example, 
$$(d_2-2)d_2\left(d_2^4 + 2d_2^3 + 4d_2^2 - 22d_2 -33\right)
=d_2^6-30 d_2^3+11d_2^2+66d_2$$
is (for $d_2>2$) the number of embeddings of a given general plane
curve of degree~$d_2$ containing 6~general points, and satisfying
the constraint of having a given general section of $\mathcal O(1)$
contained in a prescribed line.

In case (3), $\mathcal C$ is the union of two general curves of
degrees $\ge 2$, in general position. The discussion is analogous to
that given in case (2); in this case the PNC will have no component of
type~I, but a second component of type~II, with multiplicity
$2m_1$; and there will be $3d_1(d_1-2)$ new components of type~IV
corresponding to the flexes on $\mathcal C_1$, each with multiplicity
$3m_1$. A new phenomenon concerns the components of type~IV due to the
points of intersection of $\mathcal C_1$ and $\mathcal C_2$. At such
points the relevant Newton polygons:
$$\includegraphics[scale=.6]{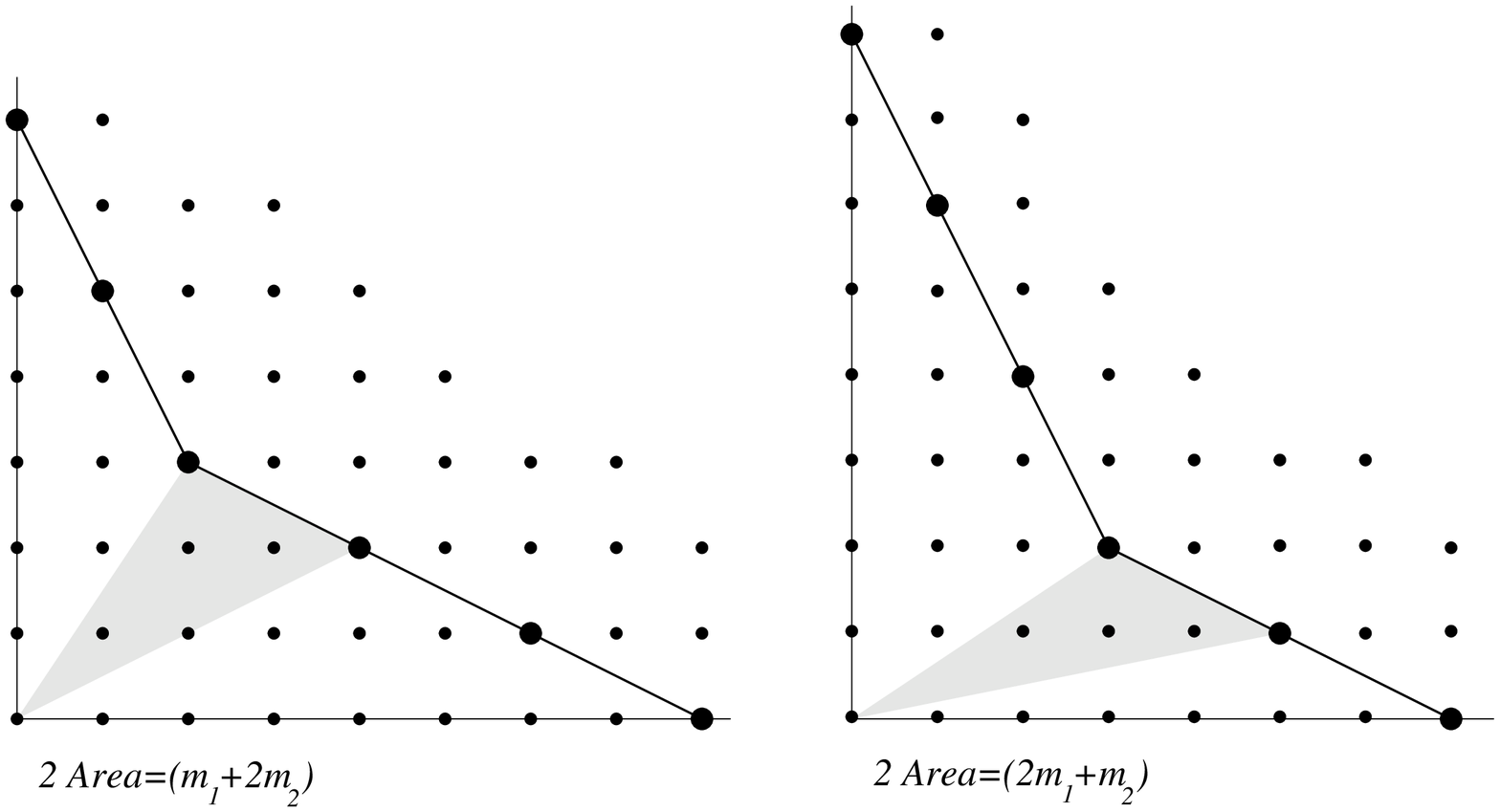}$$
have {\em two\/} sides with slope between $-1$ and $0$ (these will
join the points $(m_1,m_2)$ and $(m_1+2m_2,0)$, respectively
$(m_2,m_1)$ and $(2m_1+m_2,0)$). Thus, each of these $d_1 d_2$ points
can potentially contribute two components to the PNC. The components
are marked, as above, by 1-PS germs $\alpha(t)$ with weights $(1,2)$,
$\im\alpha(0)$ the intersection point, and highest weight line tangent
to one of the branches. The corresponding limits are schematically
represented by 
$$\includegraphics[scale=.6]{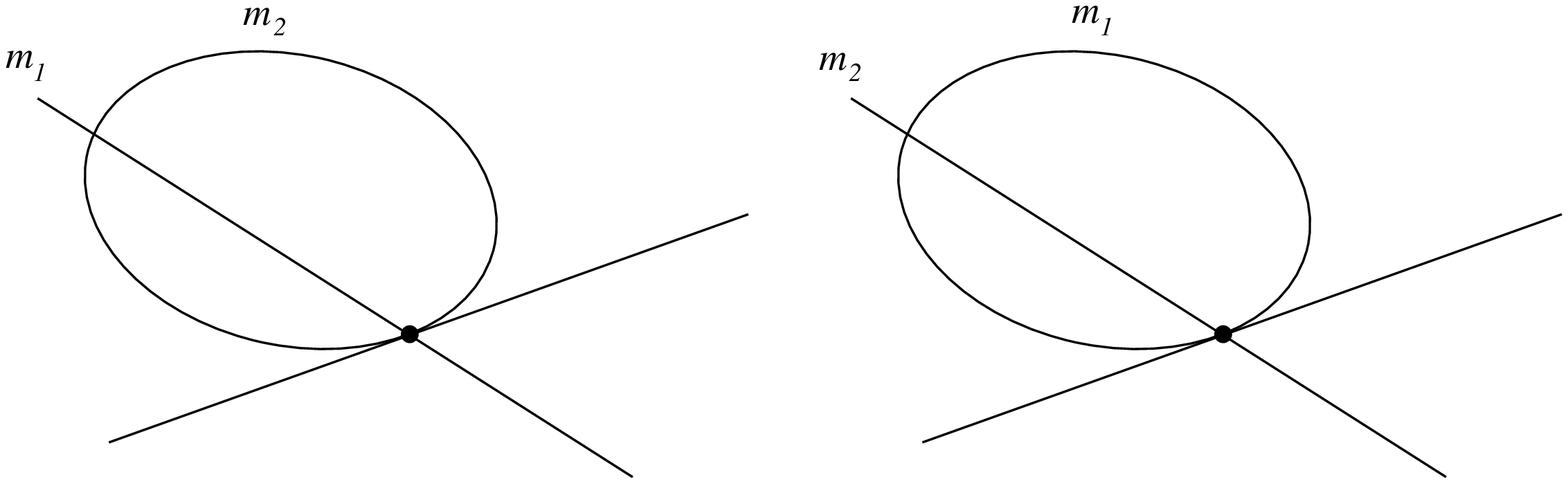}$$
indicating multiplicities. The components are distinct if and only if
$m_1\ne m_2$ (according to \S\ref{statement}). If $m_1=m_2=m$, the
single type~IV component determined by such a point has multiplicity
$6m$ in the PNC.

The a.p.p.~for $\mathcal C$ can be determined by using the results in
\cite{MR2001h:14068}, similarly to case (2), using the Newton polygon
data listed above. The predegree of $\mathcal C$ turns out to be {\small
\begin{multline*}
(m_1d_1+m_2d_2)^8-28(m_1d_1+m_2d_2)^2
(49d_1^2m_1^6+24d_1d_2m_1^5m_2+30d_1d_2m_1^4m_2^2+20d_1d_2m_1^3m_2^3 \\
+30d_1d_2m_1^2m_2^4+24d_1d_2m_1m_2^5+49d_2^2m_2^6)
+72(d_1m_1+d_2m_2)(111d_1^2m_1^7+63d_1d_2m_1^6m_2 \\
+42d_1d_2m_1^5m_2^2
+35d_1d_2m_1^4m_2^3+35d_1d_2m_1^3m_2^4+42d_1d_2m_1^2m_2^5+63d_1d_2m_1m_2^6
+111d_2^2m_2^7) \\
-(15879d_1^2m_1^8+11904d_1d_2m_1^7m_2+2688d_1d_2m_1^6m_2^2+2688d_1d_2m_1^5m_2^3
+2310d_1d_2m_1^4m_2^4 \\
+2688d_1d_2m_1^3m_2^5+2688d_1d_2m_1^2m_2^6+11904d_1d_2m_1m_2^7+15879d_2^2m_2^8)
+10638(d_1m_1^8+d_2m_2^8).
\end{multline*}
}
In the reduced case, let $d=d_1+d_2$ be the degree of the curve and
$n=d_1 d_2$ be the number of points of intersection; then this formula
evaluates the predegree of $\mathcal C$ as
$$d^8-1372 d^4+7992 d^3-15879 d^2+10638 d-24n(35 d^2-174 d+213)\quad.$$
In fact, this is the predegree of a general plane curve of degree $d$
with $n$ ordinary nodes, cf.~Example~4.1 in \cite{MR2001h:14068}.

\subsection{}
Let $\mathcal C$ be a star consisting of $d\ge 3$ distinct reduced lines
through a point:
$$\includegraphics[scale=.6]{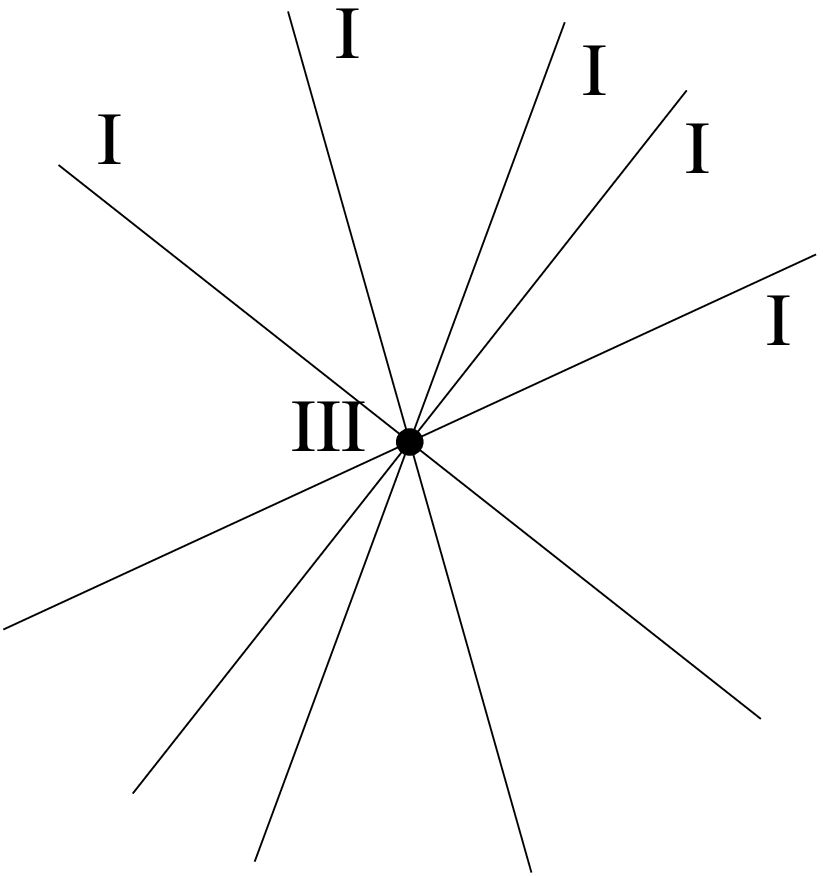}$$
and let $A$ be the number of automorphisms of the tuple of lines. Thus
$A=6$ for $d=3$, $A=4$ for general stars with $d=4$, and
$A=1$ for general stars with $d>4$.
According to \S\ref{statement} and \S\ref{multstat}, the PNC for this
curve has one reduced component of type~I for each line, and one
component of type~III, appearing with multiplicity~$dA$. If
coordinates are chosen so that one of the lines is $z=0$, the germ
$$\begin{pmatrix}
1 & 0 & 0 \\
0 & 1 & 0 \\
0 & 0 & t
\end{pmatrix}$$
marks the corresponding type~I component. In this case the limit
consists of a pair of lines, with multiplicities $1$ and $d-1$
respectively (Proposition~\ref{globallin}):
$$\includegraphics[scale=.6]{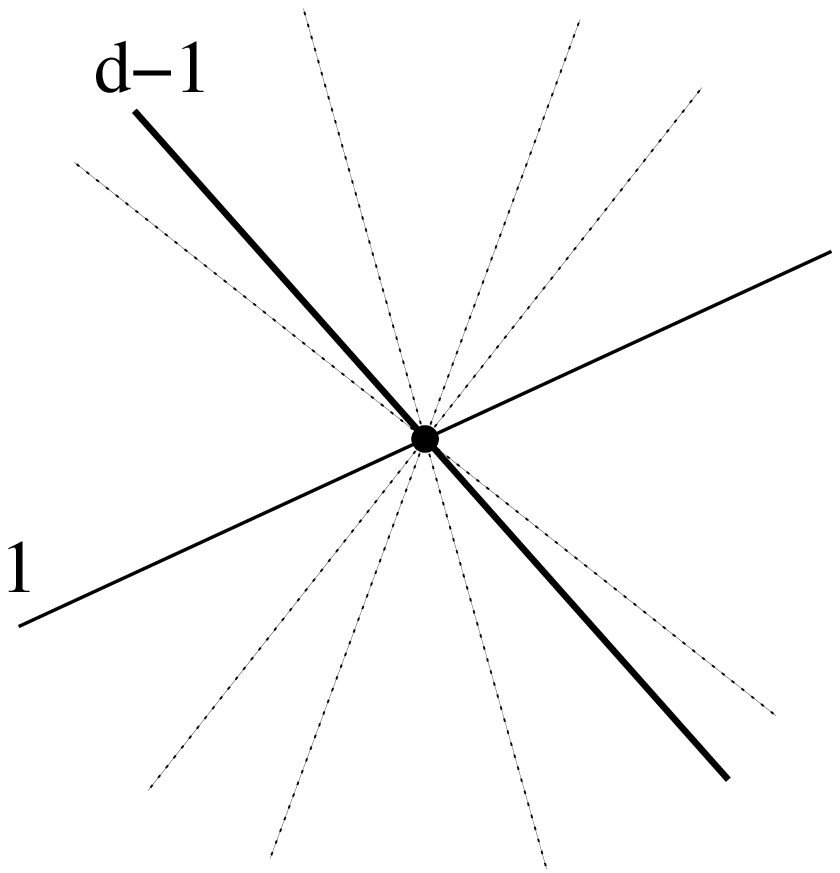}$$
The component of type~III is marked by 1-PS with image the multiple
point, and equal weights (Proposition~\ref{tgcone}):
$$\begin{pmatrix}
1 & 0 & 0 \\
0 & t & 0 \\
0 & 0 & t
\end{pmatrix}$$
if the center of the star is at $(1:0:0)$. This example is somewhat
atypical, in that {\em these limits are nothing but translates of the
  original curve.\/} 

Proposition~3.1 in \cite{MR2001h:14068} can again be invoked to
evaluate the (additive) contribution to the a.p.p.~due to the $d$
type~I components, as the antiderivative of
$$-\frac d2 \exp(-dH)H^2 \left(1+(d-1)H+\frac{(d-1)^2
  H^2}2\right)\quad.$$ 
Proposition~3.3 in loc.cit.~evaluates the contribution due to the
  type~III component:
$$-\frac{d^2(d-1)(d-2)(d^2+3d-3)}{30}\left(\frac{H^6}{24}
-\frac{d H^7}{28}+\frac{d^2 H^8}{64}\right)\quad.$$ 
Note that the factor of $A$ appearing in the multiplicity of the
type~III component is absorbed by other factors in computing this
contribution, so that the result does not depend on $A$ after all.
Also note that both contributions from type~I and~III have nonzero
coefficients for $H^6$, $H^7$, and $H^8$; however, the a.p.p.~of
$\mathcal C$ must have degree~5, because the dimension of the orbit
closure of $\mathcal C$ has dimension~5. This means that cancellations
must occur in the computation of the a.p.p., and indeed, applying
Proposition~1.1 from \cite{MR2001h:14068} gives
\begin{multline*}
1+dH+\frac{d^2 H^2}2+\frac{(d-1)d(d+1)H^3}{3!}
+\frac{(d-1)d(d^2+d-3)H^4}{4!}\\
+\frac{(d-2)(d-1)d(d^2+3d-3)H^5}{5!}\quad.
\end{multline*}
The conclusion is that the degree of the orbit closure of~$\mathcal C$
is 
$$\frac{(d-2)(d-1)d(d^2+3d-3)}{A}\quad.$$
In fact, the a.p.p.~for $\mathcal C$ is seen to equal the truncation
to $H^5$ of the power
$$\left(1+H+\frac{H^2}2\right)^d\quad,$$
a phenomenon that can be explained by `multiplicativity'
considerations such as those presented in \S4.2 of
\cite{MR2001h:14068}. Theorem~2.5~(i) in \cite{MR2002d:14084} shows
how to modify this statement in order to account for possible
multiplicities of the lines in the star.

\subsection{}
As a final example, we consider the curve $\mathcal C$ of degree $7$
from \S\ref{7ic} with equation
$$
x^3z^4-2x^2y^3z^2+xy^6-4xy^5z-y^7=0.
$$
Without difficulty, one finds that $\mathcal C$ has three singularities,
at $P=(1:0:0)$, $Q=(0:0:1)$, and $R=(1:-4:-8)$. One sees that $P$ and $Q$
are irreducible singularities, while $R$ is an ordinary node. It follows
that $\mathcal C$ is irreducible. Thus the PNC has one global component of 
type~II, with multiplicity~$2$. Clearly, there are no components of type~III.

To describe the local components of the PNC, a closer analysis of the
singularities of $\mathcal C$ is required. We begin at $P=(1:0:0)$.
It turns out that $\mathcal C$ has a very simple Puiseux expansion there:
$$\begin{cases}z=t^6+t^7,&\\y=t^4.&\end{cases}$$
(In particular, $\mathcal C$ is a rational curve.)
The singularity has two Puiseux pairs, $(2,3)$ and $(2,7)$.
In the notation of \S5 of \cite{MR2001h:14068}: $m=4$, $n=e_1=6$,
$d_1=2$, $e_2=7$, $d_2=1$, $r=2$, and the singularity absorbs $55$ flexes.

The singularity $Q=(0:0:1)$ has one Puiseux pair $(3,7)$. In the notation
of loc.~cit., $m=3$, $n=e_1=7$, $d_1=1$, $r=1$, and the singularity
absorbs $43$ flexes. The ordinary node $R$ absorbs at least $6$ flexes,
which leaves at most $1$ flex from the total number of $3d(d-2)=105$ flexes.
It turns out that $R$ is not a flecnode and that the point
$F=(7^7,2^87^3,2^{12}3)$ is a simple flex. (In the given parametrization,
$F$ corresponds to $t=-4/7$ and $R$ to $t=-1\pm i$.)

At the simple flex $F$, the relevant side of the Newton polygon joins the
points $(0,1)$ and $(3,0)$. The corresponding type~IV component appears
with multiplicity $3$. It is marked by a 1-PS $\alpha(t)$ with weights
$(1,3)$, $\im\alpha(0)=F$, and highest weight line the tangent line
to $\mathcal C$ at $F$.

At the ordinary node $R$, the two lines in the tangent cone both yield
a side of the Newton polygon joining the points $(1,1)$ and $(3,0)$.
The corresponding type~IV component appears with multiplicity $3+3=6$.
It is marked by 1-PS $\alpha(t)$ with weights $(1,2)$, $\im\alpha(0)=R$,
and highest weight line one of the two tangent lines.

At the singular point $Q=(0:0:1)$, the relevant side of the Newton polygon
joins the points $(0,3)$ and $(7,0)$. The corresponding type~IV component
appears with multiplicity $21$ and is marked by the 1-PS
$$\alpha(t)=\begin{pmatrix}
t^7 & 0 & 0 \\ 0 & t^3 & 0 \\ 0 & 0 & 1 \end{pmatrix}.$$
The limit curve has equation $x^3z^4=y^7$. The $3$ branches of
$\mathcal C$ at $Q$ do not possess a characteristic $C$, so there isn't 
a component of type~V here.

Finally, consider the singular point $P=(1:0:0)$. The relevant side
of the Newton polygon joins the points $(0,4)$, $(3,2)$, and $(6,0)$.
The corresponding type~IV component is marked by the 1-PS
$$\alpha(t)=\begin{pmatrix}
1 & 0 & 0 \\ 0 & t^2 & 0 \\ 0 & 0 & t^3 \end{pmatrix}.$$
The limit curve has equation $x^3z^4-2x^2y^3z^2+xy^6=x(y^3-xz^2)^2=0$
(a double cuspidal cubic together with its unique inflectional tangent, which
equals the kernel line). Thus the type~IV component appears with
multiplicity $12$.

In Example~\ref{7icb} we obtained the components of type~V due to $P$. The
two truncations $y^{3/2}$ and $-y^{3/2}$ are siblings (cf.~Example~\ref{sibl};
take a primitive 8th root of 1 for $\xi$). Thus we get a single 
contribution. We have $\ell=2$, $W=2\cdot\frac32+2\cdot\frac74=\frac{13}2$,
and $A=4$. Hence the multiplicity of this component equals $52$.

We conclude by computing the a.p.p.~for $\mathcal C$. The contribution
of the type~II component is
$$-\frac7{10}H^5+\frac{371}{180}H^6-\frac{71}{30}H^7+\frac{49}{30}H^8,$$
the contribution of the type~IV component due to the flex $F$ is
$$-\frac1{48}H^6+\frac3{70}H^7-\frac{197}{4480}H^8,$$
and the contribution of the type~IV component due to the node $R$ is
$$-\frac16H^6+\frac{101}{280}H^7-\frac{25}{64}H^8;$$
all three are special cases of formulas stated earlier in this section.

The contributions of the irreducible singularities $P$ and $Q$ are evaluated
using Theorem 5.1 in \cite{MR2001h:14068}. In the notation used there,
the (additive) contribution of $P$ is
$$-\left\{(24P(4,6)+2P(2,4))\cdot\left(\frac{k^2H^6}{6!}+\frac{kH^7}{7!}
+\frac{H^8}{8!}\right)\right\}_2
=-\frac{577}{30}H^6+\frac{5779}{70}H^7-\frac{6353}{35}H^8,$$
while the contribution of $Q$ equals
$$-\left\{21P(3,7)\cdot\left(\frac{k^2H^6}{6!}+\frac{kH^7}{7!}
+\frac{H^8}{8!}\right)\right\}_2
=-\frac{3059}{240}H^6+\frac{2199}{40}H^7-\frac{15775}{128}H^8$$
(note $P(1,2)=0$).

The a.p.p.~for $\mathcal C$ equals therefore the truncation to $H^8$ of
$$\exp(7H)\cdot\left(1-\frac{7}{10}H^5-\frac{5419}{180}H^6
+\frac{56939}{420}H^7-\frac{509977}{1680}H^8\right),$$
that is,
$$1+7H+\frac{49}{2}H^2+\frac{343}{6}H^3+\frac{2401}{24}H^4+\frac{16723}{120}H^5+\frac{6163}{48}H^6+\frac{119417}{1680}H^7+\frac{145139}{13440}H^8.$$
Since $P$, $Q$, $R$, and $F$ form a frame, $\mathcal C$ has trivial stabilizer.
Therefore the degree of its orbit closure equals
$$8!\cdot\frac{145139}{13440}=435417.$$



\end{document}